%% file: paper.tex
\documentclass[11pt, a4paper]{article}
\usepackage[utf8x]{inputenc}

\usepackage{amssymb}
\usepackage{amsmath}
\usepackage{amsthm}
\usepackage{comment}
\usepackage{mathrsfs}
\usepackage{mathtools}
\usepackage[british]{babel}
\newcommand{\fn}{\footnote}
\newcommand{\diag}{{\text{diag}}}

\newcommand{\FFF}{\mathbb{F}}

\newcommand{\FF}{\mathcal{F}}

\newcommand{\RR}{\mathbb{R}}

\newcommand{\Ifr}{\text{\large$\mathfrak{I}$}}

\newcommand{\KK}{\mathcal{K}}

\newcommand{\HH}{\mathcal{H}}

\newcommand{\ZZ}{\mathbb{Z}}

\newcommand{\EE}{\mathcal{E}}
\newcommand{\EEE}{{\mathbb{E}}}
\newcommand{\VVV}{{\mathbb{V}}}

\newcommand{\qv}[1]{\left<  #1  \right>}

\newcommand{\MMM}{\mathcal{M}}
\newcommand{\MMMM}{\mathfrak{M}}
\newcommand{\NNNN}{\mathfrak{N}}

\newcommand{\JJ}{g_{\pi/2}}
\newcommand{\pd}{\coloneqq} 
\newtheorem{thm}{Theorem}[section]
\newtheorem{lem}[thm]{Lemma}
\newtheorem{prop}[thm]{Proposition}
\newtheorem{cor}[thm]{Corollary} 
 \newtheorem{rem}[thm]{Remark}
\newtheorem{ex}[thm]{Example}
\newtheorem{df}[thm]{Definition}

\newenvironment{prf}{\vspace{1ex}\begin{proof}[\bf Proof]}{\end{proof}\vspace{2ex}}

\begin{document}

\title{\bf Orthogonal series for si- and related processes, Karhunen-Lo\`eve decompositions 
\\[4ex]}

\author{
Kacha Dzhaparidze\\[3ex]
%Research group Stochastics\\
Centrum Wiskunde \& Informatica (CWI)\\
Science Park 123\\
1098 XG Amsterdam\\
The Netherlands\\
kacha.dzhaparidze@gmail.com
}
\begin{comment}
\and
Harry van Zanten\\[3ex]
Department of Mathematics\\
Vrije Universiteit Amsterdam\\
De Boelelaan 1081a\\
1081 HV Amsterdam\\
The Netherlands\\
harry@cs.vu.nl\\[3ex]}
\end{comment}

\date{\today \ (this version)}

\maketitle

\begin{abstract}

This paper reproduces results from Chapter 11 of the forthcoming book \cite{dzh25}. It discusses series expansions of
processes with stationary increments (si-processes) and certain associated processes.
Making use of de Branges theory of Hilbert spaces of entire functions, it sheds new light on the existing literature and makes available some new results. In particular, it provides some new decompositions of   the Karhunen-Lo\`eve type.

\end{abstract}
%\frontmatter
%\tableofcontents
%\mainmatter

\numberwithin{equation}{section}

\section{Introduction}
This paper makes use of some material from our previous work \cite{Dzh04}--\cite{Dzh05} in collaboration with Harry van Zanten. One of the prime topics we were interested in was series expansions of {\it fractional Brownian motion} (FBM) and, in general, of {\it si-processes}. The latter acronym  stands for {\it processes with stationary increments} which are defined on  a Gaussian probability space 
$(\Omega, \FF, P)$ and are squire integrable in the metric of this space.  They are assumed to be   mean squire continuous, centered, and to emanate from the origin. Denoted by $(X_{t})_{t\ge 0}$ such process has stationary increments in the sense that  for every $u$ the increment process 
$(X_{t+u}-X_u)_{t\ge 0}$ has the same mean and covariance function as process $X$. 
The covariance function $r(s,t)=\EEE(X_s{X_t})$ is determined 
by  variance  $v(t)= \EEE |X_t|^2$ in this way 
%\begin{align}\label{eq:covsi}
$ r(s,t)= \frac12\big(v(s)+v(t)-v(t-s)\big)$
%\end{align}
for all $s,t\ge 0$, where $v(t)$ is extended to negative values of the argument by the  convention  $v(t)=v(-t)$.  Indeed, $v(t-s)=\EEE |X_{t-s}|^2= \EEE |X_t- X_s|^2 $ for $t>s$ and this in turn is equal to   $v(s)+v(t)-2 r(s, t)$. Since $X$ is mean squire continuous and starts from the origin, the variance function is a continuous function and starts from the origin, moreover $v(t)>0$ if $t>0$.
Conversely, if a centered process $X$ has a covariance function of the preceding form %\eqref{eq:covsi} 
for a continuous function $v$ starting from the origin, then $X_0=0$ a.s. and for every $u$  the process $(X_{t+u}-X_u)_{t\ge 0}$ has the same  covariance function as the process $X$. Hence $X$ is an si-process. 

\begin{comment}
A complex valued Gaussian stochastic process $X=(X_t)_{t\ge 0}$ such that $\EEE|X_t|^2<\infty$ for all $t\ge 0$ is said to have {\it stationary increments} if  the joint distribution of its real and imaginary parts $\xi=\text{Re\;}X$ and $\eta=\text{Im\;}X$ satisfies the following condition. If $x= [\xi, \eta]$ is a row-vector, then for every $u$ the process $(x_{t+u}-x_u)_{t\ge 0}$ has the same mean and matrix valued covariance function as the process $x$. By definition \eqref{eq:covrein} the latter is expressed in terms of the covariance function  $r(s,t)=\EEE(X_s\overline{X}_t)$ of process $X$ by
\[
 R(s,t)\pd \EEE (x^\top_s x_t)=\frac12\begin{bmatrix}
                               \text{Re\;} r(s,t) & \text{Im\;} r(s,t)\\
                               -\text{Im\;} r(s,t) & \text{Re\;} r(s,t)
                               \end{bmatrix}.
\]
It is determined by the matrix valued variance function $\VVV(t)\pd \EEE (x^\top_t x_t)$, since by the stationarity of increments 
$
\VVV(t-s)=\EEE (x^\top_{t-s} x_{t-s})=\EEE (x_{t}- x_{s})^\top(x_{t}- x_{s}), $
hence
\begin{align*}
 R(s,t)= \frac12\big(\VVV(s)+\VVV(t)-\VVV(t-s)\big)
\end{align*}
for all $s,t\ge 0$, with the convention $\VVV(t)=\VVV(-t)$. In conform with \eqref{eq:covsi} we now have
\begin{align}\label{eq:reimrst}
\text{Re\;} r(s,t)= \frac12\big(v(s)+v(t)-v(t-s)\big)\qquad v(t)=\EEE |X_{t}|^2\nonumber\\
\text{Im\;} r(s,t)=\frac12 \big(u(s)+u(t)- u(t-s)\big) \qquad u(t)=
\EEE (\xi_t\eta_t)
\end{align}
for all $t\ge s\ge 0$. 
 
\end{comment}
According to \cite[Theorem 2]{Kre44}
the variance function 
 permits the representation
 \begin{equation*}
 %\label{eq:Akh5.43}
  v(t)=i\gamma t+\int\Big(1-e^{i\lambda t}+\frac{i\lambda t}{1+\lambda^2}\Big)\frac{\mu(d\lambda)}{\lambda^2}
 \end{equation*}
with $v(t)=v(-t)$ for $t\in \RR$,  where $\gamma$ is a real 
constant determined uniquely by the variance function  and $\mu$ is a non-negative Borel  measure on the real line, squire integrable in the sense of 
\begin{equation}\label{eq:sqipaper}
    \int\frac{\mu(d\lambda)}{1+\lambda^2}<\infty.
\end{equation} 
The measure $\mu$ is defined essentially uniquely.
Conversely, any function $v$ defined by  Kre\u{\i}n's presentation  is a variance function of an si-process.

Due to Kre\u{\i}n's presentation, the covariance function of an si-process is represented as 
\begin{equation}\label{eq:anamova10p}
 r(s,t)=\int\big(1-e^{i\lambda s}\big)\big(1-e^{-i\lambda t}\big) \frac{\mu(d\lambda)}{\lambda^2}
\end{equation}
for $s,t\in \RR_+.$

We will be also interested in double-sided si-processes $(X_{t})_{t\in \RR}$ defined by the covariance function of the form \eqref{eq:anamova10p} but with the arguments 
$s,t$  taking their values on the whole real line $\RR$. 
A process  emanates from the origin and combines two identically distributed si-processes which run  to the left and to the right of the origin.
For convenience, the process to the left of the origin will be assumed reflected across the $x$-axes, and in this way we  deal with   two  single-sided processes  $(X_t)_{t\ge 0}$ and $(-X_{-t})_{t\ge 0}$.  Moreover, 
%like in \cite{Dzh04} 
we take the mean 
of these two processes
 at each  positive value of $t$, i.e. 
 $$X^e_t\pd \frac12\big(X_t+(-X_{-t})\big),$$ and call it the {\it even part} of the reflected si-process (this explains the upper index $e$).  The corresponding   odd part will be $$X^o_t\pd \frac12\big(X_t-(-X_{-t})\big)$$
 %The sample paths of process $X^e$ (resp. $X^o$)  are even (resp. odd).
%The sample paths of process $X^e$ (resp. $X^o$)  are even (resp. odd).
and process $X$  will be represented by its even and odd parts, in this manner
\begin{equation*}
%\label{eq:signt11}
 X_t= \text{sign}\, (t)\, X^e_{|t|}+ X^o_{|t|}.
\end{equation*}
A convenience of such a splitting  is in the form of the covariance functions  
\begin{align}\label{eq:analmp}
\EEE (X^e_s X^e_t)&=\int_\RR\frac{\sin s\lambda}{\lambda}\frac{\sin t\lambda}{\lambda}\mu(d\lambda)\nonumber\\
\EEE (X^o_s X^o_t)&= 
\int_\RR\frac{\cos s\lambda-1}{\lambda}\frac{\cos t\lambda-1}{\lambda}\mu(d\lambda)
\end{align}
for $s,t\ge 0$ and mutual independence of two processes: we have 
$
\EEE(X_{s}X_t)=\EEE(X^e_{s}X^e_{t})+\EEE(X^o_{s}X^o_{t}),
$ since 
$\EEE(X_{-s}X_t)=\EEE(X_{s}X_{-t})$. 
\\~\\
Let us overview shortly  the content of the paper. 

There are various kinds of series expansions of si-processes,  e.g. in pioneering works of  R. E. Paley and N. Wiener \cite{pal87} and P. Lévi \cite{lev92} or  popular textbooks such as M. Lo\`eve \cite[Section 34.5]{loe63} or I.I. Gikhman and A.V. Skorokhod
\cite[Section 5.1]{gikh69}. Based on their characterization of finite Fourier transforms as entire functions $f$ of exponential type at most $a>0$ which satisfy identity
\begin{equation}
\label{eq:Paleyp}
\int|f(\lambda)|^2d\lambda=\frac{\pi}{a}\sum_n\left|f\left(\frac{\theta+n\pi}{a}\right)\right|^2,
\end{equation}
Paley and Wiener   obtained a series expansion and took it as the definition for the standard Brownian motion (which they call the ``fundamental random function'', cf. \cite{lev92}, Section 13 of the complement).  They first considered the series 
$
 \displaystyle\sum_{n\in \ZZ} e^{int} Z_n
$
for $ t\in [0,2\pi]$, where $Z_n$ are i.i.d. complex valued standard Gaussian random variables. This series correspond to white noise, but the series does not converge in the usual sense. So, instead they considered its formal integral
$
\displaystyle \sum_{n\in \ZZ} \frac{e^{int}-1}{in} Z_n. %\qquad t\in [0,2\pi].
$
The latter series is shown to converge almost surely and is taken as the definition for the complex valued Brownian motion.
In \cite{Dzh04} and \cite{Dzh05} one can find extension of this result to  a complex valued fractional Brownian motion of arbitrary Hurst index, not necessarily $H=1/2$. The basic tools in the latter paper and the subsequent one \cite{Dzh05b} are Kre\u{\i}n's spectral theory and theory of vibrating strings. Combining the latter theory with de Branges' theory of reproducing Hilbert spaces of entire functions, H. Dym and H. P. McKean \cite{dym71} 
discuss  so-called ``sampling  formula'': on  p. 302 one can find  a generalization of identity \eqref{eq:Paleyp}.
Dym and McKean took this over from De Brange's theorem 22, see Section~\ref{s:samplingp}, in particular Theorem~\ref{thm:roparcop}
%that provides Parseval's identity \eqref{eq:obsthatp}. 
The theorem lay in the basis for the proofs of Theorem~\ref{thm:Pawelp77p} and its Corollary~\ref{exreop} which extend the Paley--Wiener  series expansion (PW-series) to si-processes. In Section~\ref{sub:PWFBMp} the particular case of  fractional Brownian motion receives 
special attention. Section~\ref{sub:asspp} considers moving average representations of si-processes in the form of Wiener integrals with respect to the so-called {\it fundamental martingales}, processes with martingale properties, adapted to the filtration of the given si-processes and producing equal filtration. Theorem~\ref{thm:wcmsp} provides the PW-series expansion for fundamental martingales. We briefly discuss also  
PW-series expansion for stationary processes, in particular 
the
{\it Ornstein--Uhlenbeck } and {\it autoregressive} processes.

The series of \cite[Section 34.5]{loe63} is of a completely different nature. A random process defined on a finite interval,  is decomposed into series of orthogonal functions which represent eigenfunctions of its covariance. This requires solving the Fredholm integral equation  associated with the covariance kernel of the process. There are no restrictions on processes, except  the mean squire continuity. This makes the {\it Karhunen-Lo\`{e}ve expansion}, as it is usually called, 
%\index{Karhunen-Lo\`{e}ve (KL) expansion} 
a powerful tool. But explicit solutions to the required  equations run into usual difficulties and up to now only some simple examples   are known (a list of references can be found in the recent survey \cite{naz23}).
%{\it Karhunen-Loeve Expansions and their Applications} by Liming Wang). 
We shall discuss some of these examples and investigate   more possibilities in  Section~\ref{s:KLp}.

Besides their mathematical appeal, series expansions are of a clear   practical significance. Truncated at some point, they may be used for simulation of processes in question (as regards the FBM, there is a number of good papers discussing various 
simulation techniques,  such as e.g. \cite{Mey99} or \cite{Aya02}, to mention a few).  
The level at which the series should be truncated  depends of course on the intended accuracy of the approximation.  It is therefore desirable in each concrete application to get an idea about  the rate of convergence of the series.   
We will discuss this question at hand of of a number of  examples.

\section{Preliminaries} 
\subsection{De Branges spaces}\label{s:B.2.1}
As is said in the introduction, the proofs utilize de Branges theory of Hilbert spaces of entire functions \cite{deB68}. Let us briefly review  the necessary facts from this theory. The generalization of the Paley--Wiener  space begins with replacing the  exponential functions, fundamental for  Fourier analyses, with more general entire functions which we shall call {\it de Branges functions} and denote by $E(z)$. Like the  exponential functions $e^{-iaz}$ of a complex variable $z=x+iy$,  
De Branges functions satisfy the
inequality
\begin{equation*}
%\label{eq:basicdefp}
|E(x-iy)|<|E(x+iy)|
\end{equation*}
for $y>0$ (equivalently, $|E^\sharp({z})|<|E(z)|$ for all $z$ in the upper
half-plane, where $E^\sharp({z})= \overline{E(\bar z)}$).
Each such function generates a linear space of entire functions, called the
{\it De Branges space} and denoted by   ${\mathcal H}(E). $
By definition, it consists of entire functions $f$ 
which satisfy the following conditions.\\
$-$\;The norm of an element $f$ is defined by
\begin{equation*}
%\label{eq:2ap}
\|f\|_E^2=\int\Big|\frac{f(\lambda)}{E(\lambda)}\Big|^2d\lambda<\infty.
\end{equation*}
$-$\;The restrictions of the quotients $f/E$ and $f^\sharp/E$ to the upper half-plane are of bounded type  
and have non-positive mean type in the half-plane.\\
The space ${\mathcal H}(E)$ is a Hilbert space with the inner product 
\[
 \qv{f,g}_E=\int \frac{f(\lambda)\overline g(\lambda)}{|E(\lambda)|^2}d\lambda.
\]
Moreover, it is reproducing kernel Hilbert space with the reproducing kernel
\begin{equation}
\label{eq:peofBp}
    K(w, z)=\frac{\overline{E(w)}E(z)-
\overline{E^\sharp(w)}E^\sharp(z)}{-2\pi i(z-\bar{w})}=\frac{{A(\bar w)}B(z)- {B(\bar w)}A(z)}{\pi (z-\bar{w})},
\end{equation}
which is the point evaluator in the sense that for every complex number
$w$ and every element $f$ of the space
\begin{equation*}
%\label{eq:reppp}
f(w)=\langle f, K_w\rangle_E=\int
\frac{f(\lambda)\overline{K_w(\lambda)}}{|E(\lambda)|^2}d\lambda.
\end{equation*}
In the definition of the reproducing kernel
%\eqref{eq:peofBp} 
the two entire functions  
\begin{equation*}
%\label{eq:ABoep} 
A(z)=\frac{E^\sharp(z)+E(z)}{2}\quad\quad
B(z)=\frac{E^\sharp(z)-E(z)}{2i}
\end{equation*}
are real for real $z$, and called {\it even} and {\it odd} components of the de Branges function,  since $A^\sharp(z)=A(z)$ and
$B^\sharp(z)=B(-z)$. In applications like in \cite{dym76}, an important
specific feature of the space is symmetry about
the origin, meaning that if $f(z)$  belongs to the space, then  $f(-z)$ belongs to the space as well
and its norm is equal to that of
$f(z)$. 
Any space $\mathcal{H}(E)$ generated by de
Branges function such that $E^\sharp(z)=E(-z)$ has this property, 
as
is easily verified at hand of its definition. Following \cite{dym76},
we say that a symmetric de Branges function satisfies the {\it reality condition}.
Obviously, the reality
condition is equivalent to the  condition that $A(z)$ and $B(z)$ are
even and odd entire functions, respectively, since $A^\sharp(z)=A(z)=A(-z)$
and $B^\sharp(z)=B(z)=-B(-z)$.
Theory of symmetric spaces is developed in \cite{deB62} and \cite[Section 47]{deB68}. It is shown, in particular, that if the symmetric space $\HH(E)$ contains an element with a non-zero value at the origin, then the generating de Branges function can be chosen so that $E(0)=1$. In applications that we shall discuss the si-processes will be real. Therefore their spectral measures will be symmetric about the origin (the spectral function $\mu(\lambda)=\mu(-\infty, \lambda]$ is odd  $\mu(\lambda)=-\mu(-\lambda)$)
and all de Branges spaces contained isometrically in $L^2(\mu)$ will satisfy the reality conditions (see \cite[Problem 176]{deB68} or \cite[Theorem 6.3.3]{dzh25}).

Observe that different de Branges spaces can have the same reproducing kernel. 
To see this, rewrite the numerator in the form 
$$
[A(z), B(z)] \JJ [A(w), B(w)]^*= [A(z), B(z)] U\JJ U^\top[A(w), B(w)]^*
$$ 
with the help of  matrix 
\begin{equation}\label{eq:ulabGp}
g_\theta=
\begin{bmatrix}
  \cos\theta & -\sin\theta \\
  \sin\theta & \cos\theta 
\end{bmatrix}
\end{equation} 
which defines the group of rotations about the origin of the $2$-dimensional Cartesian plane through angle  $\theta$,
see e.g. \cite[Section II.1.2]{vil68} for the unitary representation of the rotation group $SO(2)$. For  $\theta= \pi/2$,  
the latter matrix
rotates the plane through angle $\pi/2$ about the origin and   matrix $U$  in this representation of the numerator represents
any $\JJ$-unitary matrix,  i.e. one that satisfies $U \JJ U^\top =\JJ$. 
It is then clear that every space $\HH(E)$ based on the De Branges functions with the components $[A(z), B(z)] U$ have the same reproducing kernel.
For instance, the matrices \eqref{eq:ulabGp} are $\JJ$-unitary and the reproducing kernels are rotation invariant. Denote
\begin{equation}
\label{eq:ulabp}
 [A_\theta(z), B_\theta(z)]=[A(z), B(z)]\,g_\theta.
\end{equation}    
As $[A_{\pi/2}(z), B_{\pi/2}(z)]=[B(z), -A(z)].$ 
the  rotation by angle $\pi/2$ leads to interchanging the roles of $A$ and $B$,

In terms of \cite{dym76}, a de Branges space $\HH(E)$ is called {\it short} if it is invariant with respect to the backward shift operator\fn{For a fixed complex number $\alpha$
the  {\it backward 
shift operator} $R_\alpha f$ is defined at
every point $z$ by
\begin{equation*}
\label{eq:backwarshp}
 R_\alpha f(z)=
\begin{cases}
 \dfrac{f(z)-f(\alpha)}{z-\alpha} &\text{if $z\neq \alpha$}\\
f'(\alpha) &\text{if $z= \alpha$}.
\end{cases}
\end{equation*}}
for some (and hence every) complex number $\alpha$, more precisely 
\begin{equation*}
%\label{eq:ratsp}
\frac{f(z)-f(\alpha)}{z-\alpha}\in{\mathcal H}(E)
\end{equation*}
whenever $f$ belongs to ${\mathcal H}(E)$. Propositions 1 and 2 in \cite[Section 6.2]{deB68} characterize such spaces. It is stated, in particular, that
${\mathcal H}(E)$ is short if and only if
$E$ is of  exponential type and satisfies  inequality 
\begin{equation*}
%\label{eq:debfshort10p}
    \int\frac{d\lambda}{(1+\lambda^2)\,|E(\lambda)|^2}<\infty.
\end{equation*}
Such a de Branges function $E$ (and its components $A$ and $B$) does satisfy
\begin{equation*}
%\label{eq:0typeEAB-p}
    \lim_{y\nearrow \infty}\frac{1}{y}\log|E(iy)|=p
\end{equation*}
and for $0<\theta<\pi$
\begin{equation*}
%\label{eq:0typeEABR-p}
    \lim_{r\nearrow
\infty}\frac{1}{r}\log|E(re^{i\theta})|=p\sin\theta
\end{equation*}
with  some positive constant $p$.

In \cite[Sections 4.8 and 6.4]{dym76} short de Branges spaces are associated with 
{\it Kre\u{\i}n's alternative} which regards
space $L^2(\mu)$  and its 
subspace which is the  linear span of 
$(e^{i\lambda t}-1)/{\lambda}$ for $ -r\le t\le r$, 
%where $r$ is a positive number, 
closed in the metric of $L^2(\mu)$. There are only two possibilities:  the span,  denoted as
\begin{equation}
\label{eq:spanedp}
 \overline{sp}\big((e^{iz t}-1)/{z}:\,  |t|\le r\big)_\mu,
\end{equation}
$r>0$, is eider whole of $L^2(\mu)$ or a proper subspace of $L^2(\mu)$. In the latter case, it comprises  class $\mathbf{I}_{r}(\mu)$ of entire functions  which are of exponential type at most $r$ and squire integrable with respect to the spectral measure $\mu$. This class, in its turn, is equivalent to the de Branges space $\HH(E)$ of type $r>0$ which is contained isometrically in $L^2(\mu)$.

The spectral analysis of si-processes via theory of de Branges spaces is based on the following  statement which is  formulated as Problem 127 in \cite{deB68} and proved by approximation of $\HH(E)$ with finite dimensional spaces.\fn{In  \cite{dzh25}
this statement is formulated and proved as
Theorems 4.3.4 and 4.5.4.}
\begin{thm}\label{thm:thm40p}
 (i) Given an arbitrary spectral measure 
$\mu$ of property  \eqref{eq:sqipaper}, there exists a short de Branges space $\HH(E)$ of  exponential type
 which is contained isometrically in $L^2(\mu)$.
\end{thm}
\noindent
Moreover, in \cite{deB68} De Branges proves Theorem 40 which says that there exists whole {\it chain of de Branges spaces} 
\begin{equation}\label{eq:famofspp}
 \mathcal{H}(E_t),\; t\in [0, \infty),
\end{equation}
which is contained isometrically in $L^2(\mu)$. 
In general, any finite or infinite interval on the real line is allowed, but since si-processes emanate  from the origin and run unboundedly to the left or right, 
we restrict our attention to $[0, \infty)$. Also, there may be isolated points at which isometry is spoiled, but this will be excluded in the sequel. Taking this into consideration, we provide only  a weakened version of \cite[Theorem 40]{deB68}, suited to our purposes. Moreover, it will be assumed throughout that the basic de Branges function has no real zeros and is normalized  $E_t(0)=1$.

\paragraph{Theorem~\ref{thm:thm40p}} (continued.) {\it
(ii) Let a short de Branges space $\mathcal{H}(E)$ of exponential type $r>0$ be contained isometrically in $L^2(\mu)$, where $\mu$ is a 
spectral measure of a real si-process.  
Then there exists a chain of spaces 
\eqref{eq:famofspp}
and a diagonal structure function $m_t=\text{diag\,}[\alpha_t, \gamma_t]$,     
$t\in [0,\infty)$, with the following properties\\
--\; $E(z)=E_c(z)$ for some number $c$ such that  $\tau_c=r$.\\
--\; 
% $ \lim_{t \nearrow \infty}(\alpha_t+\gamma_t)=\infty.$\\
%(iv)\; 
The spaces $\mathcal{H}(E_t)$ of type $\tau_t$ are contained isometrically in $L^2(\mu)$.\\
--\; $E_t(w)$ is a continuous function of $t$ for every complex $w$, and
the integral equation      
 \begin{equation}
 \label{eq:intforABp}
     [A_{b}(w), B_{b}(w)]g_{\pi/2}-[A_{a}(w), B_{a}(w)]g_{\pi/2}=w\int_a^b [A_{t}(w), B_{t}(w)]dm_t
\end{equation}
holds for $0\le a<b<\infty$. \\
--\; The kernel
$\displaystyle
%\label{eq:iesakvirveliap}
K_{t}(w,z)=\frac{{A_{t}(\bar w)}B_{t}(z)- {B_{t}(\bar w)}A_{t}(z)}{\pi (z-\bar{w})}
$
in  space $\mathcal{H}(E_t)$ is such that 
$\displaystyle
\lim_{t \searrow 0} K_{t}(w,w)=0 
$
for every complex number $w$, and 
\begin{equation}
\label{eq:KaKbp}
K_{b}(w,z)  - K_{a}(w,z)=\frac1\pi
\int_a^b[A_t(z), B_t(z)]dm_t [A_t(w), B_t(w)]^*
\end{equation}
for $0\le a<b<\infty$ and all complex numbers $w$, $z$. }
\\~\\
\noindent 
The notions of the structure function and the type defined by this structure function trough the so-called {\it type integral} 
\begin{equation}
\label{eq:typeintegral}
 \tau_b-\tau_a=\int_a^b \sqrt{\alpha'\gamma'}d\sigma
\end{equation} 
need  further explanation (the entries $\alpha$ and $\gamma$ on the diagonal of the structure function are dominated by a certain function $\sigma$ and differentiated with respect to this dominating function). In general, the structure function  is a    non-decreasing  matrix valued function whose entries are 
continuous real valued functions of argument $t$ (of the form  as in \cite[Theorem 40]{deB68}), but the reality condition 
restricts the off-diagonal entries to constants. Therefore we can ignore the off-diagonal entries and let  the  structure function  to be diagonal. In \cite[Problem 127]{deB68} the type integral is defined as the integral of $\sqrt{\text{det}\, m'}$ with respect to the dominating function $\sigma$, which in our case reduces to the type integral mentioned above. Occasionally, the determinant may vanish but in the sequel $\tau_t$  always will be  strictly increasing, $\tau'_t>0$. 

In \cite[Section 43]{deB68}      
the structure function defines  a Hilbert space $L^2([0, \infty), m)$ of vector valued Borel measurable functions on the interval, say $h=[\varphi,\psi]$, such that $\|h\|_m^2=\int  
h\,dm\, h^*$. In our case this so-called {\it structure space} falls into  two parts: $L^2([0, \infty), \alpha)$ and $L^2([0, \infty), \gamma)$ spaces of scalar valued functions with squire norms $\|\varphi\|_\alpha^2=\int  
|\varphi|^2d\alpha$ and $\|\psi\|_\gamma^2=\int  
|\varphi|^2d\gamma$. Since 
\begin{equation}\label{eq:q}
q(t,w)\pd 
[A_t(w), B_t(w)]\,1_{(0,c)}(t)
\end{equation}
belongs to   space $L^2([0,\infty), m)$  as a function of $t$ for every  number $c>0$ and 
every complex number $w$, the integrals in \eqref{eq:intforABp} and \eqref{eq:KaKbp} are well-defined. 

The next theorem from \cite{deB68} which is essential for our purposes is Theorem 44 on the eigenfunction expansion of elements of a de Branges space which belong to the given chain of spaces \eqref{eq:famofspp} contained isometrically in $L^2(\mu)$ where $\mu$ is the spectral measure of a given si-process. It defines    generalized Fourier transforms which  are the cosine and sine transforms in the classical case of the Paley--Wiener spaces  \cite[Sections 16--18]{deB68}  and the Hankel transforms in the case of the forthcoming  Section~\ref{sub:PWFBMp}. 
%In general, an eigenfunction expansion is relative to the components $A_t(z)$ and $ B_t(z)$ of the given chain.   
We again provide only  a weakened version of \cite[Theorem 44]{deB68}, adapted to the present purposes.

\begin{thm}\label{thm:thm44p}
%Let  $m_t$ be a non-decreasing matrix-valued function of the form \eqref{eq:stfuc} whose entries are continuous real-valued functions of $t$ defined in an interval $(t_-,t_+)$. Assume  $\alpha_t>0$ for $t>t_-$ and $\alpha_t\searrow 0$ as $t\searrow t_-$. 
At  some point $c>0$, let a de Branges space $\HH(E_c)$ belong to the chain of spaces \eqref{eq:famofspp} which  is contained isometrically in $L^2(\mu)$ where $\mu$ is the symmetric spectral measure of $X$. Let its diagonal structure function with an unboundedly growing trace define a strictly increasing type through the type integral
\eqref{eq:typeintegral}.\fn{In terms of Definition~\ref{def:def}, the chain \eqref{eq:famofspp} determines the first chaos associated with an si-process.}
Then \\
%(i)\; $q_0(t,w)= [A_t(w), B_t(w)]\,1_{(t_-,c)}(t)$ belongs to   space $L^2(m)$  as a function of $t$ for every regular number $c$ and every complex number $w$. \\
--\; If for each element $h=[\varphi,\psi]$ of $L^2([0,\infty), m)$ which vanishes outside of 
$[0, c]$, a function $f$ is defined by
\begin{equation}\label{eq:eigentrp}
\pi f(w)=\int_{0}^{\infty} h_t\,dm_t\, q(t,w)^\top
\end{equation}
for all complex numbers $w$, then $f$ is an entire function, it belongs to $\mathcal{H}(E_c)$ and 
\begin{equation}\label{eq:neigentrp}
\pi \int_{-\infty}^\infty |f(\lambda)/E_c(\lambda)|^2d\lambda=\int_{0}^{\infty} h \,dm_t\, h^*.
\end{equation}
--\; If  $g$ is in $\mathcal{H}(E_c)$, then $g=f$ 
for some such choice of  $h=[\varphi,\psi]$ in $L^2(m)$.
\end{thm}
Since $
q(t,w)$ in \eqref{eq:eigentrp} is defined by \eqref{eq:q}, we have
\begin{equation}\label{eq:B.11}
 \pi f(w)=\int_0^c \varphi_tA_t(w)d\alpha_t+
 \int_0^c \psi_tB_t(w)d\gamma_t.
\end{equation}
In \eqref{eq:neigentrp} the left-hand side represents $\pi$ multiple of the squire norm of $f$ in the metric of space $\HH(E_c)$. The right-hand side is the sum of the squire norms $\|\varphi\|^2_\alpha+\|\psi\|^2_\gamma$.

In the sequel, we shall distinguish two kinds of the kernels  in the  integral representation \eqref{eq:B.11}.
\begin{df}\label{df:sing}
{\rm
A kernel  $h= [\varphi,\psi]$
in the space $L^2([0,\infty), dm)$ that does vanish outside of the interval $[0,c]$
is {\it singular} of degree $n>0$ if it is orthogonal to the first $n$ coefficients in the formal power series expansion    
of the entire function $q_0(t,z)= [A_t(z), B_t(z)]\,1_{(0,c)}(t)$. If $\langle h, c_0\rangle_m\ne 0$, then $h$ is {\it non-singular}; $c_0=[1,0]$ is the free coefficient.
} 
\end{df}

Consider again, the chain of 
 short de Branges spaces  
$\mathcal{H}(E_t),\; t\in [0,\infty),$  which  satisfies  the conditions of Theorem~\ref{thm:thm44p}, $E_t(0)=1$. On discussing single-sided si-processes we shall need to shift  each space $\mathcal{H}(E_c)$ of entire functions of type at most $r=\tau_c$ by multiplying every element $f(z)$ of the space by the exponential
$e^{-izr}$. The mapping  $f(z)\mapsto e^{-izr}f(z)$ forms a set which we shall 
denote by  $e^{-izr}\HH(E_c)$. The  eigenfunction expansion of  elements of the  latter set  is just $e^{-izr}$ multiple of 
 \eqref{eq:B.11}, of course.
 
The set $e^{-izr}\HH(E_c)$, endowed with the metric of $\HH(E_c)$ is  a reproducing Hilbert space but fails to be a de Branges space due to its asymmetry.
%satisfies the axioms (H1) and (H2) of Section~\ref{sec:1.1} but not (H3). 
If $K_c(w,z)$ is a kernel of the de Branges space $\HH(E_c)$, then the reproducing kernel in  space $e^{-izr}\HH(E_c)$ will be
\begin{equation}
\label{eq:irvelia9p}
 \mathsf{K}_c(w, z)=e^{-i(z-\bar w)r} {K}_c(w, z).
\end{equation}
As is shown in  \cite[Section 6.2.3]{Dzh05}, it is useful to represent this kernel as
 \begin{equation}\label{eq:+rvirveliap}
\pi \mathsf{K}_c(w, z)=\int_{0}^c\overline{{\mathcal E}_t(w)}{\mathcal E}_t(z) d\alpha_t  
\end{equation}
for all $c>0$, where 
$$
{\mathcal E}_c(z)=e^{-izr}
\big(
{A}_{c}(z)
+
\frac{d(\beta_c-i\,\tau_c)}
{d\alpha_c}
{B}_{c}(z)
\big)
$$
is a   de Branges function.    Here $\beta_c$ is the off-diagonal entry in the structure function $m_c$,  therefore it vanishes in the present situation. 
To see that the right-hand side of \eqref{eq:+rvirveliap} is equal to the $\pi$ multiple of the right-hand side of
\eqref{eq:irvelia9p}, check that
\begin{align*}
\int_0^{c}\!\!\overline{{\mathcal E}_{t}(w)}{\mathcal E}_{t}(z) d\alpha_t=\! \int_0^c\!\! e^{-i\tau_t(z-\bar w)} [A_t(\bar w), B_t(\bar w)] d(m_t-i\tau_t g_{\pi/2})[A_t(z), B_t(z)]^\top
\end{align*}
and that this is $\pi$ times  the integral of 
$
e^{-i\tau_t(z-\bar w)} dK_t(w, z)+ K_t(w,z) de^{-i\tau_t(z-\bar w)}$.
In the latter display  $g_{\pi/2}$ is matrix \eqref{eq:ulabGp}  that rotates the plane through angle $\theta=\pi/2$ about the origin. 
Since $E_c(0)=1$, also    
${\mathcal E}_c(0)=1$ and $\mathsf{K}_c(0, z)=\mathsf{B}_c(z)/(\pi z)$ where $\mathsf{B}_c(z)=e^{-izr}{B}_c(z)  $.
Each element $f(z)$ of the shifted space $e^{-izr}\HH(E_c)$ has  the eigenfunction expansion
 \begin{equation} \label{eq:ktfemup}
  \pi f(z)=\int_{0}^c 
  {k_t}\, {\mathcal E}_t(z) \,d\alpha_t
\end{equation}
for some choice of a kernel $k_t$ which vanishes outside of interval $[0, c]$ and is squire integrable with respect to $d\alpha$.\fn{This claim is clearly true for the specific choice 
$
f(z)= \mathsf{K}_c(0, z)
$, 
since in this case  $k_t= 1_{(0, c)}(t)$ in view of  \eqref{eq:+rvirveliap} and the normalisation ${\mathcal E}_t(0)=1$.   
It is equally simple to check the statement for any finite linear combination of the type
$\sum c_j \mathsf{K}_c(w_j, z)$ which is an element of the shifted space whose eigenfunction expansion is given by $\sum c_j {\mathcal{E}_t(w_j)}1_{(0, c)}(t)$.
The general statement regarding closed  span of such linear combinations can be  deduced by arguing as in the course of proving \cite[Corollary 8.2.6]{Dzh05}.}
\begin{comment}
by the same arguments as in the proof of Theorem~\ref{thm:thm44p}. 
The role of the set $\MMM_c$ of the latter theorem is taken over by the closed linear span in $L^2(d\alpha)$ of the linear combinations $\sum c_j {\mathcal{E}_t(w_j)}1_{(0, c)}(t)$. Note that if $k_t$ is such linear combination, then the corresponding element of  $H(\EE_c)$ is $e^{-izr}$ multiple of an element of space $\HH(E_c)$, $r=\tau_c$, since \eqref{eq:matsfktwz} implies
\[
 \pi f(z) =\int_0^c {k_t}\, {\mathcal E}_t(z) \,d\alpha_t=\sum_j \bar c_j \mathsf{K}_c(w_j, z)=e^{-izr} \sum_j c'_j {K}_c(w_j, z) 
\]
with $c'_j=\bar c_je^{irw_j}$.
\end{comment}
The norm of $f$ in the metric of the shifted space is the same as in the metric of $\HH(E_c)$ and  equals to 
\begin{equation}
\label{eq:normmindap}
\pi\| f\|^2=\int_{0}^{\infty} | k_t|^2d\alpha_t. 
\end{equation}
Let us formulate these results regarding shifted spaces as the second part of Theorem~\ref{thm:thm44p}. For the sake of brevity, the shifted space $e^{-izr}\HH(E_c)$, $r=\tau_c$, will also be denoted as $H({\mathcal E}_c)$, indicating its relationship with the de Branges function ${\mathcal E}_c$.

\paragraph{Theorem~\ref{thm:thm44p}}(continued). 
{\it
Let
$H({\mathcal E}_t), t\in [0,\infty),$ be the chain of shifted spaces which satisfy the conditions of the first part    
of the theorem. At point $c>0$, the reproducing kernel  \eqref{eq:irvelia9p} is representable in the form 
\eqref{eq:+rvirveliap}. Each element of space $H({\mathcal E}_c)$  has the eigenfunction representation 
\eqref{eq:ktfemup} for some choice of a function $k_t$ which vanishes outside of interval $[0, c]$ and is squire integrable with respect to $d\alpha$.
The norm of this element in the metric of the shifted space
 is given by \eqref{eq:normmindap}.
}
\\~\\
Analogously to Definition~\ref{df:sing},  the kernel  $k$ in the integral representation \eqref{eq:ktfemup} is called  {\it singular} of degree $n>0$ if it is orthogonal to the first $n$ coefficients in the formal power series expansion    
of the de Branges function ${\mathcal E}$ in the metric of 
$L^2([0,\infty], d\alpha)$. If $\langle k, 1\rangle_\alpha\ne 0$, then $k$ is {\it non-singular}. 
\begin{comment}
\begin{cor}
\label{cor:8.3.8p}
 In the situation of Theorem~\ref{thm:thm40} and Theorem~\ref{thm:jerar} let  \eqref{eq:famofsp}  be a chain of short de Branges spaces generated by de Branges functions $E_t(z)$ of exponential type $\tau_t$ and consisting of entire functions of exponential type at most $\tau_t$. Under the same conditions as in Corollary~\ref{cor:problem231}
let
 $H(\EE_t), t\in [0,\infty),$ be the chain of spaces 
related to  \eqref{eq:famofsp} as  above, 
 with  kernel \eqref{eq:+rvirvelia}. Then for $c>0$ each element of space $H(\EE_c)$  has the eigenfunction representation 
 \eqref{eq:ktfemu} for some choice of a function $k_t$ which vanishes outside of interval $[0, c]$ and is squire integrable with respect to $d\alpha$.
 The norm of this element in the metric of $H(\EE_c)$
 is given by \eqref{eq:normminda}. 
\end{cor}
\end{comment}

\subsection{Isometry}\label{s:isopp}
The spectral representation   \eqref{eq:anamova10p} 
gives rise to an isometric relationship between random variables $X_t$ and integrals of  exponential functions $
\widehat{1}_t(z)=\int_0^t e^{iz u} du$, the 
Fourier transforms of the indicator function $1_{(0, t)}$ of the set $(0,t)$.   
With this notation  \eqref{eq:anamova10p}
is written as the inner product $\langle\widehat{1}_s,\widehat{1}_t\rangle_\mu$ in the metric of space $L^2(\mu)$.
Let us express this isometric connection of  random variables 
$X_t$
%and  $X_t$ 
with functions of a complex variable $z$ 
 in this manner 
\begin{equation}\label{eq:sourcep}
%{\Ifr}_1\, e^{it\,\cdot}= Y_t\qquad\quad
{\Ifr}_1\,
\widehat 1_{(0,t)}= X_t.
\end{equation} 
The isometry extends in a natural way to closed linear manifolds of random variables and functions of complex variable $z$ %for which we introduce separate notations
\begin{align*}
% \NNNN_{(a, b)}=\overline{sp}\{e^{izt}: a\le t\le b\}\qquad\qquad\MMMM_{(a, b)}&=\overline{sp}\{Y_t: a\le t\le b\}\\
 \NNNN_{(a, b)}=\overline{sp}\{\int_0^te^{izu}du: a\le t\le b\}\qquad
\MMMM_{(a, b)}&=\overline{sp}\{X_t: a\le t\le b\}
\end{align*}
where the  closure of the linear span of random variables $\MMMM_{(a, b)}$ is in  mean squire
%(both processes $Y$ and $X$ are mean squire continuous by assumption)
and
the linear span $\NNNN_{(a, b)}$  is closed in $L^2(\mu)$, where $\mu$ is the spectral measure of 
 process 
 %$Y$ and 
$X$. 
Note that 
%in the case of a stationary process 
the present closed linear manifold $\NNNN_{(-r, r)}$ coincides with the linear span 
%$L_r^2(\mu)$ of  \eqref{eq:spanen}, and  in the case of an si-process  with  that of 
\eqref{eq:spanedp} mentioned earlier in connection with    Kre\u{\i}n's alternative. On discussing single-sided si-processes, it is to be taken into consideration that     
$e^{-irz}$ multiple of manifold  $\NNNN_{(-r, r)}$ 
%generated by the exponential functions 
turns into  
the shifted manifold $\NNNN_{(0, 2r)}$, symbolically $e^{-irz}\,\NNNN_{(-r, r)}= \NNNN_{(0, 2r)}$.\fn{see Note 1 to \cite[Section 3.6]{dzh25}. }

The spectral isometry can be now extended by linearity and continuity to the foregoing linear manifolds   
\begin{align}\label{eq:manisometryp}
%{\Ifr}_1:\NNNN_{(a, b)}&=\overline{sp}\{e^{izt}: a\le t\le b\}\rightarrow \MMMM_{(a, b)}=\overline{sp}\{Y_t: a\le t\le b\}\\
{\Ifr}_1:
 \NNNN_{(a, b)}&=\overline{sp}\{\int_0^te^{izu}du: a\le t\le b\}\rightarrow
\MMMM_{(a, b)}=\overline{sp}\{X_t: a\le t\le b\}
\end{align}
for all intervals $(a,b)$  within which processes are defined. 
%In these (and similar) cases, we shall say  $ {\Ifr}_1$-map of an element of  manifold $\NNNN_{(a, b)}$   yields the corresponding element of manifold $\MMMM_{(a, b)}$. \index{$ {\Ifr}_1$-map}
%\\~\\ \noindent (b)\; Let $X$ be a real valued Gaussian si-process defined on  a  probability space $(\Omega, \FF, P)$ with the covariance function  \eqref{eq:anamova10}, or  \eqref{eq:cova} for its double-sided version.
Hence, at the right-hand side of this map  
one associates  with  process $X$   
the linear span of the random  variables $\MMMM_{(0, r)}=\overline{sp}\{X_t: 0\le t\le r\}$ in single-sided case and
 $\MMMM_{(-r, r)}=\overline{sp}\{X_t: -r\le t\le r\}$ in the double-sided case.
%, where $T$ is some fixed positive number, a {\it time horizon}. 
The closure in $L^2 (\Omega, \FF, P)$ of the span will be  called  the {\it first chaos associated with the process} $X$ and denoted by $\MMMM_1(r)$.\fn{The term is borrowed from Wiener's works \cite{Wie38} and \cite{wie58}. The first chaos $\MMMM_1(r)$ is clearly related to the flow of $\sigma$-algebras $\FF_{t\ge 0}$ with which the given probability space $(\Omega, \FF, P)$ is equipped, $r=\tau_t$.} 
Clearly, all random variables in the first chaos are centered and Gaussian. Since an si-process is mean squire continuous, its first chaos is separable in the metric of  $L^2 (\Omega, \FF, P)$. The set of elements of $\MMMM_1(r)$ which are linear combinations of random variables of the form $\sum c_j X_j$ with $c_j$ and $t_j$ rational, is a countable,  
dense subset.
%For each $r>0$, the spectral isometry \eqref{eq:manisometryp} yields  first chaos $\MMMM_1(r)$ associated with an si-process $X$ by mapping the linear span $\NNNN_{r}(\mu)=\overline{sp}\{\widehat{1}_t(z):  0\le t\le r\}$ closed in the metric of $L^2(\mu)$, where $\mu$ is the spectral measure of the process. % of property \eqref{eq:sqi}. 
%the closure in $L^2(\mu)$ of the set $\{\hat{1}_t(z):  |t|\le r\}$ which is denoted by $L_r^2(\mu)$,  see \eqref{eq:spaned} (in fact, the set is denoted by  $L_r^2(\mu, \CCC)$ for complex valued $z$ and the latter notation is reserved for the restriction to the real axis, but this should not cause any confusion in the present context). 
%In Section~\ref{s:kregeneral} $L^2(\mu)$ is the space of squire integrable functions with respect to the spectral measure $\mu$  of property \eqref{eq:sqi}, so as in Theorem~\ref{thm:spofsi} which defines   spectral measures of si-processes. 
Space $\MMMM_1(r)$ endowed with the  $L^2 (\Omega, \FF, P)$ norm and $\NNNN_{r}(\mu)=\overline{sp}\{\widehat{1}_t(z):  0\le t\le r\}_\mu$  endowed with the $L^2(\mu)$ norm are by construction Hilbert spaces.
In these terms,
the spectral isometry \eqref{eq:manisometryp} is
%The spectral representation \eqref{eq:anamova10} for a real valued si-process $(X_t)_{t\ge 0}$ suggests to define at each instant $t\ge 0$ an isometry ${\Ifr}_1\, \widehat 1_{(0,t)}=X_t$, which can be extended by linearity and continuity to 
the linear map
\begin{equation}\label{eq:spisometryp}
 {\Ifr}_1: \NNNN_{r}(\mu) \rightarrow \MMMM_1(r)
\end{equation}
for each $r>0$ and
$\EEE |{\Ifr}_1  f|^2= \|f\|^2_\mu $ 
for every element $f$ of  space $\NNNN_{r}(\mu)$. 

If an si-process is double-sided, then the closed linear span $\NNNN_{r}(\mu)=\overline{sp}\{\widehat{1}_t(z):  -r\le t\le r\}$ on the  
lift hand-side of \eqref{eq:spisometryp} is the same as  \eqref{eq:spanedp} which is a proper subspace of $L^2(\mu)$,  according to Kre\u{\i}n's alternative. 
As is said above,    it comprises  class $\mathbf{I}_{r}(\mu)$ of entire functions  which are of exponential type at most $r$ and squire integrable with respect to the spectral measure $\mu$. This class, in its turn, is equivalent to the de Branges space $\HH(E)$ of type $r>0$ which is contained isometrically in $L^2(\mu)$. 
On the right-hand side of \eqref{eq:spisometryp} the first chaos   $\MMMM_1(r)=\overline{sp}\{X_t: -r\le t\le r\}$ is associated with a double-sided si-process which may be represented as the sum of two independent single-sided processes     
$X^e$ and $X^o$ with  
the  covariance functions  \eqref{eq:analmp}.  Although both components are not si-processes, there is a natural way to associate with both processes corresponding first chaoses and define the spectral isometry 
between\\ 
1)\; $\NNNN_r^{ e}(\mu)$ and  $\NNNN_r^{ o}(\mu)$ which are the closure in $L^2(\mu)$ of the linear span of  functions $\{\sin tz/z:  |t|\le r\}$ and $\{(\cos tz-1)/z:  |t|\le r\}$, respectively, and \\
2)\; $\MMMM^e_1(r)$ and $\MMMM^o_1(r)$ which are the closure
in $L^2 (\Omega, \FF, P)$ of
the linear span of the  random  variables $\{X^e_t: 0\le t\le r\}$ and $\{X^o_t: 0\le t\le r\}$, respectively.\\  
We first define  isometry
${\Ifr}_1 \widehat 1^e_t =X^e_t$ and ${\Ifr}_1 \widehat 1^o_t =X^o_t$ for every $t\ge 0$, where $\widehat 1^e_t(z)$ and $\widehat 1^o_t(z)$ are cosine and sine transforms of the indicator function $1_{(0, t)}$, respectively, and then extend this
by linearity and continuity  to
\begin{equation}\label{eq:Ifrep}
{\Ifr}_1: \NNNN_r^{ e}(\mu)\rightarrow\MMMM^e_1(r)\qquad {\Ifr}_1: \NNNN_r^{ o}(\mu)\rightarrow\MMMM^o_1(r).
\end{equation}
Since the processes $X^e$ and $X^o$ are mutually independent, the sets $\MMMM^e_1(r)$ and $\MMMM^o_1(r)$ are mutually independent. Clearly, we have the orthogonal decomposition  $\NNNN_r(\mu)=\NNNN_r^{e}(\mu)\,\text{\large$\oplus$}\,\NNNN_r^{ o}(\mu)$, orthogonality is in $L^2(\mu)$. Likewise,   decomposition 
$
 \MMMM_1(r)=\MMMM^e_1(r)\,\text{\large$\oplus$}\,\MMMM^o_1(r)
$
is orthogonal in the metric of $L^2 (\Omega, \FF, P)$.  The spectral isometry is then reformulated in this manner.

\begin{thm}\label{thm:ifrHp}
For $r>0$, let  $\MMMM_1(r)$ be the first chaos associated with a given  
double-sided si-process $X$.
There exists a de Branges space 
of entire functions of exponential type at most $r$, contained isometrically in $L^2(\mu)$ where $\mu$ is a spectral measure  of the given process. 
It can  be identified with one of the spaces in a chain of spaces \eqref{eq:famofspp}, say $E(z)=E_c(z)$ for a certain  point $c>0$ which is related to type $r$ by the
type integral  \eqref{eq:typeintegral},  $\tau_c=r$. Then \eqref{eq:spisometryp} is  equivalent to
 \begin{equation}\label{eq:isodebrcp}
  \Ifr_1 : \HH(E_c) \rightarrow \MMMM_1(\tau_c)
 \end{equation}
at point $c>0$ such that $\tau_c=r$ {\rm (in terms of Definition~\ref{def:def}, chain \eqref{eq:famofspp} {\it determines} the first chaos)}.

By equivalence, $\EEE |{\Ifr}_1 f|^2=\|f\|_E^2= \|f\|^2_\mu $ for every element $f$ of the de Branges space $\HH(E)=\HH(E_c)$.   
\end{thm}
%Here $\|f\|^2_E$ is the norm \eqref{eq:2a} of an element $f$ in the metric of de Branges space $\HH(E)$, i.e. $\|f\|^2_E=\int|f(\lambda)/E(\lambda)|^2d\lambda$. 
Since the process is real,  its spectral measure is symmetric about the origin and there exists a symmetric short de Branges space contained isometrically in $L^2(\mu)$. Therefore,  the de Branges space $\HH(E)=\HH(E_c)$ is symmetric about the origin,  $E^\sharp(z)=E(-z)$, i.e,    $A(z)=A(-z)$ and
$B(z)=-B(-z)$.
\\~\\
If an si-process is single-sided, the closed linear span $\NNNN_r(\mu)$ in \eqref{eq:spisometryp} is also single-sided $\NNNN_{r}(\mu)=\overline{sp}\{\widehat{1}_t(z):  0\le t\le r\}$. It is obtained as     
$e^{-irz/2}$ multiple of manifold  $\NNNN_{(-r/2, r/2)}$
which we  denote as $e^{-irz/2}\,\NNNN_{(-r/2, r/2)}= \NNNN_{(0, r)}$.
In this case we have to modify our arguments by introducing $e^{-izr}$ multiple of  $\HH(E_c)$ which is as before  
a de Branges space of entire functions of exponential type at most $r=\tau_c$, contained isometrically in $L^2(\mu)$ where $\mu$ is a spectral measure of the process. 
The  set obtained by multiplying elements of 
this de Branges space by $e^{-izr}$ is denoted as
 $e^{-izr}\HH(E_c)$ or, alternatively, as $H({\mathcal E}_c)$. 
It is easily seen that
\\
1)\; the  set $H({\mathcal E}_c)$ is contained 
in the single-sided span $\NNNN^\sharp_{2r}(\mu)$, $r=\tau_c$,
where 
\begin{equation}\label{eq:btinmwitnp}
\NNNN^\sharp_{r}(\mu)\pd\overline{sp}\Big(\frac{1-e^{-izt}}{z}:  0\le t\le r\Big),
\end{equation}
\\
2)\; if $\NNNN^\sharp_{2r}(\mu)$ is not the whole of $L^2(\mu)$, then it is contained in  set $H({\mathcal E}_c)$.
 
As is said above,  set $H({\mathcal E}_c)$, endowed with the metric of $\HH(E)$ is  a reproducing Hilbert space but fails to be a de Branges space due to its asymmetry.
In these terms, the statement of Theorem~\ref{thm:ifrHp} regarding the double-sided case extends to the single-sided case as follows.\fn{For more details, see \cite[Section 9.3.1]{dzh25}. 
For example, in the classical case of the exponential de Branges function 
$e^{-izr} $ 
space $\HH(E)$ 
is the Paley--Wiener space of squire integrable  entire functions of type at most $r>0$, 
with the reproducing kernel 
$\displaystyle k_r(w,z)=
\frac{\sin(z-w)}
{\pi(z-w) }$. Therefore, the kernel of the shifted space $e^{-izr}\HH(E)$ becomes $\displaystyle 
\mathsf{K}_r(w, z)=
\frac{1-e^{-i2(z-w)r}}
{\pi(z-w) }$.
}
\\~\\ 
\noindent  
{\bf Theorem~\ref{thm:ifrHp}} (continued).  {\it
For  $r>0$, let  $\MMMM_1(r)$ denote the first chaos associated with a given   
single-sided si-process $X$. Let $\HH(E_c)$ be a de Branges space 
of entire functions of exponential type at most $r/2$, contained isometrically in $L^2(\mu)$ where $\mu$ is a spectral measure  of the given process. This space, shifted by multiplying with  $e^{-izr}$, $r=\tau_c$, and denoted by $H({\mathcal E}_c)$
defines the first chaos as the following $\Ifr_1$-map
\begin{equation}\label{eq:isodesharpcp}
  \Ifr_1 :H^\sharp(\mathcal E_c) \rightarrow \MMMM_1(2r)
 \end{equation}
at any point $c$ such that $r=\tau_c$.  
}
%\\ \noindent Now the same remark as at the end of paragraph (c) is to be made. In the same situation as in the {\it amplification} we will again deal with the de Branges space which is identified with one of the spaces in a chain of de Branges spaces \eqref{eq:famofsp}. In Theorem~\ref{thm:thm40} $E(z)=E_c(z)$ for a certain regular point $c$ with respect to the structure function of the chain, which  determines the type of spaces in the chain by \eqref{eq:typeint10}. Similarly to the previous case, we will have

\subsection{Moving average w.r.t. fundamental martingales}
\label{s:fm&map}
\begin{comment}
(a)\; 
Theorem~\ref{thm:isophipsi} is divided in two parts. In the first part  
$\mathfrak{M}_1(r)$, $r>0$, is the first chaos associated with a real valued double-sided si-process $X$ with a symmetric spectral measure $\mu$. It assumes constructed a chain of de Branges spaces $\HH(E_t)$ of exponential type $\tau_t$, $0\le t<\infty$. The diagonal structure function $m_t=\text{diag} [\alpha_t,\gamma_t]$ which defines the type  has  the regular starting point $t=0$  and the trace growing  unboundedly. 
The space $L^2(\mu)$ contains isometrically the de Branges space $\HH(E_c)$ whenever $c>0$ is a regular point with respect to the structure function.  
\end{comment}
\subsubsection{Fundamental martingales}
Let the first chaos for a double-sided  si-process $X$ be {\it determined} by
a chain of de Branges spaces \eqref{eq:famofspp} in the sense of \begin{df}\label{def:def}
{\rm We say that chain \eqref{eq:famofspp} {\it determines} the first chaos for an si-process $X$, if it is contained isometrically in $L^2(\mu)$ where $\mu$ is the symmetric spectral measure of $X$ and its diagonal structure function with an unboundedly growing trace defines a strictly increasing type through the type integral}
\eqref{eq:typeintegral}.
\end{df}

Theorem~\ref{thm:thm44p} establishes one-to-one correspondence between the elements  of the first chaos $\MMMM_1(r)$, $r=\tau_c$, and the  elements of $L^2([0,\infty), m)$ that vanish outside of interval $[0,c]$. This relationship is displayed through equation \eqref{eq:eigentrp} in which $f(z)$ is an
element of $\HH(E_c)$,  $h_u=[\varphi_u, \psi_u]$ is an element of  $L^2([0,\infty),m)$ that vanishes outside of interval $0\le u\le c$ and  
$[A_u(z), B_u(z)]$ are the components of the de Branges functions $E_u(z)$, $0\le u\le c$.

In this section we are interested in two particular examples. For a fixed $t$ from  interval $[0,c]$, 
in the first case $\pi[\varphi_u, \psi_u] =[1_{[0,t]}(u),  0]$ and in the second  case $\pi[\varphi_u, \psi_u]=[0,  1_{[0,t]}(u)]$, where
$1_{[0,t]}(u)$ denotes the indicator function of  interval $[0,t]$. On the right-hand side of \eqref{eq:eigentrp} appears $\pi$ multiple of  
$ \int_0^t A_u(z)d\alpha_u=B_t(z)/z$ 
in the first case and  $ \int_0^t B_u(z)d\gamma_u= (A_t(z)-1)/z$ in the second case. For each $t\in [0,c]$, both entire functions belong to space $\HH(E_c)$ and by Theorem~\ref{thm:ifrHp} their $\Ifr_1$-map belong to the first chaos $\MMMM_1(r)$, $r=\tau_c$. 
Moreover,  the first one is an even function of $z$ and the second one is an odd function of $z$, therefore
\begin{equation}
\label{eq:alsnogp}
M^e_t\pd \Ifr_1 R_0B_t \qquad\quad M^o_t\pd -\Ifr_1 R_0A_t
\end{equation}
are in $\MMMM^e_1(r)$ and  $\MMMM^o_1(r)$, respectively, 
the even and odd parts of the first chaos.
For the sake of brevity, we make use of the  {\it backward shift operator}, i.e.
$
R_0 B_t(z)= {B_t(z)}/{z}$ and $ R_0A_t(z)=(A_t(z)-1)/{z}.
$
Note that in view of the integral representation \eqref{eq:KaKbp}, the first of definitions \eqref{eq:alsnogp} is equivalent to $M_t^e=\pi\,\Ifr_1 K_t(0,\cdot)$. 
Clearly, the definition is independent of the right endpoint $c>0$ and $(M_t^e)_{t\ge 0}$ and $(M_t^o)_{t\ge 0}$ can be viewed as squire integrable processes defined on the same 
Gaussian probability space $(\Omega, \FF, P)$ as $X$, with the second order properties
\begin{equation}
\label{eq:dubelvarp}
\EEE M^e_s\,M^e_t=\pi \alpha_{s\wedge t}\qquad \EEE M^o_s\,M^o_t=\pi \gamma_{s\wedge t}
\end{equation}
for  all non-negative $s,t$. Hence, both processes are mutually independent squire integrable martingales with respect to their own filtrations $\FF^{M^e}_{t\ge 0}$ and $\FF^{M^o}_{t\ge 0}$, the $\sigma$-algebras generated by $M^e$ and $M^o$, respectively, since
\begin{equation}\label{eq:mfundMp}
 \EEE (M^e_t-M^e_s|M^e_u: 0\le u\le s)=0 \qquad a. s.
\end{equation}
for $s\le t$ (the same is true for the odd martingale, of course). Moreover, the increments are orthogonal not only to   
of $M^e_s$, but also to all elements of the even part of the first chaos $\MMMM(\tau_s)$. Let us  define the filtrations 
\begin{align}\label{eq:filtXpe}
\FF^e_{t\ge 0}\pd \MMMM^e_1(\tau_t)\qquad \FF^o_{t\ge 0}\pd \MMMM^o_1(\tau_t)
\end{align}
for even and odd parts of a given si-process and formulate the forgoing results as
\begin{thm}\label{cor:martxp}
Let $X$ be a real valued double-sided si-process and 
$\mathfrak{M}_1(r)$, $r>0$,  its first chaos determined by a chain of de Branges spaces \eqref{eq:famofspp} in the sense of Definition~\ref{def:def}. 
Then \eqref{eq:alsnogp} 
defines two mutually independent squire integrable martingales  $(M, \FF^e)_{t\ge 0}$ and  $(M, \FF^o)_{t\ge 0}$ adapted to the filtrations 
 \eqref{eq:filtXpe}, with   the quadratic variations $\langle M^e\rangle_t=\pi \alpha_t$ and $\langle M^o\rangle_t=\pi \gamma_t$.
\end{thm}
For each fixed $t$, filtration $\FF^e_{t\ge 0}$ is larger than  $\FF^{M^e}_{t\ge 0}$ (the same is true for the odd filtrations, of course). But in fact, the filtrations coincide as is not difficult to check.\fn{see \cite[Section 10.1]{dzh25}.} 
\\~\\
\noindent
The second part of  Theorem~\ref{thm:ifrHp} regards
a single-sided si-process. The spectral isometry is displayed through \eqref{eq:isodesharpcp}.
On the right-hand side  $\MMMM_1(r)$ is the first chaos  associated with a given single-sided si-process at each instant $r>0$, and 
on the left-hand side we have $e^{-izr/2}$ multiple of a de Branges space $\HH(E_c)$ of exponential type $r/2$, $r=\tau_c$.  
The reproducing kernel in this shifted  space is given by \eqref{eq:irvelia9p}.
Like in the case of even martingale, we define process 
\begin{equation}\label{eq:eigentMtp} 
  M_t\pd \pi \Ifr_1 {\mathsf K}_{t}(0,\cdot)=\Ifr_1 R_0 \,{\mathsf B}_{t}
\end{equation}
for $t\ge 0$, where $\mathsf  B_t(z)=e^{-iz\tau_t} B_t(z)$, which is again martingale with respect to  filtration 
\begin{equation}\label{eq:filtXp}
\FF^X_{t\ge 0}\pd \MMMM_1(2\tau_t)
\end{equation}
(the first chaos on the right-hand side is associated with the given single-sided si-process $X$ so as in the second part of Theorem~\ref{thm:ifrHp}). 
The second order   moments  are 
$\EEE M_s\,\overline M_t=\pi \alpha_{s\wedge t}$
for  all non-negative $s,t$.\\~\\
{\bf Theorem~\ref{cor:martxp}} (continued). {\it If in the first part of the theorem the given si-process $X$ is single-sided, then  process   
$
(M, \FF^X)_{t\ge 0}
$
defined by \eqref{eq:eigentMtp}  possesses property
\begin{equation*}
%\label{eq:mfundMMp}
 \EEE (M_t-M_s|X_u: 0\le u\le \tau_s)=0 \qquad a.s.
\end{equation*}
for $s\le t$ and 
hence, it is a squire integrable martingale with   the quadratic variation  $\langle M\rangle_t =\pi \alpha_t$.  
}

\begin{ex}\label{ex:FBM2p}
\small{{\rm 
In  \cite{Dzh04} process $X$ is a {\it fractional Brownian motion} (FBM) of Hurst index $0<H<1$ with the
covariance function 
\begin{equation}\label{eq:covfbmp}
 r(s,t)=\frac{1}{2} \big(s^{2H} +t^{2H} -|t-s|^{2H}\big)
\end{equation}
and 
the spectral measure 
\begin{equation}\label{eq:cHstvisp}
 \mu(d\lambda)=c_H |\lambda|^{1-2H}\, d\lambda\qquad c_H=\frac{\Gamma(1+2H)\sin H\pi }{2\pi}.
\end{equation}
The arguments $s,t$ are allowed to take values  on the whole real axis, so that $X$ is a double-sided si-process. The corresponding even and odd processes have the covariance functions
\begin{align*}
 4\EEE (X_s^eX_t^e)=|s+t|^{2H}{-}|s-t|^{2H}\quad
 4\EEE (X_s^oX_t^o)=2(s^{2H}{+}t^{2H}) {- }|s+t|^{2H}{-}|s-t|^{2H}.
\end{align*}
The latter process is known in  applied literature under the name  
{\it sub-fractional Brownian motion}. 
The respective spectral representations are
\begin{align*}
\EEE (X_s^eX_t^e)=c_H\int_0^\infty\!\!\! \frac{\sin\lambda s \, \sin\lambda t}{\lambda^{1+2H}}\, d\lambda\quad
\EEE (X_s^oX_t^o)= c_H\int_0^\infty\!\!\! \frac{(\cos\lambda s-1) \, (\cos\lambda t-1)}{\lambda^{1+2H}}\, d\lambda.
\end{align*}
As is shown in \cite{Dzh04}, the even and odd fundamental martingales $M^e$ and $M^o$, adapted to the filtrations of    $(X^e_t)_{t\ge 0}$ and  $(X^o_t)_{t\ge 0}$, respectively, are given by
\begin{align}\label{eq:wclamep}
M^e_r=C\int_0^r
(r^2-t^2)^{\frac12-H} 
\,dX^e_t\qquad
M^o_r=\frac{C}{\alpha_r}\int_0^r
(r^2-t^2)^{\frac12-H} \frac{d\big(\alpha_t\,X^o_t\big) }{\alpha_t'} 
\end{align}
with constant ${1}/{C} = (1-H)\,{B(\frac12, \frac32-H)}$ ($B(x,y)$ is the beta function). For $H<1/2$, the odd martingale can be written in a less complicated form 
\begin{equation}\label{eq:Hnakl12p}
M^o_r=\frac{2r^{2H}}{B(\frac12, \frac12-H)}\int_0^r 
(r^2-t^2)^{-\frac12-H} 
\,X^o_t\,dt. 
\end{equation}
The second order moments of these two squire integrable martingales are
\begin{equation}\label{eq:kvfbmp}
\EEE|M^e_t|^2=\pi \alpha_{t}=\pi \frac{t^{2-2H}}{2-2H}\qquad \EEE|M^o_t|^2=\pi \gamma_t= \pi \frac{t^{2H}}{2H},
\end{equation}
The relations inverse to \eqref{eq:wclamep} are shown to be
\begin{comment}
It will be shown in Theorem~\ref{lem:gexpr8} that 
the relations \eqref{eq:-Xostvis12} and \eqref{eq:Xostvis12} are invertible. 
We have two equations
\begin{align}\label{eq:11.2.3ship}
 \frac{\sin rz}{z}
 &= C_1\int_0^r (r^2-t^2)^{H-\frac12} A_t(z)d\alpha_t\nonumber\\
 \frac{1-\cos rz}{z} &=-C_1\int_0^r \frac{d}{dt}\Big(\int_t^r(u^2-t^2)^{H-\frac12}du\Big) B_t(z)dt
\end{align}
which imply 
\end{comment}
\begin{align}\label{eq:11.2.3eop}
 X_r^e =C_1 \int_0^r \varphi_r(t) dM^e_t \qquad X_r^o &=C_1 \int_0^r \psi_r(t) dM^o_t
\end{align}
where  $C_1=2/B(1-H, H+1/2) $ and
\begin{align*}
\varphi_r(t)=(r^2-t^2)^{H-\frac12}\qquad
\psi_r(t)=t^{2-2H}\Big((r^2-t^2)^{H-\frac12}+ \int_t^r (u^2-t^2)^{H-\frac12}\frac{du}{u^2}\Big).
\end{align*}

Usually,  the FBM  is defined as a single-sided process, a real valued Gaussian process  
$(X_t)_{t\ge 0}$ with the covariance function \eqref{eq:covfbmp} of $s,t\ge 0$, and the power spectral measure \eqref{eq:cHstvisp}. 
It is well-known that in this case the fundamental martingale is defined as
\begin{equation}\label{eq:martkfsp}
M_{r/2}=\frac{C}{2}\int_{0}^{r} t^{\frac12-H} (r-t)^{\frac12 -H}\, dX_t
\end{equation}
and that this  relationship is invertible
\begin{align}\label{eq:sionmep}
 X_{2r}=2C_1\int_0^r \Big(r^{H-\frac12}(r-t)^{H-\frac12}-\int_t^r(u-t)^{H-\frac12}\, du^{H-\frac12}\Big) dM_t.
\end{align}
The proofs can be found e.g. in \cite{Jose02}.
}}
$\hfill$ $\qed$
\end{ex}

\subsubsection{Moving average}
%\label{subs:map}

The integral representations \eqref{eq:11.2.3eop} and \eqref{eq:sionmep} are called {\it moving average}. The general definition is  as follows. Let $(\Omega, \FF, \FFF, P)$ be  a  stochastic basis equipped with the filtration  $\FFF =\FF^X_{t\ge 0}$ which is defined by the first chaos of an si-process $X$ through \eqref{eq:filtXpe} or \eqref{eq:filtXp}.  Define also 
filtration $\FF^M_{t\ge 0}$  generated by the random variables $(M_u: 0\le u\le t)$ where 
$M$ is the squire integrable  martingale  fundamental
for the si-process $X$. The filtrations are equal, as was emphasized  above. 

\begin{df}\label{df:me}
{\rm
The integral representations of the elements of the first chaos $\MMMM_1(r)$, $r=\tau_c$, in the form of  Wiener integrals with respect to the fundamental martingales 
are called  {\it moving average}  if, and  only if, the kernels involved are {\it non-singular}. 
 }
\end{df}
The definition needs further explanation. The basic fact is that every element of
the first chaos has the required integral representation. This is a consequence of Theorem~\ref{thm:thm44p} and the spectral isometry. In the case of a double-sided si-process $X$, if $f$ is an element of the first chaos $\MMMM_1(r)$, $r=\tau_c$, with the integral representation \eqref{eq:eigentrp}, then the spectral isometry represents the corresponding element ${\Ifr}_1 f\pd X_r(h)$ of the first chaos as 
\begin{equation}\label{eq:marvoorFBM1p}
X_r(h)=\int_0^c \varphi\,dM^e+
\int_0^c \psi\,dM^o
\end{equation}
where $M^e$ and $M^o$ are the squire integrable, mutually independent martingales \eqref{eq:alsnogp} and   
$h=[\varphi,\psi]$ is the same as in \eqref{eq:eigentrp}. In the case of a single-sided si-process $X$, if $f$ is an element of the first chaos $\MMMM_1(2r)$, $r=\tau_c$, with the integral representation \eqref{eq:ktfemup}, then the spectral isometry represents the corresponding element ${\Ifr}_1 f\pd X_r(k)$ of the first chaos as 
\begin{equation}\label{eq:xrkMp}
X_r(k)=\int_0^ck\,dM
\end{equation}
where $M$ is the squire integrable martingale \eqref{eq:eigentMtp} and 
 $k$ is the same as in \eqref{eq:ktfemup}. Definition~\ref{df:me} requires {\it non-singularity} of the kernels in the sense of Section~\ref{s:B.2.1}:  for $h$ see Definition~\ref{def:def} and for $k$ the last parapraph of that section. Let us reformulate this as a separate theorem.

\begin{thm}\label{thm:mavoorxypp}
In the situation of Theorem~\ref{thm:thm44p} the integral representations \eqref{eq:eigentrp} and \eqref{eq:ktfemup} for the elements
of the first chaoses associated with si-processes {\rm (double- and single-sided, respectively)} are moving average if, and only if, $\int_0^c \varphi d\alpha\ne 0$ and $\int_0^c k d\alpha\ne 0$, respectively.

In particular, the integral representations
\begin{equation}\label{eq:mavoorxpp}
X^e_{r}=\int_0^c \varphi dM^e\qquad 
X^o_{r}=\int_0^c \psi dM^o
\end{equation}
are moving average, where the kernels  come from the eigenfunction expansions
\begin{equation*}
%\label{eq:fidapsipp}
\frac{\sin zr}{z}=\int_0^c \varphi_t A_t(z)d\alpha_t\qquad
\frac{\cos zr-1}{z}=\int_0^c \psi_t B_t(z)d\gamma_t 
\end{equation*}
in space  $\HH(E_c)$ 
of exponential type $r=\tau_c$.

Likewise, the integral representation
\begin{equation}\label{eq:mavoodxpp}
X_{2r}=\int_0^c k dM
\end{equation}
is moving average, where the kernel comes from the eigenfunction expansion
\begin{equation*}
\label{eq:mavoory-pp}
\frac{e^{-2iz r}-1}{iz}=\int_0^c k_t \EE_t(z) d\alpha_t
\end{equation*}
in 
the shifted 
space $H(\mathcal E_c)$, $r=\tau_c$.

\end{thm}
Clearly, not all integral representations  for the elements of the first chaoses with respect to the fundamental martingales are moving average. Methods for determining such representations for all elements of the first chaoses can be found in \cite[Section 10.3]{Dzh05}.   

\section{Paley-Wiener series}\label{s:PWsp}
It is shown in this section how to extend the Paley--Wiener series expansion  (PW-series)  of Brownian motion to a wide class of si-processes.
For example, if the motions are fractional with Hurst index $0<H<1$, then  the PW-series
\begin{equation}
\label{eq:pwsep}
 \sum_{n\in \ZZ} \frac{e^{2in\pi t}-1}{2in\pi} Z_n\qquad t\in [0,1]
\end{equation}
which we had for $H=\frac12$ retains  this  form for $H\neq 1/2$ too, however 
the real zeros $n\pi$ of $\sin z$ are substituted by the real zeros of the Bessel function $J_{1-H}(z)$ of the first kind  
(since $\sin z=\sqrt{{\pi z}/{2}}J_{1/2}(z)$, the Bessel function 
$J_{1/2}(z)$ and the sine function have common zeros). This result and its extension to si-processes are obtained in \cite{Dzh05} (see
also the Ph.D. theses \cite{zar07}) by methods developed along the lines of Kre\u{\i}n's theory of strings so as in  \cite{dym76}.
The  present approach is based on the De Branges theory of orthogonal sets of \cite[Section 22]{deB68}. The following section reproduces some of these basic results.

\subsection{Sampling formula}\label{s:samplingp}
\noindent
With  each de Branges function $E(z)$ one associates  the so-called phase function $\varphi$ such that
$\pi K(\lambda, \lambda)=|E(\lambda)|^2 \varphi'(\lambda)>0$ at each real point $\lambda$, where $K(w,z)$ is 
the reproducing kernel \eqref{eq:peofBp} in the de Branges space $\HH(E)$. It is shown in \cite[Section 22]{deB68} that if $(\lambda_n)$ are roots of  equation $\varphi(\lambda_n)=\theta+n\pi$ for a fixed real number $\theta$, then the entire functions 
$K(\lambda_n, z)/\bar E(\lambda_n)$  
comprise orthogonal set in $\HH(E)$ so that $\pi\|K(\lambda_n, \cdot)/\bar E(\lambda_n)\|^2=\varphi'(\lambda_n). $\fn{In \cite[Section 22]{deB68} this is formulated as Problem 49; the proof is given in  \cite[Lemma 4.2.1, item iv]{dzh25}. }
Moreover, the set is complete whenever  the number $\theta$ is chosen so that rotation  \eqref{eq:ulabp} of the components of the de Branges function yields  $B_\theta$ that is not in the space $\HH(E)$.
One can always make this choice and, moreover, one can choose $\theta=0$ without loss of generality. With the latter choice  the $\lambda_n$ become  zeros of the odd component $B(z)$ (as was observed above, rotation by angle $\pi/2$ interchanges the roles of even and odd components and makes possible  to work with zeros of  component  $A(z)$).
Therefore, it will be assumed throughout we have in our disposal a complete basis  in space $\HH(E)$, which is defined by 
the set of real zeros of  $B(z)$. 
It is then possible to expand every element $f$  in terms of this basis. 
\begin{thm}\label{thm:roparcop} {\rm (\cite[Theorem 22]{deB68})}.
Let ${\mathcal H}(E)$ be a given de Branges space and let $\varphi$ be
the phase function associated with $E$. Let $(\lambda_n)$ be a
sequence of real numbers such that $\varphi(\lambda_n)= \theta$
modulo $\pi$, with some fixed real number $\theta$. Then the
functions $ {K_{\lambda_n}}/{\bar E(\lambda_n)} $ comprise an
orthogonal set in ${\mathcal H}(E)$. The only elements of
 ${\mathcal H}(E)$ which are orthogonal to the functions $
{K_{\lambda_n}}/{\bar E(\lambda_n)} $ for every $n$ are constant
multiples of $B_\theta$. If this function does not belong to the space
${\mathcal H}(E)$, then
\begin{equation}\label{eq:obsthatp}
\int\left|\frac{f(\lambda)}{E(\lambda)}\right|^2d\lambda
  =
\sum_n\left|\frac{f(\lambda_n)}{E(\lambda_n)}\right|^2\frac{\pi}{\varphi'(\lambda_n)}
\end{equation}
for every element $f$ of the space.

\end{thm}
By Theorem \ref{thm:roparcop} series
\begin{equation}\label{eq:9.47isp}
 \sum_n f(\lambda_n)\frac{K(\lambda_n, z)}{K(\lambda_n, \lambda_n)}
\end{equation}
does converge in mean square to $f(z)$, in the 
metric of $\HH(E)$.  Parseval's identity estimates the norm 
\begin{equation}\label{eq:samplingp}
 \|f\|^2_E=\sum_n\frac{|f(\lambda_n)|^2}{K(\lambda_n, \lambda_n)}.
\end{equation}
The coefficients $f(\lambda_n)/\sqrt{K(\lambda_n,\lambda_n)}$ are calculated by making use of   the reproducing property of the kernel. This is the {\it sampling  formula} in terms of
Dym and McKean \cite[p. 302]{dym71}. 
\begin{comment}
In Section~\ref{s:f} which is 
devoted to homogeneous  de Branges spaces, 
the sampling formula \eqref{eq:normunterB} is discussed 
in the course of proving Theorem~\ref{thm:homoki};  cf. also Note 4 to Section~\ref{s:theorem22} where Theorem~\ref{thm:roparco} is applied to this homogeneous case. In Section~\ref{subs:sl2} finite dimensional spaces $\HH(E)$ are generated by polynomial de Branges functions and the sampling formula is the same as the interpolation formula \eqref{eq:theexpan}: any element $f(z)$ of a $d$-dimensional space is represented as 
\begin{equation*}
 f(z)= \sum_{n=1}^d  f(\lambda_n)\frac{K(\lambda_n, z)}{K(\lambda_n, \lambda_n)}=\sum_{n=1}^d \frac{f(\lambda_n)}{B'(\lambda_n)} \frac{B(z)}{z-\lambda_n}.
\end{equation*}
\\~\\
(b)\;
\end{comment}

Throughout this paper, the de Branges space $\HH(E)$ is identified with one in the chain of spaces $\HH(E_t)$, $0\le t<\infty$, at certain  point $t=c$ and space $\HH(E_c)$ (and the whole chain)  is contained isometrically in $L^2(\mu)$, where $\mu$ is the spectral measure of an si-process in question.  Therefore
\begin{equation}\label{eq:doublisp}
\Big\|f- \sum_{|n|\le N} f(\lambda_n)\frac{K_c(\lambda_n, \cdot)}{K_c(\lambda_n, \lambda_n)}\Big\|_\mu\rightarrow 0
\end{equation}
as $N\rightarrow \infty$ whatever element $f$ of space $\HH(E_c)$ with the reproducing kernel $K_c(w,z)$. 
The $\lambda_n$ are zeros of the odd component $B_c(z)$. 

As regards mapping  
\eqref{eq:isodesharpcp}, we note that for elements of the shifted space $e^{-izr}\HH(E_c)$, $r=\tau_c$, say $\mathsf{f}=e^{-izr} f$ where $f\in \HH(E_c)$, it holds that
\begin{equation}\label{eq:singlisp}
\Big\|\mathsf{f}- \sum_{|n|\le N} \mathsf{f}(\lambda_n)\frac{\mathsf{K}_c(\lambda_n, \cdot)}{\mathsf{K}_c(\lambda_n, \lambda_n)}\Big\|_\mu=\Big\|{f}- \sum_{|n|\le N} {f}(\lambda_n)\frac{{K}_c(\lambda_n, \cdot)}{{K}_c(\lambda_n, \lambda_n)}\Big\|_\mu\rightarrow 0
\end{equation}
as $N\rightarrow \infty$, where 
$\mathsf{K}_c(w, z)$ is the kernel in the shifted space $e^{-izr}\HH(E_c)$ and the $\lambda_n$ are zeros of the odd component $B_c(z)$. This is clear,  as the metric does not change by shifting and $\mathsf{K}_c(\lambda, \lambda)={K}_c(\lambda, \lambda)$ on the real diagonal. 

\subsection{PW-series}\label{sub:pwserp}
The sampling formula  can be used  for  calculating  covariance functions. %\eqref{eq:covrp}.
%In general, the calculations can be conducted in  $L^2(\mu)$ via the spectral representation \eqref{eq:anamova10} or  in the corresponding structure space of Theorem~\ref{thm:isophipsi} via the eigenfunction expansion of the Fourier transforms of the indicator function \widehat{1}_t$, cf. the chain of identities \eqref{eq:identnorm}. Yet another possibility arises 
If one has the set of zeros $(\lambda_n)$, obtained so as is described  above, then the covariance function of a double-sided si-process $X$  is  calculated in this manner 
\begin{equation}\label{eq:extraposp}
 \EEE (X_s \overline X_t) = \sum_n \frac{\widehat{1}_s(\lambda_n)\overline{\widehat{1}_t(\lambda_n)}}{K_c(\lambda_n,\lambda_n)}
\end{equation}
for all $-r\le s,t\le r$. This follows directly from  \eqref{eq:samplingp} applied to the particular element $\widehat{1}_t$ of the space $\HH(E_c)$,
as $|t|\le r$.
If $X$ is single-sided, then we have again formula  \eqref{eq:extraposp} for all $0\le s,t\le 2r$ in this case, since $\mathsf{K}_c(\lambda, \lambda)={K}_c(\lambda, \lambda)$, as we know.  
\begin{comment}
In  the next section we  apply the sampling formula \eqref{eq:sampling} to basic functions  
 of the single- or double-sided span $\overline{sp}\{\widehat{1}_t(z):  0\le t\le r\}_\mu$ or $\overline{sp}\{\widehat{1}_t(z):  -r\le t\le r\}_\mu$, the Fourier transforms of the indicator functions $1_{(0, t)}$. 
 
 We will need  not only the convergence of the series to each entire function $\widehat{1}_t(z)$ for a fixed point $0\le t\le r$, but 
 also the uniform convergence  over the full set of functions.
\end{comment}

%The preceding shows that if one has in his disposal  the chain of de Branges functions determining an si-process $X$, then one can take a linear combination  of
Now, $\Ifr_1$-map of
the entire functions $K_c(\lambda_n,z)/K_c(\lambda_n, \lambda_n)$ defines the 
independent Gaussian random variables $Z_n$ with zero mean and variance $\EEE|Z_n|^2=1/K_c(\lambda_n, \lambda_n)$. Since $X_t$ is $\Ifr_1$-map of $\widehat{1}_t(z)$ by definition,
the spectral isometry gives the covariances 
\[
 \EEE ({X}_t \overline Z_n)= \frac{\overline{\widehat{1}_t(\lambda_n)}}{K_c(\lambda_n,\lambda_n)}. 
\]
But then the finite linear combination 
$
 \Sigma_t^N \pd 
 \sum_{|n|\le N }
 \widehat{1}_t(\lambda_n)\,Z_n
$ 
% = \sum_{|n|\le N }\frac{e^{i\lambda_n t}-1}{i\lambda_n}\, Z_n\]
is so that
\begin{equation}\label{eq:sfzetebip}
 \EEE\Big|X_t-\sum_{|n|\le N }
 %\overline{\widehat{1}_t(\lambda_n)}
 \frac{e^{i\lambda_n t}-1}{i\lambda_n}
 \, Z_n\Big|^2\rightarrow 0 
\end{equation}
as $N\rightarrow\infty$ for each $t\in [-r,r]$ or $t\in [0,2r]$ depending whether $X$ is double- or single-sided.
For definiteness, let us focus our attention to the case of a double-sided $X$, the singles-sided case will be  handled analogously.
We call $
 \Sigma_t^N = 
 \sum_{|n|\le N }
 \widehat{1}_t(\lambda_n)\,Z_n
$  the {\it partial sum} and consider  partial processes $\Sigma^N=(\Sigma_t^N)_{t\in [-r,r]} $. 
It follows from the preceding that the finite dimensional distributions of  $\Sigma^N$ converge weakly to the finite dimensional distributions of $X$. 
To extend the Paley and Wiener result  beyond the Brownian motion, one need to show that the partial processes converge in $C[-r,r]$ with probability one, where $C[-r,r]$ is the space of continuous functions on $[-r,r]$, endowed with the supremum metric. But this is guaranteed by the Lévy-It\^o-Nisio Theorem \cite[Theorem 4.1]{Ito68} on the convergence of sums of symmetric Banach space valued random variables, cf. also 
\cite[Theorem 2.4]{Led91}. Hence, we have the following extension of the   Paley--Wiener theorem.

\begin{thm}\label{thm:Pawelp77p} 

Let   $\cdots <\lambda_{-1}< \lambda_{0}=0< \lambda_{1}<\cdots$ be the real valued zeros of the  component $B_r (z)$ of  the de Branges function generating space $\HH(E_r)$ which is contained isometrically in the symmetric space $L^2(\mu) $ where $\mu$ is the spectral measure of the double-sided si-process $X$ in question. Then the process, restricted to  interval $-r\le t\le r$,  has the series representation
\begin{equation*}
 X_t= \sum_n 
 \frac{e^{i\lambda_n t}-1}{i\lambda_n}
 \, Z_n
\end{equation*}
where the $Z_n$ are mutually independent Gaussian random variables with zero mean and variance
\[
\EEE |Z_n|^2=\sigma^2_n=\frac{1}{K_r(\lambda_n,\lambda_n) }.
\]
The series converges in mean square for each $t\in [-r,r]$. If process $X$ admits a continuous version, then with probability one, the series converge uniformly in $[-r, r]$.

\end{thm}
%By definition  $L^2(\mu)$ is symmetric if the generating spectral measure is odd, cf. Section~\ref{s:seriesABp}.  By Theorem~\ref{thm:problem175p}  space $\HH(E_r)$ is symmetric about the origin. The {\it reality condition} of Section~\ref{subs:symetryp} is satisfied and  $B_r(z)$ is an odd function of $z$,  has a zero at the origin and 
%$-\lambda_{-k}=\lambda_{k}$. 
Due to the definitions  $X^e_t=(X_t-X_t)/2$ and $X^o_t=(X_t+X_t)/2$, $t\ge 0$, for the even and odd parts of a double sided si-process $X$, the theorem implies
\begin{cor}\label{exreop}
The even and odd parts of the si-process of Theorem~\ref{thm:Pawelp77p} expand in the series
\begin{align*}
  X^e_t= t Z + 
  2\sum_{n>0} \frac{\sin \lambda_n t}{\lambda_n}\,Z^{e}_n \qquad
X^o_t=2\sum_{n>0} \frac{\cos \lambda_nt-1}{\lambda_n}\,Z^{o}_n  
\end{align*}
for each $t\in [0, r]$, in terms of independent Gaussian random variables with variances 
\[
 \EEE|Z|^2=\sigma^2_0\qquad \EEE|Z_n^e|^2=\EEE|Z_n^o|^2=\frac{1}{2}\sigma^2_n
\]
where $\sigma_n^2=1/K_c(\lambda_n,\lambda_n)$ is the same as in the theorem.

\end{cor}
\begin{prf}
One only needs to check that on calculating even and odd parts there occur random variables 
\[
 Z=Z_0 \qquad Z_n^e=\frac{Z_n+Z_{-n}}{2}\qquad Z_n^o=\frac{Z_n-Z_{-n}}{2}
\]
for $n>0$ which are independent and Gaussian with zero mean and indicated variances.  
\end{prf}
For applications, for instance, for simulation purposes,  it is important to know how quickly the remainder term
\begin{equation}\label{eq:remtp}
 R_t^N\pd X_t-\sum_{|n|\le N} 
 \frac{e^{i\lambda_n t}-1}{i\lambda_n}\, Z_n=\sum_{|n|> N} 
 \frac{e^{i\lambda_n t}-1}{i\lambda_n}
 \, Z_n
\end{equation}
vanishes as $N\rightarrow \infty$. 
It depends, of course, on the behaviour of the variances of the random variables $Z_n$ for large indices $n$. To characterized this behaviour it will be  assumed in the sequel that there exists a positive  number $a>-1$ such that  
\begin{equation}\label{eq:me78p}
 \EEE|Z_n|^2\le c\, n^{-a}
\end{equation}
for all $n\ge N$, with $N$  large enough; constant $c$ can be any positive number.
\paragraph{Remark.}
Jumping ahead, we refer to 
Lemma~\ref{thm:samplefbmp} which will  prove   that in the case of an FBM with Hurst index $H$ it holds that $\EEE|Z_n|^2\sim  n^{1-2H}$. Therefore the variances decline for $H>1/2$, stay constant for an ordinary BM and increase for $H<1/2$.
The proof is based on the well-known fact
that the $n^{\text{th}}$ zero of the Bessel function $J_{\nu+1}(rz)$ (and hence of the odd component   $B_r(z)$, cf. \eqref{eq:Jentp}) is such that $|\lambda_n|=(n\pi/r)(1+o(1))$ far out; see e.g. \cite[p. 509]{Wat44}.
By definition   $\EEE |Z_n|^2=1/K_c(\lambda_n,\lambda_n)$, 
 condition \eqref{eq:me78p} is equivalent to 
\begin{equation}\label{eq:Pawelp78p}
 K_c(\lambda_n,\lambda_n)\ge c|\lambda_n|^a
\end{equation}
for $n\ge N$. The exponent 
$a$ will determine the speed at which the remainder term $R^N$ vanishes as $N\rightarrow\infty$. 
It varies, in principal,   
with the choice of the metric with which the difference 
between the process $X$ and the partial sum $\Sigma^N$ is measured.  For instance, one may calculate  the expected values   of squired pointwise deviations, i.e. $\EEE|R^N_t|^2$,  and take maximum over the whole interval $[-r,r]$. The following estimate rather straightforward
\begin{equation}\label{eq:pointestp}
 \sup_{t\in[-r,r]}\EEE|R^N_t|^2\simeq N^{-(1+a)},
\end{equation}
where the relation symbol $\simeq$ is used to designate  there exists a positive   constant such that the left-hand side is less or equal to the constant multiple of the right hand-side. Indeed, the variances of the remainder terms for ${-r< t\le r}$ are uniformly bounded  
\[
 \EEE|R^N_t|^2\le \sum_{|n|>N} \frac{1}{|\lambda_n|^2 K_c(\lambda_n,\lambda_n)}.
  \]
When   $|\lambda_n|$ is of order $(2n\pi/r)$ for $n$ large enough, so as in the case of an FBM,  condition \eqref{eq:Pawelp78p} yields the  result
\[
\sup_{t\in[-r,r]}\EEE|R^N_t|^2\simeq  \sum_{n>N} n^{-(2+a)}\simeq\int_N^\infty  x^{-(2+a)}dx= N^{-(1+a)}/(1+a).
\]
which shows how quickly the remainder term vanishes. $\hfill$ $\qed$
\\~\\
Turning back to the general case, we shall follow  \cite{dzh05c} and make
use of the uniform norm $\|R^N\|_\infty$. The objective is to obtain an upper bound
\begin{equation}\label{eq:pawlisp}
 \EEE \|R^N\|_\infty\simeq (N^{-(a+1)}\,\ln N)^{1/2}
\end{equation}
(\cite{dzh05c} discuses  fractional Brownian sheets, however the method of proof is applicable to the present more general case too).

\paragraph{
Proof  of \eqref{eq:pawlisp}.}
For a positive integer $m$, define the difference of the partial sums
\begin{equation*}
 \Delta_m(t)\pd\Sigma^{2^m}_t-\Sigma_t^{2^{m-1}}
\end{equation*}
and estimate the expectation of the supremum over  interval $[-r,r]$ in this manner. Let the interval be covered with smaller intervals $I_1,\cdots, I_l$ with centers $t_1,\cdots, t_l$  
and length $\epsilon$,  where the length and the covering number are related as $l\simeq 1/\epsilon$ (length $\epsilon $ will   later go down with  increase of the number $m$). Clearly, $I_k=(t_{k}-\epsilon/2,t_k+\epsilon/2)$.  The supremum   over  interval $[-r,r]$ of the $m^{\text{th}}$  difference of the partial sums can be bounded by two terms  
\[
 \sup_{t\in[-r,r]} |\Delta_m(t)|\le \sup_{k=1,\cdots, l} |\Delta_m(t)|+ \sup_{k=1,\cdots, l}\sup_{s,t\in I_k} |\Delta_m(t)-\Delta_m(s)|.
\]
The expectation of the first term can be estimated by aid of a standard maximal inequality for Gaussian sequences which can be found, e.g., in \cite[Lemma 2.2.2]{vaa1}. It says
\[
\EEE \sup_{k=1,\cdots, l} |\Delta_m(t)|\simeq \sqrt{1+\ln l}
\sup_{k=1,\cdots, l} \sqrt{\EEE |\Delta_m(t_k)|^2}.
\]
The variance of the sum of independent random variables is easily calculated and bounded as follows
\[
 \EEE |\Delta_m(t_k)|^2= 
 \sum_{2^{m-1}< j\le 2^m}
 \frac
 {
 |\widehat 1_{(0,t_k]}
 (\lambda_j)|^2
 }
 {K_c(\lambda_j,\lambda_j)}
 \simeq \sum_{2^{m-1}< j\le 2^m}\frac{1}{|\lambda_j|^2K_c(\lambda_j,\lambda_j)}.
\]
For the rest, we can argue as in the case of \eqref{eq:pointestp} and obtain
\[
 \sum_{2^{m-1}< j\le 2^m}\frac{1}{|\lambda_j|^2K_c(\lambda_j,\lambda_j)}\simeq 1/ 2^{m(a+1)/2},
\]
therefore the expectation of the first term on the right-hand side of the inequality for $\sup_{t\in[-r,r]} |\Delta_m(t)|$ is bounded by
\[
 \EEE \sup_{k=1,\cdots, l} |\Delta_m(t)|
 \simeq \sqrt{1+\ln l}/ 2^{m(a+1)/2}.
\]
The expectation of the second term is bounded by
\[
 \EEE \sup_{k=1,\cdots, l}\sup_{s,t\in I_k} \sum_{2^{m-1}<j<2^m}|Z_j| \frac{|\widehat 1_{(0,t)}(\lambda_j)-\widehat 1_{(0,t)}(\lambda_j)|}{K_c(\lambda_j, \lambda_j)}
\]
where $|\widehat 1_{(0,t)}(\lambda_j)-\widehat 1_{(0,t)}(\lambda_j)|
=|\int_s^t e^{i\lambda_j u} du|\le |t-s|\le \epsilon $ for $s,t\in I_k$. By \eqref{eq:Pawelp78p}, kernel $K_r(w,z)$ at $w=z=\lambda_n$ is bounded below by a multiple of $n^a$, therefore 
\[
 \EEE \sup_{k=1,\cdots, N}\sup_{s,t\in I_k} |\Delta^m(t)-\Delta^m(s)|
 \simeq
\epsilon \sum_{2^{m-1}<j\le 2^m}j^{-a/2}\le \epsilon\, 2^{\frac{m-1}{1-a/2}}.
\]
We have estimated the expectations of both terms in the foregoing inequality for $\EEE\sup_{t\in[-r,r]}|\Delta_m(t)|$. By combining these two results, we obtain
\[
 \EEE\sup_{t\in[-r,r]}|\Delta_m(t)|\simeq
\sqrt{1+\ln N} \,2^{-m(a+1)/2} 
 +\epsilon\, 2^{\frac{m-1}{1-a/2}}.
\]
The covering number and the length of intervals are in relation $l\simeq 1/\epsilon$. We shall now specify the length so as to make the second term on the right-hand side of the latter display
of lower order than the first. This is achieved by $\epsilon=2^{-2m}$.  Hence
\begin{equation}\label{eq:paw4.1.11p}
\EEE\sup_{t\in[-r,r]}|\Delta_m(t)|\simeq
\sqrt{m}\,2^{-m(1+a)/2}.
\end{equation}
Now, let us fix the number $N$ of terms in the partial sum and select $m$ so that $2^{m-1}<N\le 2^m$. At each point $t\in [-r,r]$ the remainder term  
is bounded by
the sum of two terms
\begin{equation}\label{eq:eeap}
 |R^N_t|\le \Big|\sum_{N<j\le 2^m} \frac{\widehat 1_{(0, t]}(\lambda_j)}{K_c(\lambda_j,\lambda_j)}\,Z_j\Big|+\sum_{j>m} |\Delta_j(t)|.
\end{equation}
The second one is bounded as follows
\begin{align*}
\sum_{j>m}\EEE \sup_{t\in[-r,r]}|\Delta_j(t)| &\simeq \sum_{j>m} \sqrt{j}\,2^{-j(1+a)/2}\simeq \sqrt{m+1}\, 2^{(m+1)(a+1)/2}\\
& \simeq N^{-(a+1)/2} \,\sqrt{\ln N}
\end{align*}
(number $m$ is selected so that $2^{2m+1}<4N$). To prove this we first apply \eqref{eq:paw4.1.11p} and  then the fact that
$\sum_{j>m} \sqrt{j}\,2^{-j\, \alpha}\simeq \sqrt{m+1}\, 2^{(m+1)\,\alpha}$
for any positive number $\alpha$. As for the first term on the right-hand side of \eqref{eq:eeap}, we can make again use of 
\eqref{eq:paw4.1.11p} to get
\[
 \EEE \sup_{t\in[-r,r]}\Big|\sum_{N<j\le 2^m} \frac{\widehat 1_{(0, t]}(\lambda_j)}{K_c(\lambda_j,\lambda_j)}\,Z_j\Big|\simeq\sqrt{m}\,2^{-m(a+1)/2}\simeq N^{-(a+1)/2} \,\sqrt{\ln N}.
\]
Hence
\[\EEE \sup_{t\in[-r,r]} |R^N_t|
\simeq N^{-(a+1)/2} \,\sqrt{\ln N}
\]
which is the required estimate \eqref{eq:pawlisp}. $\hfill$ $\qed$
\\~\\
As is pointed out  earlier, Theorem~\ref{thm:Pawelp77p} regards double-sided si-processes, but  is easily extendable to single-sided processes. Because, if the 
independent  random variables $\mathsf Z_n$ are defined
as $\Ifr_1$-map of  element
 $\mathsf K_c(\lambda_n,z)/\mathsf K_c(\lambda_n, \lambda_n)$
of the space $e^{-izr}\HH(E_c)$, $r=\tau_c$ (which is  $e^{-i(z-\lambda_n)r}$ multiple of $ K_c(\lambda_n,z)/K_c(\lambda_n, \lambda_n)$),
then their Gaussian distribution is the same as that of the previous $Z_n$. Clearly, they are again centered 
and have variances $\EEE|\mathsf Z_n|^2=1/K_c(\lambda_n, \lambda_n)$. We have mentioned already that if $X$ is single-sided, then the convergence in \eqref{eq:sfzetebip} holds for $t\in [0,2r]$, and this is true with  $Z_n$  substituted by $\mathsf Z_n$. 
\\~\\ 
\noindent  
{\bf Theorem~\ref{thm:Pawelp77p} } (continued).  {\it
In the situation of the first part of the theorem,
let $\mu$ be the spectral measure of a single-sided si-process $X$. Then the process, restricted to  interval $0\le t\le 2r$,  has the series representation
\begin{equation*}
 X_t= \sum_n 
 \frac{e^{i\lambda_n t}-1}{i\lambda_n}
 \, Z_n
\end{equation*}
where the $Z_n$ are mutually independent Gaussian random variables with zero mean and variance 
\[
 \EEE |Z_n|^2=\sigma^2_n=\frac{1}{K_c(\lambda_n,\lambda_n)}.
\]
The series converges in mean square for each $t\in [0,2r]$. \\
If process $X$ admits a continuous version, then with probability one, the series converge uniformly in $[0, 2r]$.\\
If process $X$ is real, then the series turns into
\begin{align*}
  X_t= t Z + 
  2\sum_{n>0} \frac{\sin \lambda_n t}{\lambda_n}\,Z^{(1)}_n +2\sum_{n>0} \frac{\cos \lambda_nt-1}{\lambda_n}\,Z^{(2)}_n  
\end{align*}
for each $t\in [0,2r]$, with independent Gaussian random variables, centered and having  variances} 
\[
\EEE|Z|^2= \sigma^2_0 \qquad\EEE|Z_n^{(1)}|^2=\EEE|Z_n^{(2)}|^2=\frac{1}{2}\sigma^2_n.
\]
\begin{prf}
The proof follows that of the first part of theorem, except the statement regarding the real case. To proof the latter, note that in formula 
\eqref{eq:extraposp} for the covariance function the terms indexed by $n$ and $-n$ are complex conjugated. Therefore 
\begin{equation*}
 \EEE X_s \overline X_t =\sigma^2_0+ 2 \sum_{n>0} \sigma_n^2{\text{Re}\;}\big({\widehat{1}_s(\lambda_n)\overline{\widehat{1}_t(\lambda_n)}}\big)
\end{equation*}
with
\begin{align*}
{\text{Re}\;}\big({\widehat{1}_s(\lambda_n)\overline{\widehat{1}_t(\lambda_n)}}\big)=\frac{\sin \lambda_ns}{\lambda_n} \frac{\sin \lambda_nt}{\lambda_n}+\frac{\cos \lambda_ns-1}{\lambda_n} \frac{\cos \lambda_nt-1}{\lambda_n}.
\end{align*}
This shows that the series for real $X$ converges in mean squire sense and that the resulting process has the same distribution as Brownian motion. On calculating covariances we make use of the relations 
\[
 Z=Z_0 \qquad Z_n^{(1)}=\frac{Z_n+Z_{-n}}{2}\qquad Z_n^{(2)}=\frac{Z_n-Z_{-n}}{2}
\]
for $n>0$. 
\end{prf}
\subsection{PW-series for  FBM}\label{sub:PWFBMp}

(a)\, In this section  the general results of the preceding one are applied to fractional Brownian motions of Hurst index $0<H<1$. It will be shown below that these particular si-processes  are {\it determined} by 
the chain of homogeneous\fn{This term is used on the basis of the identities
\begin{align*}%\label{eq:Jentpropp}
[ A_a(z), B_a(z) ]&=[A_1(az), \,
a^{1+2\nu}B_1(az)]\nonumber\\ K_a(w, z)&=a^{2+2\nu}K_1(aw, az)          
\end{align*} 
which show  homogeneous nature  of these functions.}  de Branges spaces 
which are based on de Branges functions $E_t(z)$, $ 0\le t<\infty,$ with the components
\begin{equation}
\label{eq:Jentp}
[ A_t(z), B_t(z) ]=\Gamma(1-H)\Big(\frac{tz}{2}\Big)^{H} [J_{-H}(tz), 
t^{1-2H}J_{1-H}(tz)].         
\end{equation} 
Let us show that these components indeed define a de Branges function. 
\begin{prop}\label{thm:besdebp}
For $\nu>-1$ and $r>0$ the entire  function $E_r(z)=A_r(z)-iB_r(z)$  with $-H=\nu$ 
is a De Branges function,   normalized at the origin $E_r(0)=1$. It  is of exponential type $r$.   
\end{prop}
\begin{prf}
The normalization is clear, since
$
  z^{-\nu}J_\nu(z)\rightarrow(1/2)^\nu/\Gamma(\nu+1)$ {as}
 $
z\rightarrow 0
$
for $\nu>-1$.
For future use,  we shall prove 
more general statement:   for each $r>0$ the following entire function of 
two complex variables $w$ and $z$ has the integral representation
\begin{align}
\label{eq:kornerip}
 K_r(w,z)\pd&\frac{B_r(z)A_r(\bar w)-A_r(z)B_r(\bar w)}
{\pi (z-\bar w)}\nonumber\\
=&\frac{1}{\pi} \int_0^r  [A_t(\bar w),B_t(\bar w)]\,  dm_t\,  [A_t(z),
  B_t(z)]^\top
\end{align}
where the integration is carried out with respect to a diagonal matrix-valued function
$m_t=\text{diag}\,[\alpha_t, \gamma_t]$ with the 
entries $
 \alpha_t={t^{2\nu+2}}/(2\nu+2)$ and $ \gamma_t={t^{-2\nu}}/(-2\nu)
$
and the 
density  
\begin{equation}
\label{eq:ab1c0p}
 m'_t=
  \text{diag}\, \big[
      t^{1+2\nu},\;
      (1/t)^{1+2\nu}
     \big]
\end{equation}
for  $t>0$. To prove \eqref{eq:kornerip}, calculate the integrals with the help  of the integration and differentiation rules, cf. e.g.  \cite[formulas 5.3.5-5.3.7]{Leb72}.
%\eqref{eq:1indint3}. 
Since $\nu>-1$, one obtains  
\begin{align}\label{eq:cornerstonep}
 \int_0^r &A_t(\bar w) A_t(z) t^{1+2\nu}\,dt +\int_0^r B_t(\bar w) B_t(z) (1/t)^{1+2\nu}\,dt\nonumber
 \\&=\Gamma^2(\nu+1)\Big(\frac{4}{{\bar w}z}\Big)^\nu\,  \frac{J_{\nu+1}(rz)J_{\nu}(r \bar w)-J_{\nu}(rz)J_{\nu+1}(a \bar w)}{ (z-\bar w)/r}.
\end{align}
It is easily seen that  the right hand-side is equal to $\pi K_a(w,z)$. 
It follows that $K_r(w,w)>0$ and
\begin{align*}
\pi K_r(w,w)=
\frac{|E_r(w)|^2-|E_r^\sharp(w)|^2}{4 y} >0
\end{align*}
for each $w$ in the upper half-plane, $r>0$. Hence $|E_r(w)|>|E_r^\sharp(w)|$  in the upper half-plane and  $E_r(z)$ is a de Branges function.

We deal here with a generalization of
the classical example $E_r(z)=e^{-irz}$, for if  $\nu=-1/2$,
then $A_r(z)=\cos rz$ and $B_r(z)=\sin rz$.
As in this particular case, the positive number $r$ is a type parameter.
It follows from the inequalities  from Section 3.3 of Watson's book \cite{Wat44}  that
 $E_r(z)$ is of exponential type $r$. Indeed, those inequalities  show  
\[
|E_r(w)|\le \big(1+|w|\,{r^{2\nu+2}}/{(2\nu+2)}\big)e^{r|y|} 
\]
for every $w=x+iy$.  
The proof is complete.
\end{prf}
The proposition determines the reproducing kernels of the present chain of spaces  and the structure function through \eqref{eq:kornerip} and \eqref{eq:ab1c0p}, respectively. Indirectly, it defines also the spectral measure $\mu$ such that the whole chain is contained isometrically in $L^2(\mu)$. To see this,  associate  
with the present diagonal structure function the {\it structure space} $L^2([0,\infty), m)$  which 
consists of all pairs $h=[\varphi,\psi]$ of functions  squire integrable with respect to $dm$. Each such pair that vanishes outside of  interval $[0,r]$ defines an element $f$ of the de Branges space $\HH(E_r)$, since the expansion theorem \cite[Section 44]{deB68} says
\begin{align}
\label{eq:defogpifp}
 \pi f(w) 
=\int_0^\infty h_t\, dm_t q(t,w)^\top
\end{align}  
for all complex numbers, where $q(t,w)=[ A_t(w), B_t(w)]1_{(0,r)}(t)$ so as in \eqref{eq:eigentrp}. Note use of the identity $\tau_t=t$ which is clear in view of \eqref{eq:ab1c0p}. The norm of element $f$ in the metric of $\HH(E_r)$ is given by Parseval's identity \eqref{eq:neigentrp} with $c=r$. Now, it is to be shown that 
the right-hand side of the latter identity, which is the squire norm of $h$ in the metric of $L^2([0,\infty), m)$, is calculated also as  
\begin{align}\label{eq:spectridaphipsip}
\int_0^\infty h\,dm\, h^\top=2\pi\int_0^\infty |f(\lambda)|^2\mu(d\lambda)
\end{align}
where $\mu(d\lambda)$ is given by 
\begin{equation}\label{eq:spshicp}
 \mu(d\lambda)= \mu'(1) |\lambda|^{1+2\nu}d\lambda
\qquad \qquad
\mu'(1)=\frac{\pi}{2^{1+2\nu} \Gamma^2(1+\nu)}.
\end{equation}
This is achieved by considering the eigenfunction representation \eqref{eq:defogpifp} (and its particular case \eqref{eq:kornerip} for $f(z)=K_a(w,z)$ and $h_t=q(t,w)$) as a certain Hankel transform by spelling it out as follows 
\begin{equation}
\label{eq:gagrdzp}
  \frac{\pi(z/2)^\nu}{\Gamma(\nu+1)} f(z)=
\int_0^r\varphi_t J_\nu(tz)t^{\nu+1}dt+\int_0^r\psi_t J_{\nu+1}(tz)t^{-\nu}dt.
\end{equation}
This shows that the entire function 
$
{\pi(z/2)^{\nu+\frac12} f(z)}/{\Gamma(\nu+1)} 
$
forms the pair of Hankel transforms of order $\nu$ and $\nu+1$ with the functions which vanish outside of the interval $[0,r]$ and are defined within the interval as 
$t^{\nu+\frac12}\varphi_t/{\sqrt 2}$ and  
$
{\psi_t}/(\sqrt 2\; t^{\nu+\frac12})$, respectively. By Parseval's theorem for Hankel transforms one obtains \eqref{eq:spectridaphipsip}.

Hence the squire norm in the metric of de Branges spaces $\HH(E_r)$ is equal to that of $L^2(\mu)$ with the spectral measure $\mu$ defined by \eqref{eq:spshicp}. For Hurst index $0<H<1$, parameter $\nu$ is confined to interval $(-1, 0)$ and so the spectral measure is of property \eqref{eq:sqipaper}. Upon comparing the latter with the spectral measure \eqref{eq:cHstvisp} for the standard  FBM, we see the difference of constant multipliers. Actually, the spectral measure \eqref{eq:cHstvisp} defines a non-standard FBM, with  
$
 r(1,1)=\mu'(1) \Gamma(1-2H) {\cos H\pi}/{H}.
$
The necessary changes caused by this difference are made by virtue of rule 1 in \cite[Section 6.6]{dym76} and \cite[Section 7.2.1]{Dzh05}.  
But this will  be unnecessary in the sequel.

The  space $\HH(E_r)$ is defined by the components \eqref{eq:Jentp} with $t=r$ and the homogeneity order $\nu>-1$ is related to  Hurst index $0<H<1$ by $\nu=-H$. In the sequel, the zeros $\lambda_n$   will be that of the second component in \eqref{eq:Jentp}, i.e.  zeros of the Bessel function $J_{1-H}(z)$  (we remark that alternatively one can choose to work with the zeros of the first component -- one can easily verify that all consequent results of this chapter remain true, like in the case $H=1/2$ in which   expansions are obtainable in terms of the zeros of the sine or cosine, cf. \cite[Section 26.1]{yag87}). 
%The forthcoming theorem verifies condition \eqref{eq:Pawelp78}. 
%For this purpose let us recall     
Regarding the zeros of the Bessel function of the first kind of order $\nu>-1$, see \cite[Section 7.9]{Erd53_2}. The Bessel function $J_\nu(z)$ has a countable number of positive zeros that can be arranged in ascending order of magnitude. We denote them by $\lambda_{\nu,1}<\lambda_{\nu,2}<\cdots$. For positive $\nu$, the Bessel function vanishes at the origin $J_\nu(0)=0$, and its negative zeros are given by $-\lambda_{\nu,1}>-\lambda_{\nu,2}>\cdots$. Hence, for $\nu\ge 0$ the zeros can be ordered as  $\cdots <\lambda_{\nu,-1}<\lambda_{\nu,0}=0< \lambda_{\nu,1}<\cdots$.

Making use of  these properties of the zeros and the form of the reproducing kernel \eqref{eq:cornerstonep},  we can  verify condition \eqref{eq:Pawelp78p}. 

\begin{lem}\label{thm:samplefbmp}
If 
the de Branges space $\HH(E_r)$ is defined by the components \eqref{eq:Jentp} 
%  and its kernel by \eqref{eq:nustvis}, where $a=r$ and   $\nu=-H$, 
and if $\cdots \lambda_{-1}< \lambda_0=0< \lambda_{1} <\cdots $ are the zeros of $J_{1-H}(rz)$, then  
condition \eqref{eq:Pawelp78p} 
 is satisfied with $a=1-2H$. 
\end{lem}
\begin{prf}
The  reproducing kernel in \eqref{eq:Pawelp78p}  is
$K_r(w,z)$, since  in the present case $\tau_t=t$.
On the diagonal $w=z=\lambda$, with $\lambda$ real, 
the  reproducing kernel is 
equal to $A_r(\lambda) B_r'(\lambda)-A'_r(\lambda) B_r(\lambda)$ divided through $\pi$, e.g.
\begin{align*}
%\label{eq:nulshich+p}
 \frac{K_r(\lambda,\lambda)}{K_r(0,0)}&=(2-2H)\Gamma^2(1-H) (r\lambda/2)^{2H}\nonumber 
 \\
 &\times\Big(J^2_{1-H}(r\lambda)-\frac{1-2H}{r\lambda}J_{-H}(r\lambda)J_{1-H}(r\lambda)+J^2_{-H}(r\lambda)\Big)
\end{align*}
by  \eqref{eq:cornerstonep} and the properties %\eqref{eq:1AppdifBp} and \eqref{eq:Hrecnulshi1p} 
of the Bessel functions.
Substituting  $\lambda$ by the $n^{\text{th}}$ zero of the Bessel function $J_{1-H}$ we get 
\begin{align}
\label{eq:nulshichp}
 \frac{K_r(\lambda_n,\lambda_n)}{K_r(0,0)}=(2-2H)\Gamma^2(1-H) (r\lambda_n/2)^{2H}
 J^2_{-H}(r\lambda_n).
\end{align}
It is shown in \cite[Section 7.1]{Wat44} that
the initial terms  in the expansion of the Bessel functions of the first kind are such that for order $\nu>-1$ and $\nu+1$  it holds that
  \[
   J^2_{\nu}(z)+J^2_{\nu+1}(z)  \sim  \frac{2}{\pi z},
  \]
see Lommel's formula on p. 200. 
This  asymptotic behaviour for large $|z|$ shows that if we substitute $z$ by the $n^{\text{th}}$ positive zero $\lambda_n$ of the Bessel function $J_\nu(z)$, then    $J_{\nu+1}(\lambda_n) \sim 2/(\pi \lambda_n)$ for 
$n\rightarrow \infty$, as the zeros 
$\lambda_n$ tend to infinity. Hence \eqref{eq:nulshichp} implies  $ K_r(\lambda_n,\lambda_n)\sim c \, |\lambda_n|^{2H-1} $ with $c$ a constant.
The proof is complete.
\end{prf}

The lemma provides a key argument for extending  the PW-expansion to fractional Brownian motions. It proves that the series of the next theorem converges towards an FBM with  
rate that is optimal in the sense of K\"uhn and Linde \cite[Theorem 7.4]{Kuh02}.
\begin{comment}
of the preceding proof is based on  estimates  of the asymptotic behavior of the Bessel functions and their roots.
\begin{lem}\label{lem:kvevit}
 Let the $\lambda_n$ be zeros of the Bessel function of the first kind
 of order $\nu>-1$. Then for all $p>0$
 \[
  \sum\frac{1}{|\lambda_n|^{p+2} J^2_{\nu+1}(\lambda_n)}<\infty.
 \]
 \end{lem}
 \begin{prf}
As is shown in \cite[Section 7.1]{Wat44}
the initial terms  in the expansion of the Bessel functions of the first kind are such that for order $\nu>-1$ and $\nu+1$  it holds that
  \[
   J^2_{\nu}(z)+J^2_{\nu+1}(z)  \sim  \frac{2}{\pi z},
  \]
see Lommel's formula on p. 200. 
This  asymptotic behaviour for large $|z|$ shows that if we substitute $z$ by the $n^{\text{th}}$ positive zero $\lambda_n$ of the Bessel function $J_\nu$, then    $J_{\nu+1}(\lambda_n) \sim 2/(\pi \lambda_n)$ for 
$n\rightarrow \infty$, as the zeros 
$\lambda_n$ tend to infinity. Hence it suffices to show 
\[
  \sum\frac{1}{|\lambda_n|^{p+2} }<\infty.
 \]
But this follows by calculations of the large zeros of cylindrical functions in \cite[Section 15.53]{Wat44} showing  the $n^{\text{th} }$ positive zero $\lambda_n$ of $J_\nu$ is assymptotically of order $n\pi$.  
 \end{prf}
 
 This is rather difficult problem, in general, and we focus our attention to relatively simple case of FBM's.      

%of Hurst index $0<H<1$ 
We shall show that the latter condition is satisfied  if the chain of de Branges spaces is homogeneous in the sense of Section~\ref{s:f} 
\end{comment}

\begin{thm}\label{thm:4.5PTRFp} 
For an arbitrary  Hurst index $0<H<1$, let   $\cdots <\lambda_{-1}< \lambda_{-1}=0< \lambda_{1}<\cdots$ be zeros of the Bessel function $J_{1-H}(rz)$. Let the $Z_n$ be     
independent, complex valued Gaussian random variables with mean zero and variance $\EEE|Z_n|^2=\sigma_n^2$ where the $\sigma_n^2=1/K_r(\lambda_n, \lambda_n)$ are given by
 \eqref{eq:nulshichp}.
Then, with probability one, the series 
 \begin{equation*}
 \sum_n 
 \frac{e^{i\lambda_n t}-1}{i\lambda_n}
 \, Z_n
\end{equation*}
converges uniformly in $t\in [-r,r]$ {\rm (or $t\in [0,2r]$)} and defines a complex valued double-sided {\rm (or single-sided)} FBM with Hurst index $H$.

The remainder term \eqref{eq:remtp}
satisfies \eqref{eq:pointestp} and \eqref{eq:pawlisp} with $a=2H-1$, i.e.
\begin{align*}
 \sup_{t}\EEE|R^N_t|&\simeq N^{-H}\\ 
 \EEE \|R^N\|_\infty&\simeq \sqrt{\ln N}\,N^{-H},
\end{align*}
hence the convergence is rate-optimal. 
\end{thm}
With increase of Hurst index,  sample paths of an FBM become smoother and the above estimates tell us that
the convergence rate does improve and less number of terms are needed to achieve a desired level of approximation.

\subsection{PW-series for fundamental martingales}\label{sub:asspp}

%In Section~\ref{s:fm&ma}  an si-process $X$ and its first chaos $\MMMM_1(r)$, $ r>0$, is associated with a number of processes whose series expansions  can be  derived by the general formulas \eqref{eq:doublis} and \eqref{eq:singlis}.  The  expansions discussed below will have the same convergence properties as process $X$ itself: the series will converge in mean square and  if process in question admits a continuous version, then with probability one the convergence will be uniform. For brevity, we do not point this out  each time,  providing only formal pointwise series.   

%\subsubsection{Fundamental martingales} 
Section~\ref{s:fm&map} introduces so-called {\it fundamental martingales}  associated with  si-processes $X$. In the case of a double-sided $X$, $\Ifr_1$-map \eqref{eq:alsnogp}  defines two martingales, the even and odd one, belonging to the even and odd part of the first chaos  
$\MMMM^e_1(r)$ and $\MMMM^o_1(r)$, respectively.
Both processes $M^e$ and $M^o$, restricted to  interval $[0, r]$, expand in  PW-series. This is clear from the general considerations of Section~\ref{sub:pwserp}.  
The same can be said regarding  the single-sided case in which  
the fundamental martingale is defined  as $\Ifr_1$-map 
\eqref{eq:eigentMtp}. Restricted to interval $[0, 2r]$, process $M$ expands in PW-series, as well.
In view of  \eqref{eq:doublisp} and \eqref{eq:singlisp} we can make the following statement.
\begin{thm}\label{thm:wcmsp}
In the situation of Theorem~\ref{thm:Pawelp77p}\\ 
(i)\; The even and odd fundamental martingales associated with the given double-sided si-process $X$ by \eqref{eq:alsnogp}, expand in series
\begin{align*}
 M_t^e =\sum_n \frac{B_t(\lambda_n)}{\lambda_n}\, Z_n\qquad
 M_t^o =\sum_n \frac{A_t(\lambda_n)-1}{\lambda_n}\, Z_n
\end{align*}
for $t\in [0,r]$,
with the zeros $\lambda_n$ and the random variables $Z_n$ as defined in Theorem~\ref{thm:Pawelp77p}. \\
(ii)\, If $X$ is single-sided, then the associated fundamental martingale defined by \eqref{eq:eigentMtp} expands in series
\begin{align*}
 M_{t} =\sum_n \frac{\mathsf B_t(\lambda_n)}{\lambda_n}\, Z_n
\end{align*}
for $t\in [0,2r]$,
with the zeros $\lambda_n$, the random variables $Z_n$ as above, 
and  ${\mathsf B_t(z)}/{z}=e^{iz\tau_t} B_t(z)/z=\pi\,\mathsf K_t(0,z)$. 
%cf. \eqref{eq:eR0AR0Bp}.

\end{thm}

In %\cite{gas07} Section~\ref{sub:Mb}  
the forthcoming Section~\ref{s:KLp} we shall be interested in so-called  {\it martingale bridges} associated with fundamental martingales like in  \cite{gas07}. For instance, the bridge associated with an even fundamental martingale is defined and expanded in PW-series as follows
\begin{align*}
 M(t;r)\pd M^e_{t} -\frac{\langle M^e\rangle_t}{\langle M^e\rangle_r}M^e_r=2\sum_{n=1}^\infty \frac{B_t(\lambda_n)}{\lambda_n}\, Z^e_n
\end{align*}
where  $Z^e_n=(Z_n+Z_{-n})/2$ and $\langle M^e\rangle_t=\alpha_t.$  To deduce this from  Theorem~\ref{thm:wcmsp}, note that in the expansion for $M_t^e$ the term  indexed $n=0$ is equal to
$\alpha_t Z_0$, since $B_t(\lambda_0)/\lambda_0=\alpha_t$, as $\lambda_0=0$. Now, 
evaluate both sides of the expansion at $t=r$ to  get 
$M_r^e=\alpha_r Z_0$. Hence, in the expansion the random variable $Z_0$ can be substituted by $M_r^e/\alpha_r $.  
The right side of the display follows by the symmetry about the origin of the zeros $\lambda_0=0$, $\lambda_{-n}=\lambda_n$, as the entire function $B_t(z)/z$ is even.

\subsection{Stationary processes}\label{s:stPWep}
The sampling method applies to 
stationary processes $(Y_t)_{t\in \RR}$  without essential changes. Since the covariance function $\EEE(Y_s Y_t)=r(s,t)$ in this case depends only on the difference $t-s$ and has the spectral representation 
\[
 r(s,t)=\int e^{-i\lambda (t-s) } \mu(d\lambda) \qquad s,t\in \RR
\]
with a finite spectral measure $\mu$, the random variable $Y_0$ has a non-zero variance equal to $\mu(\RR)<\infty$. 

In analogy to \eqref{eq:sourcep}, 
the random variable $Y_t$ is brought in   connection with exponential $e^{itz}$
in the metric of $L^2(\mu)$ by isometry $\Ifr_1 e^{it\,\cdot}=Y_t$, and the first chaos 
 $\MMMM_1(r)$ in connection  with $\Ifr_1$-map of 
the closed span  $\overline{sp}\{e^{iz t}:  0\le t\le r\}_\mu$.  This single-sided 
span of exponentials 
coincides with the shifted space $e^{-izr/2}\HH(E_{r/2})$, where $\HH(E_r)$ is a de Branges function of exponential type $r>0$ contained isometrically in $L^2(\mu)$ 
(since spectral measures are finite, constants belong to  space $\HH(E_r)$ and type $r=0$   is not excluded; but this is irrelevant in the current case). Let $B_{r}(z)$ be the odd component of the de Branges function generating space $\HH(E_r)$.
Let $\mathsf K_r(w,z)$ be kernel \eqref{eq:irvelia9p} of the shifted space $e^{-izr}\HH(E_{r})$ (for simplicity, we restrict out attention to the particular case $  
\tau_t=t$, since this will be the only case in the examples we have in mind). In this situation, the sampling formula applies to each element $e^{-izt}$, $t\in [0,2r]$, of the shifted space $e^{-izr}\HH(E_{r})$ 
\[
e^{-izt}=\sum_n e^{-it\lambda_n}\frac{\mathsf K_r(\lambda_n, z)}{\mathsf K_r(\lambda_n,\lambda_n)} 
\]
where the $\lambda_n$ are zeros of $B_{r}(z)$. Now we can make use  of spectral isometry like in Theorem~\ref{thm:Pawelp77p}
and  obtain the series expansion 
\begin{equation}\label{eq:expstp}
 Y_t= \sum_n 
 e^{i\lambda_n t}
 \, \mathsf Z_n
\end{equation}
where the $\mathsf Z_n$ are mutually independent Gaussian random variables with zero mean and variance
\[
\EEE |Z_n|^2=\sigma^2_n=\frac{1}{\mathsf K_r(\lambda_n,\lambda_n) }.
\]
The series converges in mean square for each $t\in [0,2r]$.
Let us consider  two examples.

\subsubsection{OU process}
Let $Y$ be the {\it Ornstein-Uhlenbeck} (OU) process,  centered stationary Gaussian process $Y$ having the covariance function
\begin{equation}\label{eq:covOUp}
 \EEE(Y_sY_t)=\frac{\sigma^2}{2\theta} e^{-\theta |t-s|}
\end{equation}
and the density of the spectral measure
\begin{equation}\label{eq:10ship}
\mu(d\lambda)=\frac{\sigma^2}{2\pi}\frac{d\lambda}{\lambda^2+\theta^2}
\end{equation}
with the parameters $\sigma^2, \theta>0$. 
The space $L^2(\mu)$ generated by this finite spectral measure contains isometrically all spaces based on the de Branges functions 
\begin{equation}\label{eq:parelispp} 
E_t(z)=\sqrt{2\pi/\sigma^2}(\theta-iz)e^{-izt}\qquad t\ge 0.
\end{equation}
The identity of the norms in the metric of all spaces $\HH(E_t)$ for $t\ge 0$ and $L^2(\mu)$ is obvious.
At the origin space $\HH(E_0)$ is non-trivial and to construct the chain of spaces which starts from a trivial one,  
\cite[Section 7.2.4] {dzh25} lets the chain to start from a certain negative point $t_->0$ and lets $E_t(z)$ to be $\sqrt{2\pi/\sigma^2}$ multiple of 
$(\theta-i z\,\frac{t-t_-}{0-t_-})\;1_{(t_-,0)}(t)$. This part of the chain is not contained isometrically in $L^2(\mu)$, but this is beyond  our concern, since our objective is to specify series \eqref{eq:expstp} for positive values of $t$. For this purpose one needs to show
 that the shifted space $e^{-izr}\HH(E_{r})$ which is contained isometrically in $L^2(\mu)$ has
the reproducing   kernel 
\begin{align}
 \label{eq:fktwzsp}
\pi {\mathsf K}_{a}(w,z)= \pi {\mathsf K}_{0}(w,z)+ 
\int_0^{a}\overline{{\mathcal E}_{t}(w)}{\mathcal E}_{t}(z) d\alpha_t
\end{align}
%\eqref{eq:fktwzs}  
(in contrast with \eqref{eq:+rvirveliap} there occur the additional term $ \pi {\mathsf K}_{0}(w,z)$, due to non-triviality of the de Branges space at the origin) and that in the present case this gives
\begin{equation}
\label{eq:+Harryp.42p}
\mathsf{K}_{r/2}(w,z)=\frac{2\theta}{\sigma^2} +\frac{(\theta+i\bar w)(\theta-iz)}{\sigma^2}\frac{e^{-i(z-\bar w)r}-1}{-i(z-\bar w)}
\end{equation} 
for $r\ge 0$. Check this by substituting \eqref{eq:parelispp} into definition \eqref{eq:peofBp}.
%given by \eqref{eq:+Harryp.42}. 
It follows that 
\[
\sigma^2\, \mathsf K_r(\lambda,\lambda)=2\theta^2+(\theta^2+\lambda^2)\, r.
\]
Due to the form  of the odd component $B_r(z)$ of the de Branges function \eqref{eq:parelispp}, the zeros $\lambda_n$ are the roots of  equation 
$$\theta\tan (\lambda r)+\lambda=0.$$ 
With these zeros, the  OU-process  expands  in  series \eqref{eq:expstp} for $t\in [0,2r]$ where
the $Z_n$ are mutually independent Gaussian random variables with zero mean and variance 
\[
\sigma_n^2=  \, \frac{{\sigma^2}}{2\theta^2+(\theta^2+\lambda_n^2)\, r}.
\]
The variances vanish with rate $n^{-2}$.  This gives estimates
\begin{align*}
 \sup_{t}\EEE|R^N_t|&\simeq N^{-1/2}\\ 
 \EEE \|R^N\|_\infty&\simeq \sqrt{\ln N}\,N^{-1/2}
\end{align*}
and shows  
the rate is optimal, see Section~\ref{s:n&a}, Note 1. 

\subsubsection{Autoregressive process }
The
{\it autoregressive processes} $Y$ of $n^{\text{th}}$-order  
is a cantered stationary Gaussian process with the 
spectral measure whose  density is 
\begin{equation}\label{eq:muar1p}
 \mu(d\lambda)= \frac{\sigma^2}{2\pi}\,\frac{d\lambda}{|\Theta(i\lambda)|^2}
\end{equation}
where 
\begin{equation}
\label{eq:Theta11p}
\Theta(iz)=(iz-\phi_1)\cdots (iz-\phi_n)  
\end{equation} is  
a polynomial of degree $n$. Its zeros have negative real part and, moreover, 
 $\phi_k>0$, since process $Y$ is real. Like in the previous example, it is easily seen that
space $L^2(\mu)$ generated by this finite spectral measure contains isometrically all spaces based on the de Branges functions 
\begin{equation}\label{eq:parelisap} 
E_t(z)=\sqrt{2\pi/\sigma^2}\Theta(iz) e^{-izt}\qquad t\ge 0.
\end{equation}
For $r>0$,  the odd component $B_r(z)$ is such that its zeros $\lambda_l$ are the roots  of equation
\[
 \tan(r\lambda)\, \text{Re}\, \Theta(i\lambda) = \text{Im}\, \Theta(i\lambda). 
\]
The reproducing   kernel \eqref{eq:fktwzsp} in the shifted space $e^{-izr}\HH(E_{r})$ is calculated like \eqref{eq:+Harryp.42p}  for $n=1$. The later  generalises  to all positive integers as  
\begin{align*}
%\label{eq:++Harryp.42}
\sigma^2\mathsf{K}_0(w,z)&=\frac{\Theta^\sharp(iz)\Theta(i\overline w)-\Theta(iz)\Theta^\sharp(i\bar w)}{-i(z-\bar w)}\\
\sigma^2\mathsf{K}_r(w,z)&=\sigma^2\mathsf{K}_0(w,z) +{\Theta(i\overline w)\Theta^\sharp(iz)}\,\frac{e^{-i(z-\bar w)r}-1}{-i(z-\bar w)},
\end{align*}
therefore
\begin{equation*}
%\label{eq:++Harryp.42}
\sigma^2\mathsf{K}_r(\lambda_l,\lambda_l)=\sigma^2\mathsf{K}_0(\lambda_l,\lambda_l) +|\Theta(i\lambda_l)|^2\,r
\end{equation*} 
for $r\ge 0$. Since  $\phi_k>0$ by assumption, we have
$
|\Theta(i\lambda_l)|^2
=
(\lambda_l^2+\phi_1^2)\cdots (\lambda_l^2+\phi_n^2)
$
 and 
\[
\sigma^2\mathsf{K}_0(\lambda_l,\lambda_l)= 4\,\text{Re\;}\Theta(i\lambda_l)\, \text{Im\;}\Theta'(i\lambda_l).
\]
The series expansion
\begin{equation*}
 Y_t= \sum_l 
 e^{i\lambda_l t}
 \, \mathsf Z_l
\end{equation*}
for $0\le t\le r $ is in terms of  mutually independent Gaussian random variables $Z_l$ with zero mean and variance
\[
 \sigma_n^2= \frac{\sigma^2}{4 \,\text{Re\;}\Theta(i\lambda_n) \text{Im\;}\Theta'(i\lambda_n)+ |\Theta(i\lambda_n)|^2\,r}.
\]

\section{
The Karhunen--Lo\`{e}ve  expansion}\label{s:KLp}
The generalized  PW-series for fundamental martingales provided in Section~\ref{sub:asspp} posses an extra feature which will be brought forward in the present section. To get the idea, let us apply Theorem~\ref{thm:wcmsp} to a standard Brownian motion. In the theorem the $\lambda_n$ are zeros of $\cos(rz)$, therefore the positive zeros are   
$\lambda_n=(n+\frac12)\pi/r$. The independent Gaussian  random variables $Z_n$ are $N(0,1/r)$, since in this particular case $\pi K_r(\lambda,\lambda)=1/r$ for all real $\lambda$. As the zeros of $\cos(rz)$ are symmetrically spread about the origin,    
Theorem~\ref{thm:wcmsp} gives
\begin{align}\label{eq:V.2.17p}
   W_t
=\sqrt{2/r}\sum_{n=0}^{\infty}\frac{\sin((n+\frac{1}{2})\pi t/{r})}{(n+\frac{1}{2})\pi/r}\,\xi_n
\end{align}
for $0\le t\le r$, where the $\xi_n$ are mutually independent standard normal random variables. If  $\lambda_n$ are zeros of $\sin(rz)$, i.e.    
$\lambda_n=n\pi/r$, then the remark following Theorem~\ref{thm:wcmsp} entails 
\begin{align}\label{eq:V.2.18p}
 W_t
 -\frac{t}{r} W_r
=\sqrt{2/r}\sum_{n=1}^{\infty}\frac{\sin(n\pi t/{r})}{n\pi/r}\,\xi_n
\end{align}
for $0\le t\le r$, 
with the same $\xi_n$ as above. The same results follow from Theorem~\ref{thm:4.5PTRFp} with $H=1/2$, of course. 
The extra feature we are speaking about is {\it bi-orthogonality} of  these expansions as $H=1/2$,  in the sense that  set $(\xi_n)$ of random variables form an orthogonal basis in $L^2(\Omega,\FF, P)$ and
the  trigonometric functions in both expansions form an orthogonal basis in $L^2[0,r]$.   This is well-known fact,   on p. 37 of  \cite{Hig77} e.g. one can find two examples of such sets of functions
\begin{align*}
&\{1/\sqrt{r},\;  \sqrt{2/r} \cos (n\pi t/r),\; n=1,2,\cdots\}\\
 &\{ \sqrt{2/r} \sin (n\pi t/r), \; n=1,2,\cdots\}.
\end{align*}
The expansions \eqref{eq:V.2.17p} and \eqref{eq:V.2.18p} are also well-known. They are  classical examples of  the {\it Karhunen--Lo\`{e}ve  expansion}, 
cf. for instance, 
\cite[Section V.2, formulas (17) and (18)]{gikh69}. This is how Lo\`{e}ve formulates his theorem.
\begin{thm}{\rm (Lo\`{e}ve \cite[Section 34.5 (b)]{loe63}).}\label{thm:Loevep} 
A random process $X_t$, continuous in the mean squire on a closed interval $[a,b]$, has an orthogonal decomposition 
\begin{equation}\label{eq:Loevep}
 X_t =\sum_n \omega_n \varphi_n(t)\,\xi_n
\end{equation}
if and only if $|\omega_n|^2$ and $\varphi_n(t)$ are eigenvalues and eigenfunctions of its covariance function. The series converge in the mean squire uniformly on $[a,b]$.
\end{thm}
The random variables $\xi_n$ are uncorrelated $\EEE \xi_m\bar \xi_n=\delta_{mn}$. If $k(s,t)$ is the covariance function, then
\[
 \int_a^b k(s,t)\varphi_n(s) ds=|\omega_n|^2\,\varphi_n(t)
\]
and the series 
\[
 k(s,t)=\sum_n |\omega_n|^2 \varphi_n(s)\overline{\varphi_n(t)}
\]
converges on $[a,b]\times [a,b]$ uniformly and absolutely. The proof is based on the so-called {\it Mercer's theorem} (cf. e.g. \cite{lim08}),  which will be discussed in the next section.\\~\\
\noindent
Let us turn back to the examples \eqref{eq:V.2.17p} and \eqref{eq:V.2.18p}. The covariance functions of the Brownian motion and Brownian bridge are $k(t,s)=s\wedge t$ and $k(t,s)=s\wedge t- st/r$, respectively. Since these real processes are  defined on interval $[0,r]$, the integral equation of Theorem~\ref{thm:Loevep} with these covariance functions is 
\[
 \int_0^r k(s,t)\varphi_n(s) ds=\omega^2_n\,\varphi_n(t)
\]
where the $\omega_n$ are positive numbers. The functions $\varphi_n$ obey the boundary condition $\varphi_n(0)=0$. Another boundary condition is found upon  
differentiating both sides with respect to $dt$. It will be seen that 
if $k(t,s)=s\wedge t$, then $\varphi'_n(r)=0$, and if  $k(t,s)=s\wedge t-ts/r$, then $\varphi_n(r)=0$. Differentiating once more we arrive at the second order differential equation $\omega^2_n\,\varphi''_n(t)=-\varphi_n(t)$. The normalized solution of this equation with the preceding boundary conditions are well-known
\begin{align*}
\begin{matrix}
  \sqrt{2/r} \sin \big((n+\tfrac12)\frac{\pi t}{r}\big), &\omega_n^{-1}=(n+\tfrac12)\pi, & n=0,1,\cdots\\
 \sqrt{2/r} \sin \big(n\frac{\pi t}{r}\big), & \omega_n^{-1}=n\pi, & n=1,2,\cdots 
\end{matrix}
\end{align*}
respectively. This obtains the Karhunen--Lo\`{e}ve  expansions \eqref{eq:V.2.17p} and \eqref{eq:V.2.18p}.\\~\\
\noindent
The arguments which  extend the preceding  results to a class of fundamental martingales and martingale bridges are quite similar. As is said at the beginning of this section, we have already generalized PW-series for these real valued processes and it is to be shown they are  KL-expansions. To describe this class, we first invoke the definitions of Section~\ref{s:fm&map}.  For a given  si-process $X$ and its first chaos $\MMMM_1(r)$, $r>0$, an associated 
fundamental martingale is  defined in terms of a chain of de Branges spaces of exponential type  which  {\it determines} the given si-process and its first chaos. The chain is assumed to be contained isometrically in $L^2(\mu)$, where $\mu$ is the symmetric spectral measure of $X$.
It will be assumed throughout the present section that the following proposition holds true. 
\begin{thm}\label{thm:convwedgep}
Let $\alpha_t$ be a continuously differentiable strictly increasing function, $\alpha_0=0$ and $\alpha'_t>0$.  Then the system of integral equations
\begin{align}\label{eq:difchp}
1-A_t(z) =z\int_0^t B_u(z) d\gamma_u\qquad
B_t(z)=z\int_0^t A_u(z) d\alpha_u
\end{align}
with $\gamma'_t=1/\alpha'_t$, 
defines 
a chain of de Branges functions $E_t(z)
=A_t(z)-iB_t(z), 
t\ge 0,$ starting from $E_0(z)=1$, normalized at the origin $E_t(0)=1,$ and
such that the reproducing kernel in the corresponding chain of spaces  satisfies
\begin{align}
\label{eq:korneri11p}
 \pi K_t(w,z)&=\frac{B_t(z)A_t(\bar w)-A_t(z)B_t(\bar w)}
{z-\bar w}\nonumber\\
=& \int_0^t  A_u(\bar w)A_u(z)\,  d\alpha_u
+\int_0^t  B_u(\bar w)B_u(z)\,  d\gamma_u
\end{align}
for $t\ge 0$ and all complex numbers $w$ and $z$. Both integrals on the right-hand side of the latter equation satisfy the so-called 
{\it Lagrange identities} 
\begin{align}\label{eq:orilagp}
\int_0^t A_u(\bar w) A_u(z)\,d\alpha_u= \frac{zA_t(\bar w)B_t(z)-\bar w  A_t(z) B_t(\bar w)}{z^2-{\bar w}^2} \nonumber\\
\int_0^t B_u(\bar w) B_u(z)\,d\gamma_u= \frac{\bar w  B_t(z) A_t(\bar w)-zB_t(\bar w)A_t(z)}{z^2-{\bar w}^2}
\end{align}
for all complex numbers $w$ and $z$. If these numbers are real and equal, then
\begin{align}\label{eq:orinorp}
\int_0^t |A_u( w)|^2 \,d\alpha_u= \frac{w\big(A_t( w)\dot B_t(w)-  \dot A_t(w) B_t( w)\big)+A_t(w) B_t(w)}{2w} \nonumber\\
\int_0^t |B_u( w)|^2 \,d\gamma_u= \frac{w\big(A_t(w)\dot B_t(w)-  \dot A_t(w) B_t( w)\big)-A_r(w) B_t(w)}{2w}
\end{align}
for $w\ne 0$. 
\end{thm}
The derivatives $\dot A_t(w)$ and $\dot B_t(w)$ are with respect to the variable $w$.
\begin{prf}

For the integral equations \eqref{eq:difchp} and 
%, see \cite[Section 37]{deB68} and 
representation \eqref{eq:korneri11p}
of the reproducing kernel in the de Branges space $\HH(E_r)$, see \eqref{eq:intforABp} and \eqref{eq:KaKbp}.
%\cite[Problem 138]{deB68}. The proof of the latter problem can be found in \cite[Theorem 6.2.1]{dzh25}. 
It is taken here into consideration that under the present conditions 
the chain of symmetric de Branges spaces %defined by \eqref{eq:difchp}
starts from a trivial space and $E_t(0)=1$ for all $t\ge 0$.\fn{
In \cite{dym71} H. Dym discusses the equivalence of  the integral equations \eqref{eq:difchp} to the second order differential equations 
\begin{align}
\label{eq:stluzp}
 \frac{d^2A_t(z)}{d\alpha_t\,d\gamma_t} 
 =-z^2 A_t(z)\qquad 
 \frac{d^2B_t(z)}{d\gamma_t\,d\alpha_t} =-z^2 B_t(z)
\end{align}
relating the latter to classical differential equations of Sturm--Liouville type. We shall utilize this relationship in the course of proving Theorem~\ref{thm:wedgemercp}.
} 
%The representation \eqref{eq:korneri11p}
%of the reproducing kernel in the de Branges space $\HH(E_r)$ is seen from Theorem~\ref{thm:thm40p}, items $(v)$ and $(iv)$; 
%cf. \eqref{eq:korneri} and \eqref{eq:problem139}.

For the
{\it Lagrange identities}, see \cite[Lemma 1.1]{Kack74}. 
It suffices to provide the proof of the first of the equations \eqref{eq:orilagp}, since the second one follows from \eqref{eq:korneri11p}. 
Making use of $[dB_t(z), -dA_t(z)]=z [A_t(z)\,d\alpha_t, B_t(z)\,d\gamma_t]$ and integrating by parts, we get 
\begin{align*}
 -z^2 \int_0^t A_u(\bar w) A_u(z)\,d\alpha_u&=-z \int_0^t A_u(\bar w) \,dB_u(z)\\&=-z A_u(\bar w) B_u(z)\big|_0^t+ z \int_0^t B_u(z)dA_u(\bar w)
 \\&=-z A_u(\bar w) B_u(z)\big|_0^t+ 
 z\bar w \int_0^t B_u(z)B_u(\bar w)d\gamma_u.
\end{align*}
Obviously, we also have 
\begin{align*}
 -\bar w^2 \int_0^t A_u(\bar w) A_u(z)\,d\alpha_u=-\bar w A_u(z) B_u(\bar w)\big|_0^t+ 
 z\bar w \int_0^t B_u(z)B_u(\bar w)d\gamma_u 
\end{align*}
which follows from the preceding equation by interchanging the arguments $z$ and $\bar w$. It remains to subtract the first of these equations \eqref{eq:orilagp} from the second one.
The formulas \eqref{eq:orinorp} follow from \eqref{eq:orilagp}   by l'H\^opital's rule.
\end{prf}
Clearly,
the present chain  of de Branges functions is of exponential type and its diagonal structure function defines type $\tau_t$, $0\le t<\infty$,  so that  $\tau_t'=\alpha'_t\gamma'_t=1$. 
If $X$ is a fractional Brownian motion of Hurst index $0<H<1$, for instance,  the chain  is defined by \eqref{eq:Jentp} and $\alpha'_t=1/\gamma'_t=t^{1-2H}$.

\begin{rem}\label{rem:rem}
 {\rm
 For later use we make connection with the De Branges' theorem \cite[Section 37]{deB68} which proves relationship  \eqref{eq:difchp} between the even and odd components of  the de Branges functions and more. It regards the chain of matrices of entire functions which are real for a real $z$, called usually as
 {\it De Branges matrices} and denoted as
 \[M_t(z)=
  \begin{bmatrix}
   A_t(z) & B_t(z)\\
C_t(z) &  D_t(z)
  \end{bmatrix}
  \qquad 0\le t<\infty.
 \]
They are normalized so that $M_t(0)=I$. The chain of matrices is   
related to the structure function of
the chain of de Branges spaces \eqref{eq:famofspp} by $m_t= \dot M_t(0)\JJ$. Here
the chain of de Branges spaces \eqref{eq:famofspp} is generated by the upper row of the de Branges matrices and
$\dot M_t(z)$ denotes the derivative of $M(z)$ with respect to  variable $z$. 
De Branges \cite{deB68} Theorem 37 provides the following integral equation
\begin{equation}\label{eq:viap}
 M_r(z)\JJ-\JJ =z \int_0^rM_t(z)dm_t
\end{equation}
which in our particular case implies  \eqref{eq:difchp},  as well as
\fn{Parallel to \eqref{eq:stluzp}, we have 
\begin{align}
\label{eq:dym1.3+p}
 \frac{d^2C_t(z)}{d\alpha_t\,d\gamma_t} 
 =-z^2 C_t(z)\qquad 
 \frac{d^2D_t(z)}{d\gamma_t\,d\alpha_t} =-z^2 D_t(z).
\end{align}}
\begin{align}\label{eq:difch+p}
D_t(z)-1 =z\int_0^t C_u(z) d\alpha_u\qquad
-C_t(z)=z\int_0^t D_u(z) d\gamma_u. 
\end{align}
The Lagrange identities \eqref{eq:orilagp} turn now into
\begin{align}\label{eq:orilag+p}
\int_0^t D_u(\bar w) D_u(z)\,d\gamma_u= \frac{\bar w  D_t(z) C_t(\bar w)-zD_t(\bar w)C_t(z)}{z^2-{\bar w}^2} \nonumber\\
\int_0^t C_u(\bar w) C_u(z)\,d\alpha_u= \frac{  zC_t(\bar w)D_t(z)-\bar wC_t(z)D_t(\bar w)}{z^2-{\bar w}^2}
\end{align}
for all complex numbers $w$ and $z$.

In conclusion, we  point out that
by  rotation of the de Branges matrices 
\begin{align}
\label{eq:rot2p}
\JJ^\top M_t(z) \JJ=\begin{bmatrix} 
 D_t(z) & -C_t(z)\\
- B_t(z) &  A_t(z)\end{bmatrix}
\end{align}
the roles of the components $A_t(z)$ and $B_t(z)$ are given to  $D_t(z)$ and $-C_t(z)$. respectively. The latter even and odd entire functions  form  kernel 
\begin{equation}\label{eq:Qexlap}
 Q_t(w,z) =\frac{D_t(z)C_t(w)-C_t(z)D_t(w)}{\pi (z-\bar w)}
\end{equation}
which is introduced in \cite[Sections 27 and 28]{deB68} for characterizing bilinear $J$-forms of de Branges matrices. Applied to the present particular situation, the results of these sections  make available  %besides  chain $E_t(z), t>0$, one has 
a new chain of de Branges functions $\tilde E_t(z)=D_t(z)+iC(z), t>0$,  of the same exponential type as the former chain $E_t(z), t>0$. The reproducing kernel in this new space is \eqref{eq:Qexlap}.
\begin{comment}
we have
\begin{align*}
\text{diag} [\gamma_t,\alpha_t] 
=\begin{bmatrix} 
-\dot C_t(0) & -\dot D_t(0)  \\
\dot A_t(0) & \dot B_t(0)  \end{bmatrix} .
\end{align*}
In this connection recall Lemma~\ref{lem:bilJform} on bilinear $J$-forms of de Branges matrices. Presently, we are working with short de Branges spaces, i.e. $s(z)=1$ identically, and Corollary~\ref{cor:shortM} is in force saying that parallel to  chain $E_t(z), t>0$, one has the chain of de Branges functions $\tilde E_t(z), t>0$,  of the same exponential type as the former. 

Upon replacing the system of integral equations \eqref{eq:difchp} with  \eqref{eq:difch+p}, the statement of  Theorem~\ref{thm:convwedge}  will regard kernel \eqref{eq:Qexlap} (which takes the place of \eqref{eq:korneri11p}) and the Lagrange identities \eqref{eq:orilagp}  turn into
\eqref{eq:orilag+p}
for all complex numbers $w$ and $z$.

 Theorem~\ref{thm:thm37} which says that the structure function determines via equation \eqref{eq:mie} a chain of de Branges matrices $M_t(z)$ of  unit determinant, $0\le t\le r.$  The entries in the upper row are related to each other by \eqref{eq:difch}. Similarly, the entries in the lower row are related to each other by   
\eqref{eq:difch+p}
\end{comment}
}
\end{rem}

\subsection{Mercer's theorem}
\label{s:marcerip}

A key argument proving  the Karhunen--Lo\`{e}ve theorem is provided
by { Mercer's theorem}, as is said above.   
In the present section the latter theorem  is restricted to  symmetric positive definite kernels  
\[
 k: [a,b]\times [a,b]\rightarrow \RR,
\]
symmetry means that $k(x,y)=k(y,x)$ and positive definiteness that
\[
 \sum_{i,j} c_jc_k k(x_i,x_j)\ge 0
\]
for all finite number of points $x_j \in [a,b]$ and all choices of finite real numbers $c_j$.  With each such kernel, one associates the Hilbert-Schmidt operator 
$$
\KK\phi(t)\pd \int_a^b k(s,t)\varphi(s)ds
$$
on functions $\varphi\in L^2[a,b]$,
and defines the homogeneous Fredholm integral equation
\[
\int_a^b k(s,t)\varphi(s)ds= \epsilon \,\varphi(t),
\]
with squire integrable kernels 
\[
 \int_a^b  \int_a^b |k(s,t)|^2dsdt<\infty.
\]

Since $\KK$ is a linear operator, there is an orthonormal basis $(\varphi_n)$ in $L^2[a,b]$ where  $\varphi_n$ are eigenfunctions of the operator which correspond to non-negative eigenvalues $(\epsilon_n) $
that solve the integral equation
\[
\int_a^b k(s,t)\varphi_n(s)ds= \epsilon_n \,\varphi_n(t).
\]
The eigenvalues are assumed to be arranged in ascending order of magnitude. Mercer's theorem states an expansion of  kernel $k$
in terms of this basis.
%In the present text the kernels $K$ are represented by covariance functions of random processes. 
\begin{thm}\label{thm:mercerp}
Let $k$ be a continuous symmetric positive definite kernel. Its eigenfunctions $\varphi_n$ corresponding to positive eigenvalues $\epsilon_n$ are continuous on $[a,b]$ and  $k$ expands as
 \begin{equation}\label{eq:mercthmp}
  k(s,t)=\sum_{n=1}^\infty \epsilon_n\,\varphi_n(s)\varphi_n(t)
 \end{equation}
where the series converges absolutely and uniformly.
\end{thm}
\subsubsection{Kernels \eqref{eq:wedgianp} and \eqref{eq:wedgiamp}}

It is assumed in this and subsequent subsections that a chain of De Branges spaces $\HH(E_t)$, $t\ge 0,$ is given and satisfies the conditions of  Section~\ref{s:fm&map}. The objective is to formulate and prove  Mercer's type theorems for a number 
of positive definite kernels which are defined in terms of a diagonal structure function of the chain $m_t=\text{diag}[\alpha_t,\gamma_t]$, $t\ge 0.$ The next theorem and Theorem~\ref{thm:barczyp} regard the upper entry $\alpha$ of the structure function.  Theorem \ref{thm:wedgemergp} and Theorem~\ref{thm:barczygp}  regard the lower entry $\gamma$.
\begin{thm}\label{thm:wedgemercp}
Let a given chain of de Branges spaces $\HH(E_t), 0\le t\le r$, satisfy the conditions of Theorem~\ref{thm:convwedgep}.
\\
(i)\; If points on the real axes    $0<\lambda_1<\lambda_2<\cdots$ are such that $A_r(\lambda_n)=0$, then
the kernel 
 \begin{equation}
 \label{eq:wedgianp}
k(s,t)=
  \frac{
  \alpha_{s\wedge t}
  }
  {
  \sqrt{\alpha'_s}
  \sqrt{\alpha'_t}
  }  
 \end{equation}
is expanded in the series \eqref{eq:mercthmp} with $\epsilon_n=\lambda^{-2}_n$ and
\begin{equation*}
%\label{eq:varwegd}
 \varphi_n(t)=\frac{B_t(\lambda_n)}{\sqrt{\alpha'_t}}\sqrt{\frac{2}{\pi K_r(\lambda_n, \lambda_n)}}
\end{equation*}
where $\pi K_r(\lambda_n, \lambda_n)=-\dot A_r(\lambda_n)B_r(\lambda_n)$.\\
(ii)\, If points on the real axes    $0<\lambda_1<\lambda_2<\cdots$ are such that $B_r(\lambda_n)=0$, then
the kernel 
\begin{equation}
 \label{eq:wedgiamp}
k(s,t)=
  \frac{
  \alpha_{s\wedge t}-\dfrac{\alpha_s\alpha_t}{\alpha_r}
  }
  {
  \sqrt{\alpha'_s}
  \sqrt{\alpha'_t}
  }
\end{equation}
is expanded in the series \eqref{eq:mercthmp} with $\epsilon_n$ and 
$\varphi_n(t)$ of the same form as above but with $\pi K_r(\lambda_n, \lambda_n)=A_r(\lambda_n)\dot B_r(\lambda_n)$.
\end{thm}
In both items the expressions for $K_r(\lambda_n, \lambda_n)$ stem from %Lemma~\ref{lem:dym} which says that 
the following estimate of the reproducing kernel at the real diagonal 
\begin{align}\label{eq:lem:dymp}
\pi K_r(w, w)=A_r(w)\dot B_r(w)-\dot A_r(w)B_r(w), 
\end{align}
the derivatives are with respect to  variable $w$. If $w$ is not real, then  $K_r(w, w)>0$, but if is real, then the latter inequality holds true if, and only if, $ E_r(w)\ne 0$.

It follows that
%As we know from Theorem~\ref{thm:omzeros}
the zeros of $A$ and $B$ are always real and separate each other,  $A$ and $B$
cannot have common zeros, since $E$ has no zeros at the origin. Moreover, $A$
has no  zero at the origin and $B$ has a simple zero. Their non-zero zeros   are always real,  simple and interlacing.

Indeed,  the zeros are
necessarily real, because for non-real $w$ it holds that $K(w,w)
>0$, i.e. ${A(\bar w)B(w)-
    B(\bar w)A(w)}>0$ and neither $A$ nor $B$ can vanish.
If $E$ has no real
zeros, then $|E(x)|^2=A^2(x)+B^2(x)>0$ on the real line, 
so that $A$ and $B$ cannot have common zeros. It follows
from preceding that we again have
$K(x,x)
>0$, i.e. ${A(x)B'(x)-A'(x)}B(x)>0$ 
and only simple zeros are admitted. Moreover,
\begin{equation*}
%\label{eq:nomore}
\frac{d}{ dx}\frac{B(x)}{ A(x)}
=\frac{A(x)B'(x)- A'(x)B(x)}{A^2(x)}>
0
\end{equation*}
and the quotient $B/A$ is strictly increasing when finite. Hence
its zeros and poles are interlacing.
\fn{If $E$ does have real zeros, then these are zeros of $A$ and $B$ as
well, since $|E(x)|^2=A^2(\lambda)+B^2(\lambda)$ for real
$x$. As is noticed above, in this case 
inequality $K(w,w)
>0$ is not necessarily strict and only
conclusion we can draw is that zeros separate each other and
multiplicities can differ by at most  $1$. }

\begin{prf}
$(i)$\;In view of Mercer's Theorem \ref{thm:mercerp}, it is required to prove $\phi_n(t)$ are eigenfunctions of the present linear operator that correspond to eigenvalues  $ \epsilon_n$, i.e. the eigenfunctions and eigenvalues of  
the homogeneous Fredholm integral equation
\begin{equation}\label{eq:fredhp}
\int_0^r k(s,t)\phi(s)\,ds= \epsilon \,\phi(t)
\end{equation}
(assume $\epsilon\ne 0$ to exclude trivial solutions). The present kernel \eqref{eq:wedgianp}  turns this equation into
\[
\int_0^r \alpha_{s\wedge t}\,b_s\,d\gamma_s= \epsilon \,b_t, 
\]
where $\sqrt{\alpha_t'}\phi(t)=b_t$ and $d\gamma_s=ds/\alpha'_s$ by assumption. Upon integrating the left-hand side by parts 
\[
\int_0^r \alpha_{s\wedge t}\,b_s\,d\gamma_s=\int_0^t\alpha_{s}\,b_s\,d\gamma_s+\alpha_{t}\int_t^r \,b_s\,d\gamma_s= \int_0^t d\alpha_u\int_u^rb_sd\gamma_s
\]
we obtain
\begin{equation}\label{eq:mestlup}
\int_0^t d\alpha_u\int_u^rb_sd\gamma_s= \epsilon \,b_t. 
\end{equation}
Differentiating the latter integral equation first with respect to $d\alpha$ and then  $d\gamma$, we obtain   
\begin{equation}\label{eq:stlup}
 \frac{d^2 b_t}{d\gamma_td\alpha_t}=-\frac{1}\epsilon\, b_t,
\end{equation}
a differential equation 
of the Sturm--Liouville type. %, cf. Note 3 to Section~\ref{s:f}. 
The solutions to this equation are sought under the boundary conditions $b_0=0$ and $b'_r=0$.
At hand  of the differential equations \eqref{eq:stluzp} 
the solutions are available in terms of the underlying chain of de Branges functions $E_t(z)=A_t(z)-iB_t(z)$, $0\le t\le r$,  in the form  
$b_t=c B_t(z)$, with a suitable choice of argument $z$ and rescaling constant $c\ne 0$. Any choice will satisfy the boundary condition $b_0=0$, since $E_0(z)=1$ by assumption. The condition $b_r=1$ will be satisfied by suitable choice of rescaling.  

In order to determine  required roots of the differential equation \eqref{eq:stlup}, restrict  equation \eqref{eq:stluzp} to the real axis  
and select points $0<\lambda_1<\lambda_2<\cdots$ such that  $A_r(\lambda_n)=0$. With this choice, \eqref{eq:stlup} turns into
\begin{equation}\label{eq:stlunp}
 \frac{d^2 b_n(t)}{d\gamma_td\alpha_t}=-\frac{1}{\epsilon_n}\, b_n(t)
\end{equation}
where $1/\epsilon_n=\lambda_n^2$ and $b_n(t)=c_n B_t(\lambda_n)$, with non-zero constants $c_n$ such that $b_n(r)=1$. Hence $c_n=1/B_r(\lambda_n)$. It can be  verified that 
\[
 \int_0^r b_m(t)b_n(t)\,d\gamma_t= -\delta_{mn} \frac{c_n^2}{2} B_r(\lambda_n)\dot A_r(\lambda_n)=  -\frac{\delta_{mn}}{2} \frac{\dot A_r(\lambda_n)}{B_r(\lambda_n)}.
\]
Indeed, the second  equations in \eqref{eq:orilagp} and \eqref{eq:orinorp} show that since $\lambda_n$ are zeros of $A_r(z)$, the  elements  $B_t(\lambda_n)$ of $L^2([0,r],d\gamma)$ are mutually orthogonal, with the squire norm  
\[
 \|B(\lambda_n)\|^2_\gamma\pd\int_0^r |B_t(\lambda_n)|^2\,d\gamma_t= -\frac{1}{2} B_r(\lambda_n)\dot A_r(\lambda_n).
\]
Recall now  relation $\sqrt{\alpha_t'}\phi(t)=b_t$ between the solutions to the integral equation \eqref{eq:fredhp} and the differential equation \eqref{eq:stlup}. Since we are looking for a countable number of solutions and therefore make use of index $n$, the relation is to be rewritten as  $\sqrt{\alpha_t'}\phi_n(t)=b_n(t)$. The orthogonality property of the set of functions on the right-hand side of the latter identity, just verified, shows that the set of functions $\phi_n(t)$ are mutually orthogonal in $L^2[0,r]$, with the squire norm
\[
 \|\phi_n\|^2\pd\int_0^r |\phi_n(t)|^2\,dt=- \frac{c_n^2}{2} B_r(\lambda_n)\dot A_r(\lambda_n)= \frac{K_r(\lambda_n, \lambda_n)}{2B^2_r(\lambda_n)}= -\frac{\dot A_r(\lambda_n)}{2B_r(\lambda_n)}.
\]
Upon normalization $\varphi_n(t)=\phi_n(t) /\|\phi_n\|$, we get
\[
\int_0^r k(s,t)\varphi_n(s)ds= \lambda^{-2}_n \,\varphi_n(t).
\]
The proof is complete by verifying that the $\varphi_n(t)$ are of the required form, since
\[
 \varphi_n(t)=\frac{\phi_n(t)}{\|\phi_n(t)\|}=\frac{b_n(t)/\sqrt{\alpha'_t}}{\|\phi_n(t)\|}=\frac{c_n B_t(\lambda_n)}{\|\phi_n(t)\|\sqrt{\alpha'_t}}
\]
with $c_n=1/B_r(\lambda_n)$.\\
$(ii)$\; Differentiating the Fredholm equation \eqref{eq:fredhp} with kernel \eqref{eq:wedgiamp}, first with respect to $d\alpha$ and then $d\gamma$, leads again to \eqref{eq:stlup}.
If $0<\lambda_1<\lambda_2<\cdots$ are such that $B_r(\lambda_n)=0$, then equation~\eqref{eq:stlunp}  retains its form but with $\epsilon_n$ and $b_n(t)$  depending on these new points, i.e. $\epsilon_n=\lambda_n^{-2}$ and $b_n(t)=c_n\,B_t(\lambda_n)$. The boundary conditions to be satisfied are $b_n(0)=b_n(r)=0$. They are automatically satisfied in view of the present choice of the points $\lambda_n$. The rescaling constants $c_n\neq 0$ will be used for normalization. 
Invoke again the Lagrange identities \eqref{eq:orilagp} to obtain
\[
 \int_0^r b_m(t)b_n(t)\,d\gamma_t= \delta_{mn} \frac{c_n^2}{2} A_r(\lambda_n)\dot B_r(\lambda_n).
\]
Hence, in this case
\[
 \|\phi_n\|^2\pd\int_0^r |\phi_n(t)|^2\,dt=\frac{c_n^2}{2} A_r(\lambda_n)\dot B_r(\lambda_n)= \frac{c_n^2}{2}{K_r(\lambda_n, \lambda_n)}.
\]
Choosing $c_n^{-2}= {K_r(\lambda_n, \lambda_n)}/2$, we identify  $\varphi_n(t)=\phi_n(t)$ and write
\[
\int_0^r k(s,t)\varphi_n(s)ds= \lambda^{-2}_n \,\varphi_n(t).
\]
It is easily verified  that the $\varphi_n(t)$ are of the required form
\[
 \varphi_n(t)={b_n(t)/\sqrt{\alpha'_t}}={c_n B_t(\lambda_n)}/\sqrt{\alpha'_t}.
\]
The proof is complete.
\end{prf}
\noindent
The expansions of Theorem~\ref{thm:wedgemercp} can be rewritten entirely in terms of the reproducing kernel of the underlying chain of spaces which under the present conditions is given by \eqref{eq:korneri11p}. 
\begin{cor}\label{entinkp}
In the situation of Theorem~\ref{thm:wedgemercp}
the reproducing kernels $K_t(w,z)$, $0\le t\le r$, of
 the underlying chain of spaces 
satisfy
\begin{equation}\label{eq:Kwedge00p}
 K_{s \wedge t}(0,0)=
 \sum_{n} 
 \frac{
 K_s(0,\lambda_n)
 K_t(\lambda_n, 0)
 }
 {
 K_r(\lambda_n,\lambda_n)}
\end{equation}
where the $\lambda_n$ are either zeros of $A_r(z)$ or $B_r(z)$  symmetrically spread about the origin.
\end{cor}
\begin{prf}
%cf. \eqref{eq:korneri} and \eqref{eq:problem139}.
The form \eqref{eq:wedgianp} of the kernel $k(s,t)$ entails
\[
 \alpha_{s\wedge t}= 2\sum_{n=1}^\infty\frac{1}{\lambda^{2}_n} \frac{B_s(\lambda_n)B_t(\lambda_n)}{\pi K_r(\lambda_n, \lambda_n)}
\]
with zeros $\lambda_n$ of $A_r(z)$, and since 
$\pi K_t(0,0)=\alpha_t$ and $ \pi K_t(0,z)=B_t(z)/z$, this is equivalent to \eqref{eq:Kwedge00p}. To see that the latter holds true also with zeros  of $B_r(z)$, one only needs to notice that $0$ is one of the zeros and the term corresponding to this zero   equals to $\alpha_s\alpha_t/\alpha_r$.   
\end{prf}
\begin{comment}

Write $V_t=\hat{m}_{t}(0)$ and consider the kernel $V_{t\wedge s},
s,t\in [0,r]$, which is the covariance function for the
fundamental martingale $EM_sM_t=V_{t\wedge s}$. In 
order to prove
expansion (\ref{eq:vnaxo}), we look for the eigenfunction
expansion of the normalized kernel
\[
R(t,s)=\frac{V_{t\wedge s}}{\sqrt{V'_t}\sqrt{V'_s}}
\]
To this end, let us fix $\epsilon>0$ and solve the following integral
equation:
\[
\int_0^r R(t,s) \psi(s)ds=\epsilon\psi(t)
\]
subject to the boundary conditions $\psi(0)=\psi(1)=0$. 
\begin{align*}
 \int_0^t (V_s/\sqrt{V_s'})
 \psi_s ds+
 V_t
 \int_t^r(\psi_s/\sqrt{V_s'}) ds
 &=\int_0^t dV_u\int_u^r (\psi_s/\sqrt{V_s'})ds\\&
 =a\sqrt{V_t'}\psi_t
\end{align*}
Denote $\sqrt{V_t'}\psi_t=b_t$ to get
\[
 \int_0^t dV_u\int_u^r(b_s/V_s') ds=-\epsilon b_t
\]
Differentiate with respect to $d\alpha_t = V'_t dt$. We get 
\[
 \int_t^r b_s d\gamma_s=\epsilon \frac{db_t}{d\alpha_t}
\]
where  $d\gamma_t = dt/V'_t $. Hence the Fredholm equation reduces to
\[
 \frac{d^2 b_t}{d\gamma_t  d\alpha_t}=\frac{1}{\epsilon} b_t.
\]
Let $\lambda_n$ be zeros of $A_r(z)$. Take $b_t=B_t(\lambda_n)$.
As we know
\[
 \frac{d^2 B_t(\lambda_n)}{d\gamma_t  d\alpha_t}=-\lambda_n^2 B_t(\lambda_n).
\]
Then $\epsilon_n=\lambda_n^{-2}$ and $\phi_n(t)=B_t(\lambda_n)/\sqrt{\alpha'_t}$. But
\[
 \int_0^r B_t(w) B_t(z) d\alpha_t =\frac{w A_r(w)B_r(z)- z A_1(z)B_r(w)}{z^2-w^2}
\]
\end{comment}

\subsubsection{Kernel \eqref{eq:barczyp}}
%\label{sub:mthanex}

In the present subsection Theorem~\ref{thm:wedgemercp} will be generalized. Working under the conditions of the latter theorem, we extend the statement to the kernels of the form
\begin{equation}\label{eq:barczyp}
 \sqrt{\alpha'_s}\sqrt{\alpha'_t}\;k(s,t)=\alpha_{s\wedge t} - \frac{<1_s, \kappa>_\alpha <1_t, \kappa>_\alpha}{\|\kappa\|^2_\alpha}
\end{equation}
where $\kappa$'s are non-constant continuously differentiable kernels in $L^2([0,r],d\alpha)$ (as  $\kappa=0$ brings us  back to \eqref{eq:wedgianp} and a non-zero constant $\kappa$  to  \eqref{eq:wedgiamp}). In the sequel we make use  of the notation 
\begin{equation}\label{eq:noytp}
 g_t\pd <1_t, \kappa>_\alpha \qquad g'_t \pd \frac{dg_t}{d\alpha_t}
\end{equation}
just to simplify the exposition. On proving main Theorem ~\ref{thm:barczyp}, we shall follow the same succession of arguments as before. First we shall show that in the present case of kernel \eqref{eq:barczyp} the linear operator $\KK$ of Mercer's Theorem~\ref{thm:mercerp} is determined by a Fredholm equation which is reducible to the non-homogeneous second order differential equation \eqref{eq:mestmp}, 
subject to the boundary conditions 
\eqref{eq:mestmbp}. Then we shall show how the eigenvalues and eigenfunctions of the latter determine those of  operator $\KK$. 
To prepare the way for this, we begin with characterizing equation \eqref{eq:mestmp}.
\\~\\
The next theorem will provide necessary information. The first part of this theorem gives   
the {\it Lagrange type} identity \eqref{eq:kidevlagp} for two independent solutions. The second part will 
describe eigenvalues and eigenfunctions in terms of underline chain of de Branges spaces, the same as in Theorem~\ref{thm:wedgemercp}. 

\begin{thm} \label{thm:lagrangivitp}
Let $m_t=\diag[\alpha_t,\gamma_t]$, $t\in [0,r]$, be a diagonal structure function.
Given a twice differentiable function $g$ on interval $[0,r]$, define the non-homogeneous second order differential equation 
\begin{equation}
\label{eq:mestmp}
\frac{d^2b_t}{d\gamma_t\,d\alpha_t}=- w^2\Big(b_t +\frac{d^2g_t}{d\gamma_t\,d\alpha_t}\frac{<g, b>_\gamma}{\|g'\|^2_\alpha}\Big)
\end{equation}
subject to the boundary conditions 
\begin{equation}
\label{eq:mestmbp}
 b_0=0 \qquad \int_0^r b_t\,d\gamma_t=0.
\end{equation}
If $b_t(w)$ and $b_t(z)$ are two eigenfunctions of \eqref{eq:mestmp} corresponding to the two different {\rm (real)} eigenvalues $w$ and $z$,  then
\begin{equation}\label{eq:kidevlagp}
 \int_0^r
 b_t(w)\, b_t(z)\,d\gamma_t = 
 \frac{b_r(z)\,\beta_r(w)- b_r(w)\,\beta_r(z)}{z^2-w^2}
\end{equation}
where $ \beta_r(w)=b'_r(w)+w^2 g'_r \dfrac{<g, b(w)>_\gamma}{\|g'\|^2_\alpha}$ with $\displaystyle b'_r=\frac{db_t}{d\alpha_r}$ and $\displaystyle g'_r=\frac{dg_t}{d\alpha_r}$.

If $w=z$, then
\[
 \|b(w)\|^2_\gamma=\frac1{2w} b_r(w) \dot\beta_r(w)
\]
where $\dot\beta_r(w)=db_r(w)/dw$. 
\end{thm}

\begin{prf}
 It follows from \eqref{eq:mestmp} that
 \[
  b_t'(w) =-w^2 \Big(\int_0^t b_u(w)d\gamma_u +g'_t \dfrac{<g, b(w)>_\gamma}{\|g'\|^2_\alpha}\Big).
 \]
Therefore
\begin{align*}
-w^2 \int_0^t b_u(w)\, d\gamma_u=b'_t(w)+w^2 g'_t \dfrac{<g, b(w)>_\gamma}{\|g'\|^2_\alpha} =:  \beta_t(w)
\end{align*}
and
\begin{align*}
 -z^2\int_0^rb_t(w)\, b_t(z)\,d\gamma_t& =\int_0^rb_t(w)\,d  \beta_t(z)\\
 &= b_t(w)\, \beta_t(z)\Big|_0^r
 -\int_0^r\beta_t(z) b_t'(w) d\alpha_t
\end{align*}
by integrating by parts.
Invoking 
the boundary conditions, we obtain
\begin{align*}
 -z^2\int_0^rb_t(w)\, b_t(z)\,d\gamma_t= b_r(w)\, \beta_r(z)
 -\int_0^r  b'_t(z) b_t'(w) d\alpha_t.
\end{align*}
Indeed, identity $<\beta(z), b'(w)>_\alpha= <b'(z), b'(w)>_\alpha$ comes from 
the second of boundary conditions in this manner
\begin{align*}
<g', b'(w)>_\alpha=-w^2 \Big(\int_0^r g_t' d\alpha_t \int_0^t b_u(w) d\gamma_u +<g, b(w)>_\gamma\Big)\\
=-w^2 \Big(\int_0^r b_u(w) (g_r-g_u) d\gamma_u +<g, b(w)>_\gamma\Big)=0.
\end{align*}
Interchanging $w$ and $z$ we also obtain
\begin{align*}
 -w^2\int_0^rb_t(z)\, b_t(w)\,d\gamma_t =
 b_r(z)\, \beta_r(w) 
 +\int_0^r b'_t(w)   b'_t(z) d\alpha_t,
\end{align*}
so that \eqref{eq:kidevlagp} follows by taking the difference of these  two  equations.
The norm of $b_t$ in the metric of $L^2([0,r], d\gamma)$ is obtained from \eqref{eq:kidevlagp} 
by l'H\^opital's rule.
\end{prf}
By definition of $\beta_t(w)$, the second  of boundary conditions \eqref{eq:mestmbp} is equivalent to
\begin{equation}
\label{eq:wstwesp}
 \beta_r(w)=b'_r(w)+w^2 g'_r \dfrac{<g, b(w)>_\gamma}{\|g'\|^2_\alpha}=0.
\end{equation}
 
As is said, the second part of the preceding theorem will provide eigenvalues and eigenfunctions in terms of the chain of de Branges spaces $\HH(E_t), 0\le t\le r,$ the same as in the preceding section.  The density of its diagonal structure function $m$
is such that $\alpha'_t=1/\gamma'_t$. 
Moreover, it was pointed out
%But now we need more knowledge regarding the chain. We need to invoke 
%Theorem~\ref{thm:thm37} which says 
that the structure function determines 
a chain of de Branges matrices $M_t(z)$ of  unit determinant, $0\le t\le r,$ cf. equation \eqref{eq:viap}, whose entries in the upper and lower row are related to each other by \eqref{eq:difchp} and \eqref{eq:difch+p}, respectively.
%cf. \eqref{eq:4of0}, and parallel to \eqref{eq:dym1.3} we have \eqref{eq:dym1.3+p}.

The proof of the next statement makes use of the second equations
in \eqref{eq:stluzp} and \eqref{eq:dym1.3+p} regarding  the linearly independent entire functions $B_t(z)$ and $D_t(z)$, real for a real $z$. Due to these equations, both satisfy the homogeneous part of equation \eqref{eq:mestmp} 
%\begin{align*}
%\label{eq:dym1.3} \frac{d^2B_t(z)}{d\gamma_t\,d\alpha_t} =-z^2 B_t(z)\qquad\frac{d^2D_t(z)}{d\gamma_t\,d\alpha_t}  =-z^2 D_t(z)
%\end{align*}
\begin{comment}
subject to two different initial conditions $B_0(z)=0$, $\dot{B}$
the first under the boundary 
\begin{align}\label{eq:Mtebi}
 \frac{dD(w)}{d\alpha}=wC(w)\qquad \frac{dB(w)}{d\alpha}=wA(w)
\end{align}
\end{comment}
and  their linear combination  \eqref{eq:Gtwp}  satisfies  the second order inhomogeneous equation \eqref{eq:Ginhomp}. 

\begin{lem}\label{lem:Gtwp}
If $M_t(z), 0\le t \le r,$ is the chain of de Branges matrices associated as above with the underlying chain of de Branges spaces, then the linear combination 
\begin{equation}\label{eq:Gtwp}
G_t(z)=D_t(z) \int_0^t A_u(z) dg_u-B_t(z) \int_0^t C_u(z) dg_u,
\end{equation}
satisfies
\begin{equation}\label{eq:Ginhomp}
\frac{d^2G_t(z)}{d\gamma_td\alpha_t}=-z^2 G_t(z)+\frac{d^2g_t}{d\gamma_td\alpha_t}.
\end{equation}
\end{lem}
\begin{prf}
Since $\det M_t(z)=1$ and the derivatives of $B_t(z)$ and $D_t(z)$ with respect to $d\alpha_t$ are equal to $zA_t(z)$ and $zC_t(z)$, respectively, the derivative of $G_t(z)$ with respect to $d\alpha_t$ is equal to 
%Calculating $G'_t(w)=dG_t(w)/d\alpha_t$, we get
\begin{align}
\label{eq:Gprimtwp}
g'_t &+z \Big(C_t(z)\int_0^t A_u(z)dg_u - A_t(z) \int_0^t C_u(z)dg_u\,\Big)\nonumber\\&=g'_t -z \Big(C_t(z)\int_0^t g_udA_u(z) - A_t(z) \int_0^t g_u dC_u(z)\,\Big).
\end{align}
Differentiating the latter with respect to $d\gamma$, one  obtains
$\displaystyle\frac{d^2g_t}{d\gamma_td\alpha_t}$ minus $z$ times
\begin{align*}
\frac{dC_t(z)}{d\gamma_t}\!\!
\int_0^t g_udA_u(z) - \frac{dA_t(z)}{d\gamma_t} 
\!\!\int_0^t g_udC_u(z) -g_t 
\Big(
C_t(z)\frac{dA_t(z)}{d\gamma_t}- A_t(z)\frac{dC_t(z)}{d\gamma_t}
\Big).
\end{align*}
The last of these tree terms is equal to $z g_t \big(C_t(z)B_t(z)-A_t(z)D_t(z)\big)=-zg_t $, since the determinant of the de Branges matrix is equal to $1$. The sum of the other two terms is equal to $z$ times 
\begin{align*}
-&D_t(z)
\int_0^t g_udA_u(z) + B_t(z) 
\int_0^t g_udC_u(z)\\=&-D_t(z)\Big(g_t A_t(z)-\int_0^t A_u(z)dg_u\Big)
 + B_t(z)\Big(g_t C_t(z)-\int_0^t C_u(z)dg_u\Big)\\=&
 -g_t + G_t(z).
\end{align*}
Combining these results, we obtain
\eqref{eq:Ginhomp}. 
\end{prf}
The next statement makes use of the following
\begin{cor}
In the situation of the lemma
\begin{align*}
 -z&  <1, G(z)>_\gamma= C_r(z) \int_0^r A_t(z)dg_t-A_r(z) \int_0^r C_t(z)dg_t\\
 -z&  <g, G(z)>_\gamma= \int_0^r g_t dC_t(z) \int_0^t A_u(z)dg_u-\int_0^r g_t dA_t(z) \int_0^t C_u(z)dg_u.
\end{align*} 
\end{cor}
\begin{prf}
By \eqref{eq:Ginhomp} 
$$
-z^2 < 1, G(z)>_\gamma=<1, \frac{d^2G(z)}{d\gamma d\alpha}>_\gamma
-<1, \frac{d^2g}{d\gamma d\alpha}>_\gamma
= G'_r(z) -g'_r
$$
where $G'$ is the derivative of $G$ with respect to $d\alpha$. The first of the required formulas follows from \eqref{eq:Gprimtwp}. The second one follows from
\begin{align*}
-z <g, G(z)>_\gamma=\int_0^r g_t d\Big(C_t(z)\int_0^t A_u(z) dg_u -
A_t(z)\int_0^t C_u(z) dg_u\Big)\\
=-z \Big( D_t(z) \int_0^t A_u(z) dg_u-B_t(z) \int_0^t C_u(z) dg_u\Big)
\end{align*}
by definition \eqref{eq:Gtwp}.
\end{prf}

Now, we are in a position to complete Theorem~\ref{thm:lagrangivitp}.

\paragraph{Theorem~\ref{thm:lagrangivitp}}(continued). 
{\it Let $w_1, w_2,\cdots$ be positive roots of the determinantal equation
\begin{align}\label{eq:detertalp}
w^2 \begin{vmatrix}
 <1, G(w)>_\gamma & <g, G(w)>_\gamma\\<1, B(w)>_\gamma &<g, B(w)>_\gamma    
    \end{vmatrix}
=\|g'\|_\alpha^2 <1, B(w)>_\gamma   
\end{align}
arranged in ascending order of magnitude. They define  eigenvalues $\epsilon_n=1/w^2$  for  
equation \eqref{eq:mestmp},
%\begin{equation*}
%\label{eq:mestlus}
%\frac{d^2b_t}{d\gamma_t\,d\alpha_t}=-w^2\Big( b_t + \frac{d^2g_t}{d\gamma_t\,d\alpha_t}\frac{<g,b>_\gamma}{\|g'\|_\alpha^2}\Big)
%\end{equation*}
subject to the boundary conditions \eqref{eq:mestmbp}.
%Eigenvalues
%\begin{align*}
%<1, G(w)>_\gamma <g, B(w)>_\gamma-<1, B(w)>_\gamma <g, G(w)>_\gamma\\=w^{-2}\|g'\|_\alpha^2 <1, B(w)>_\gamma   
%\end{align*}
The corresponding eigenfunctions $b_t(w_n)$ are given by
\begin{equation}\label{eq:eigfp}
 b_t(w)=c(w) \Big( B_t(w)-\frac{<1, B(w)>_\gamma}{<1, G(w)>_\gamma} G_t(w)\Big)
\end{equation}
where $c$ is a normalization constant, making $b_t$ of unit norm in the metric of $L^2([0,r], d\gamma)$, i.e.
\[
 1/c^{2}(w)=\int_0^r \big( B_t(w)-\frac{<1, B(w)>_\gamma}{<1, G(w)>_\gamma} G_t(w)\big)^2d\gamma_t.
\]

}
\begin{prf}
The determinantal equation \eqref{eq:detertalp} is well-defined, since $B_t(w)$ and $D_t(w)$ are linearly independent, as well as $1$ and $g$ (by assumption made in the beginning of the present section that $g$ is not a constant). %Regarding the roots, see Note 1 to the present section. 
Note that $w=0$ is not a root, since $M_t(0)=I$ and $G_t(0)=g_t$ identically. Indeed, the right-hand side of \eqref{eq:detertalp}, divided through $w^2$, is a positive number, since $w^{-2}<1,B(w)>_\gamma$ tends to $\int_0^r\alpha_td\gamma_t$ as $w\searrow 0$, while the left-hand side vanishes.
The boundary conditions \eqref{eq:mestmbp} are easily verified.
Due to \eqref{eq:Ginhomp}, differentiating twice both sides of \eqref{eq:eigfp} yields 
 \begin{align*}
 \frac{d^2b_t(w)}{d\gamma_td\alpha_t}&=c(w) \Big( \frac{d^2G_t(w)}{d\gamma_td\alpha_t} -\frac{<1, B(w)>_\gamma}{<1, G(w)>_\gamma} \frac{d^2G_t(w)}{d\gamma_td\alpha_t}\Big)\\
 &=-w^2 b_t(w) -c(w)  \frac{<1, B(w)>_\gamma}{<1, G(w)>_\gamma}\frac{d^2g_t}{d\gamma_td\alpha_t}.
\end{align*}
Compare this with \eqref{eq:Ginhomp}. It is seen we need prove identity
\[
 w^2 \frac{<g, b>_\gamma}{\|g'\|_\alpha^2}=c(w) \frac{<1, B(w)>_\gamma}{<1, G(w)>_\gamma}
\]
whenever $w$ satisfies \eqref{eq:detertalp}.  
But this is easily obtained by calculating $<g, b>_\gamma$ from \eqref{eq:eigfp} in this manner  
\[
\frac{<g, b>_\gamma}{c(w)}=<g,B(w)>_\gamma -\frac{<1, B(w)>_\gamma}{<1, G(w)>_\gamma}<g,G(w)>        
       \]
and this, combined with \eqref{eq:detertalp}, proves the identity.       
\end{prf}
 
Theorem~\ref{thm:lagrangivitp} provides a key argument proving the next statement of Mercer's type.
 \begin{thm}\label{thm:barczygp}
The kernel of the current subsection that is  defined by \eqref{eq:barczyp} on $[0,r]\times [0,r]$ does expand
in  series \eqref{eq:mercthmp} with $\epsilon_n=1/w_n^2$, where $w_n$ are positive roots of the determinantal equation \eqref{eq:detertalp}, and with $\varphi_n(t)=b_t(w_n)/\sqrt{\alpha'_t}$, where $b_t(w)$ is given by \eqref{eq:eigfp}.  The series converges absolutely and uniformly.
 \end{thm}
\begin{prf}
 The linear operator $\KK$ of Mercer's Theorem~\ref{thm:mercerp} is now determined by means of
the  Fredholm integral equation  
\[
 \int_0^r {k(s,t)}\varphi_s ds= \epsilon \varphi_t
\]
with kernel 
\begin{equation*}
%\label{eq:barczy}
 \sqrt{\alpha'_s}\sqrt{\alpha'_t}\;k(s,t)=\alpha_{s\wedge t} - \frac{<1_s, g'>_\alpha <1_t, g'>_\alpha}{\|g'\|^2_\alpha},
\end{equation*}
where $g$ is so as in \eqref{eq:noytp}.
%Theorem~\ref{thm:lagrangivit}, cf. \eqref{eq:barczy}.
With the notation $b_t\pd \sqrt{\alpha'_t}\, \varphi_t$,  rewrite the latter equation as 
\[
 \int_0^r {\sqrt{\alpha'_s}\sqrt{\alpha'_t}}\,{k(s,t)}\,b_s \,d\gamma_s= \epsilon b_t.
\]
By making use of \eqref{eq:mestlup} %and the notations \eqref {eq:noyt} 
we obtain
\begin{align*}
%\label{eq:mestlus}
\int_0^t d\alpha_u\int_u^rb_sd\gamma_s
%&= \epsilon \,b_t+ \frac{<1_t, \kappa>_\alpha}{\|\kappa\|^2_\alpha}\int_0^r <1_s, \kappa>_\alpha \,b_sd\gamma_s\nonumber\\&
= \epsilon \,b_t+ \frac{g_t}{\|g'\|^2_\alpha}<g ,b>_\gamma.
\end{align*}
Set $t=0$ to verify 
the first of boundary conditions $b_0=0$ ($g_0=0$ by assumption). 
Differentiate both sides with respect to $d\alpha_t$
\begin{equation*}%\label{eq:socalled2cond}
%\label{eq:mestlus}
\int_t^rb_sd\gamma_s= \epsilon \,b'_t+  g'_t\frac{<g,b>_\gamma}{\|g'\|^2_\alpha}
\end{equation*}
where $b'_t={db_t}/{d\alpha_t}$ as before. 
%Set $t=r$ to verify that this  reduces to \eqref{eq:wstwes} if $\epsilon=w^{-2}$. 
Set $t=r$ to verify the second boundary condition \eqref{eq:wstwesp}.

Differentiate once more with respect to $d\gamma_t$, this time. We obtain
\begin{equation*}
%\label{eq:mestlus}
-b_t= \epsilon \,\frac{d^2b_t}{d\gamma_t\,d\alpha_t}+ \frac{ d^2g_t}{d\gamma_td\alpha_t}\frac{<g, b>_\gamma}{\|g'\|^2_\alpha} 
\end{equation*}
or
\begin{equation*}
%\label{eq:mestlus}
\frac{d^2b_t}{d\gamma_t\,d\alpha_t}=-w^2\Big( b_t + \frac{d^2g_t}{d\gamma_t\,d\alpha_t}\frac{<g,b>_\gamma}{\|g'\|^2_\alpha}\Big)
\end{equation*}
with $\epsilon=w^{-2}$. 
This second order inhomogeneous differential equation  is 
\eqref{eq:mestmp}, of course.
Therefore, Theorem \ref{thm:lagrangivitp} determines the eigenvalues and corresponding eigenfunctions for the present linear operator $\KK$ and our claim becomes a particular case of the general  Theorem~\ref{thm:mercerp}.
\end{prf}

In general, the entire functions \eqref{eq:Gtwp} depend  rather complicatedly on kernel $\kappa$ defining the positive definite kernel \eqref{eq:barczyp}. As is assumed at the beginning of this subsection, constant $\kappa$'s  are excluded and this excludes  the case $dg=d\alpha$.  
Matters are relatively simple  when $\kappa_t=\gamma_t$.

\begin{ex}\label{ex:kappap}
{\rm Consider a positive definite kernel on $[0,r] \times [0,r]$
 
 \begin{equation*}
 %\label{eq:barczy}
 \sqrt{\alpha'_s}\sqrt{\alpha'_t}\;k(s,t)=\alpha_{s\wedge t} - \frac{<1_s, \gamma>_\alpha <1_t, \gamma>_\alpha}{\|\gamma\|^2_\alpha}.
\end{equation*}
By definition \eqref{eq:noytp}, in this case $dg_t/d\alpha_t=\gamma_t$. Therefore the second term on the right-hand side of  \eqref{eq:Ginhomp} is independent of $t$ and equals to $1$.
Definition \eqref{eq:Gtwp} reduces to
\begin{equation}\label{eq:Ggammp}
 G_t(z)= \frac{1-D_t(z)}{z^2},
\end{equation}
which certainly satisfies the differential equation \eqref{eq:Ginhomp}. 
Indeed,   calculate
\begin{align*}
 G_t(z)&=D_t(z)\int_0^t \gamma_u A_u(z) d\alpha_u -B_t(z)\int_0^t \gamma_u C_u(z) d\alpha_u\\
 &=\frac{1}{z}\Big(D_t(z)\int_0^t \gamma_u dB_u(z) -B_t(z)\int_0^t \gamma_u d(D_u(z)-1) \Big)
\end{align*}
by integrating by parts 
\begin{align*}
 z & G_t(z)
 =D_t(z)\big(B_t(z)\gamma_t-\int_0^t B_u(z)d\gamma_u\big)\\& -B_t(z)\big((D_t(z)-1)\gamma_t -\int_0^t  (D_u(z)-1) d\gamma_u\big)\\
 &=D_t(z)\big(B_t(z)\gamma_t-\frac{1-A_t(z)}{z}\big)
 -B_t(z)\big((D_t(z)-1)\gamma_t +\frac{C_t(z)}{z}+\gamma_t\big).
\end{align*}
This proves \eqref{eq:Ggammp}. Since $g_t=\int_0^t\gamma_ud\alpha_u$ by definition, we have
\[
 g_t= \frac{D_t(z)-1}{z^2}\Big|_{z=0}=-G_t(0)
\]
and the determinantal equation \eqref{eq:detertalp} takes the form
\begin{align*}
%\label{eq:detertalg}
w^2 \begin{vmatrix}
 <1, G(w)>_\gamma & -<G(0), G(w)>_\gamma\\
 <1, B(w)>_\gamma & -<G(0), B(w)>_\gamma    
    \end{vmatrix}
=\|\gamma\|_\alpha^2 <1, B(w)>_\gamma.   
\end{align*}
In the upper row  $w^2 <1, G(w)>_\gamma=\gamma_r+C_r(w)/w$ 
and $w^2 <G(0), G(w)>_\gamma=\int_0^r(G_t(0) (1-D_t(w))d\gamma_t $.
 }
 $\hfill$ $\qed$
\end{ex}

\subsubsection{Kernels \eqref{eq:wedgiagp} and \eqref{eq:wedgiagip}}
In order to apply similar arguments also to the lower entry of the structure function, the roles of $\alpha$ and $\gamma$ are to be interchanged. This is  described with the help of rotations through angle $\pi/2$ about the origin so that
%carried out by the signature matrix \eqref{eq:mathsfJ} as 
\begin{align}
\label{eq:rot1p}
 \text{diag} [\gamma_t,\alpha_t] =\JJ^\top \text{diag} [\alpha_t, \gamma_t]\JJ.
\end{align}
%\cite[Theorem 37]{deB68} which  we have already referred to in 
As is pointed out in Remark~\ref{rem:rem}, %before equations \eqref{eq:difch+}, relates the structure function with the chain of de Branges matrices by $m_t= \dot M_t(0)\JJ$. Here  $\dot M_t(z)$ denotes the derivative of $M(z)$ with respect to the variable $z$. Since the similar rotations of the de Branges matrices yield
%\begin{align*}
%\label{eq:rot2p}
%\JJ^\top M_t(z) \JJ=\begin{bmatrix} 
% D_t(z) & -C_t(z)\\
%- B_t(z) &  A_t(z)\end{bmatrix},
%\end{align*}
% rotation  through angle $\pi/2$ about the origin
%\begin{align*}
%\text{diag} [\gamma_t,\alpha_t] 
%=\begin{bmatrix} 
%-\dot C_t(0) & -\dot D_t(0)  \\
%\dot A_t(0) & \dot B_t(0)  \end{bmatrix} .
%\end{align*}
such rotation \eqref{eq:rot2p} of the de Branges matrices transfers
the role of the entries in the upper row $A_t(z)$ and $B_t(z)$  to those of the lower row $D_t(z)$ and $-C_t(z)$.  The latter entire functions  form the chain of de Branges functions $\tilde E_t(z)=D_t(z)+iC_t(z), t>0,$ and the reproducing kernels 
\eqref{eq:Qexlap} in the corresponding spaces.
%In this connection recall Lemma~\ref{lem:bilJform} on bilinear $J$-forms of de Branges matrices. Presently, we are working with short de Branges spaces, i.e. $s(z)=1$ identically, and Corollary~\ref{cor:shortM} is in force saying that parallel to  chain $E_t(z), t>0$, one has the chain of de Branges functions $\tilde E_t(z), t>0$,  of the same exponential type as the former. 
With the system of integral equations \eqref{eq:difch+p} in stead of
\eqref{eq:difchp}, the statement of  Theorem~\ref{thm:convwedgep}  will regard kernel \eqref{eq:Qexlap},  in stead  of \eqref{eq:korneri11p}. The Lagrange identities \eqref{eq:orilagp} will turn into
\eqref{eq:orilag+p}
for all complex numbers $w$ and $z$.
\\~\\
\noindent
In terms of this new kernel \eqref{eq:Qexlap} one may rephrase    Corollary~\ref{entinkp} in this manner.\\
{\it In the situation of Theorem~\ref{thm:wedgemercp}
the  kernels $Q_t(w,z)$, $0\le t\le r$, associated with
 the underlying chain of spaces by \eqref{eq:Qexlap}, 
satisfy
\begin{equation}\label{eq:Kwedgegp}
 Q_{s \wedge t}(0,0)=
 \sum_{n} 
 \frac{
 Q_s(0,\lambda_n)
 Q_t(\lambda_n, 0)
 }
 {
 Q_r(\lambda_n,\lambda_n)}
\end{equation}
where the $\lambda_n$ are either zeros of $D_r(z)$ or $C_r(z)$  symmetrically spread about the origin.}\\~\\
This claim will be confirmed at hand of the next theorem. 
In both items of the theorem  kernel \eqref{eq:Qexlap} is evaluated at the real diagonal, say at a real point $w$, and  it holds that
\begin{equation}\label{eq:Wexla+p}
\pi Q_r(w, w)= \dot D_r(w)C_r(w)-D_r(w) \dot C_r(w)
\end{equation}
where the derivatives are with respect to variable $w$, like in  \eqref{eq:lem:dymp}. Since in the present situation the de Branges functions $\tilde E_t(z),$ $ t>0,$ have no zeros at the origin, the characterization of the zeros of $A$ and $B$ just after formula \eqref{eq:lem:dymp} applies to  
the zeros of $D$ and $C$ as well. The former has no zero at the origin and the latter has  a simple zero. Their non-zero zeros   are always real,  simple and interlacing.

\begin{thm}\label{thm:wedgemergp}
Let a given chain of de Branges spaces $\HH(E_t), 0\le t\le r$, satisfy the conditions of Theorem~\ref{thm:convwedgep}.
\\
(i)\; If points on the real axes    $0<\lambda_1<\lambda_2<\cdots$ are such that $D_r(\lambda_n)=0$, then
kernel 
 \begin{equation}
 \label{eq:wedgiagp}
k(s,t)=
  \frac{
  \gamma_{s\wedge t}
  }
  {
  \sqrt{\gamma'_s}
  \sqrt{\gamma'_t}
  }  
 \end{equation}
is expanded in  series \eqref{eq:mercthmp} with $\epsilon_n=\lambda^{-2}_n$ and
\begin{equation*}
%\label{eq:varwegd}
 \varphi_n(t)=\frac{C_t(\lambda_n)}{\sqrt{\gamma'_t}}\sqrt{\frac{2}{\pi Q_r(\lambda_n, \lambda_n)}}
\end{equation*}
where $\pi Q_r(\lambda_n, \lambda_n)= \dot D_r(\lambda_n)C_r(\lambda_n)$.\\
(ii)\, If points on the real axes    $0<\lambda_1<\lambda_2<\cdots$ are such that $C_r(\lambda_n)=0$, then
kernel 
\begin{equation}
 \label{eq:wedgiagip}
k(s,t)=
  \frac{
  \gamma_{s\wedge t}-\dfrac{\gamma_s\gamma_t}{\gamma_r}
  }
  {
  \sqrt{\gamma'_s}
  \sqrt{\gamma'_t}
  }
\end{equation}
is expanded in series \eqref{eq:mercthmp} with $\epsilon_n$ and 
$\varphi_n(t)$ of the same form as above but with $\pi Q_r(\lambda_n, \lambda_n)=-D_r(\lambda_n)\dot C_r(\lambda_n)$.
\end{thm}
In fact, this is a reformulation of Theorem~\ref{thm:wedgemercp}
in accord with the rotations \eqref{eq:rot2p} and \eqref{eq:rot1p}. 
The proof requires  same 
adaptations and is not difficult to carry out. %For completeness, however, we shall follow the course of changes.

\begin{prf}
$(i)$\; 
Upon differentiating first with respect to $d\gamma$ and then $d\alpha$, the homogeneous Fredholm integral equation
\eqref{eq:fredhp}
with kernel \eqref{eq:wedgiagp} turns into
the differential equation 
of the Sturm--Liouville type
\begin{equation}\label{eq:stlu+p}
 \frac{d^2 b_t}{d\alpha_td\gamma_t}=-\frac{1}\epsilon\, b_t,
\end{equation}
where  $b_t$ is rescaling $\sqrt{\gamma_t'}\phi(t)=b_t$ of   function $\phi(t)$ involved in the Fredholm equation. 
The solutions to this equation are sought under the boundary conditions $b_0=0$ and $b'_r=0$.
The first of equations \eqref{eq:stluzp} offers  solution $A_t(z)$, but this fails to satisfy the boundary conditions. We shall see in a moment that another possibility, the first of equations  \eqref{eq:dym1.3+p}, gives a required solution. 
Looking  for
 a solution to \eqref{eq:stlu+p}  in the form 
$b_t=c\, C_t(z)$, with $\epsilon =1/z^2$, we shall choose  suitably argument $z$ and the rescaling constant $c\ne 0$.  Any choice will satisfy the boundary condition $b_0=0$, since $\tilde E_0(z)=1$ by assumption. Another condition $b_r=1$ is to be guaranteed by rescaling. 
In the first of equations  \eqref{eq:dym1.3+p} restrict argument $z$  to real axis   
and select points $0<\lambda_1<\lambda_2<\cdots$ such that  $D_r(\lambda_n)=0$. With this choice, \eqref{eq:stlu+p} turns into
\begin{equation}\label{eq:stlunCp}
 \frac{d^2 b_n(t)}{d\alpha_td\gamma_t}=-\frac{1}{\epsilon_n}\, b_n(t)
\end{equation}
where $1/\epsilon_n=\lambda_n^2$ and $b_n(t)=c_n C_t(\lambda_n)$, with non-zero constants $c_n$ such that $b_n(r)=1$. Hence $c_n=1/C_r(\lambda_n)$. It follows from \eqref{eq:orilag+p} that
\[
 \int_0^r b_m(t)b_n(t)\,d\alpha_t=\delta_{mn} \frac{c_n^2}{2} C_r(\lambda_n)\dot D_r(\lambda_n)= \frac{\delta_{mn}}{2} \frac{\dot D_r(\lambda_n)}{C_r(\lambda_n)}.
\]
%Indeed, the second  equations in \eqref{eq:orilag} and \eqref{eq:orinor} show that since $\lambda_n$ are zeros of $D_r(z)$, the  elements  $C_t(\lambda_n)$ of $L^2([0,r],d\alpha)$ are mutually orthogonal, with the squire norm  
%\[
 %\|C(\lambda_n)\|^2_\alpha\pd\int_0^r |C_t(\lambda_n)|^2\,d\alpha_t= -\frac{1}{2} C_r(\lambda_n)\dot D_r(\lambda_n).
%\]
Recall now relation $\sqrt{\gamma_t'}\phi(t)=b_t$ between the solutions to the integral equation \eqref{eq:fredhp} and the differential equation \eqref{eq:stlu+p}. Since we are looking for a countable number of solutions and therefore make use of index $n$, the relation is to be rewritten as  $\sqrt{\gamma_t'}\phi_n(t)=b_n(t)$. The orthogonality property of the set of functions on the right-hand side of this identity (seen from the latter display) shows that the set of functions $\phi_n(t)$ are mutually orthogonal in $L^2[0,r]$, with the squire norm
\[
 \|\phi_n\|^2=\int_0^r |\phi_n(t)|^2dt=\int_0^r |b_n(t)|^2\,d\alpha_t= \frac{\dot D_r(\lambda_n)}{2C_r(\lambda_n)}.
\]
Upon normalization $\varphi_n(t)=\phi_n(t) /\|\phi_n\|$, we get
\[
\int_0^r k(s,t)\varphi_n(s)ds= \lambda^{-2}_n \,\varphi_n(t).
\]
The proof is complete, since the $\varphi_n(t)$ are of the required form
\[
 \varphi_n(t)=\frac{\phi_n(t)}{\|\phi_n(t)\|}=\frac{b_n(t)/\sqrt{\gamma'_t}}{\|\phi_n(t)\|}=\frac{ C_t(\lambda_n)/\sqrt{\gamma'_t}}{C_r(\lambda_n)}\sqrt{\frac{2C_r(\lambda_n)}{ \dot D_r(\lambda_n)}}.
\]
\noindent
$(ii)$\; Differentiating the Fredholm equation \eqref{eq:fredhp} with kernel \eqref{eq:wedgiagip}, first with respect to $d\alpha$ and then $d\gamma$, leads again to \eqref{eq:stlu+p}.
If $0<\lambda_1<\lambda_2<\cdots$ are such that $C_r(\lambda_n)=0$, then equation~\eqref{eq:stlunCp}  retains its form  but with $\epsilon_n$ and $b_n(t)$  depending on these new points, i.e. $\epsilon_n=\lambda_n^{-2}$ and $b_n(t)=c_n\,C_t(\lambda_n)$. The boundary conditions  $b_n(0)=b_n(r)=0$ are  automatically satisfied. The rescaling constants $c_n\neq 0$ will be used for normalization. 
Invoke again the Lagrange identities \eqref{eq:orilag+p} to obtain
\[
 \int_0^r b_m(t)b_n(t)\,d\alpha_t=- \delta_{mn} \frac{c_n^2}{2} D_r(\lambda_n)\dot C_r(\lambda_n).
\]
Hence, in this case
\[
 \|\phi_n\|^2\pd\int_0^r |\phi_n(t)|^2\,dt=-\frac{c_n^2}{2} D_r(\lambda_n)\dot C_r(\lambda_n)= \frac{c_n^2}{2}{Q_r(\lambda_n, \lambda_n)}.
\]
Choosing $c_n^{-2}= {Q_r(\lambda_n, \lambda_n)}/2$, we identify  $\varphi_n(t)=\phi_n(t)$ and write
\[
\int_0^r k(s,t)\varphi_n(s)ds= \lambda^{-2}_n \,\varphi_n(t).
\]
It is easily verified  that the $\varphi_n(t)$ are of the required form
\[
 \varphi_n(t)={b_n(t)/\sqrt{\gamma'_t}}={c_n C_t(\lambda_n)}/\sqrt{\gamma'_t}.
\]
The proof is complete.
\end{prf}

\subsubsection{Kernel \eqref{eq:barczygp}}
As in the preceding subsection, the idea is to substitute in \eqref{eq:barczygp} the upper entry of the structure function with the lower entry. So, Mercer's Theorem~\ref{thm:barczyp} proved below will regard the kernel

\begin{equation}\label{eq:barczygp}
 \sqrt{\gamma'_s}\sqrt{\gamma'_t}\;k(s,t)=\gamma_{s\wedge t} - \frac{<1_s, \kappa>_\gamma <1_t, \kappa>_\gamma}{\|\kappa\|^2_\gamma}
\end{equation}
where $\kappa$'s are now non-constant continuously differentiable kernels in $L^2([0,r],d\gamma)$. Accordingly, the notations \eqref{eq:noytp} will turn into 
\begin{equation}\label{eq:noyt+p}
 g_t\pd <1_t, \kappa>_\gamma \qquad g'_t \pd \frac{dg_t}{d\gamma_t}.
\end{equation}
Moreover, in the forthcoming proof all differentials with respect to $d\gamma$ will receive the similar abbreviation, say $db_t/d\gamma_t=b_t'$ and so forth. 

 \begin{thm}\label{thm:barczyp}
The kernel of the current subsection defined by  \eqref{eq:barczygp} on  $[0,r]\times [0,r]$ does  expand
in the absolutely and uniformly convergent series \eqref{eq:mercthmp} with $\epsilon_n$ and $\varphi_n(t)$ defined as follows:\\
(a)\; $\epsilon_n=1/w_n^2$, where the $w_n$ are positive roots of the determinantal equation
\begin{align}\label{eq:detertal+p}
w^2 \begin{vmatrix}
 <1, \tilde G(w)>_\alpha & <g, \tilde G(w)>_\alpha\\<1, C(w)>_\alpha &<g, C(w)>_\alpha    
    \end{vmatrix}
=\|g'\|_\gamma^2 <1, C(w)>_\alpha 
\end{align}
arranged in ascending order of magnitude and the entire functions $\tilde G_t(z)$ are defined by
\begin{equation}\label{eq:Gtw+p}
\tilde G_t(z)=A_t(z) \int_0^t D_u(z) dg_u-C_t(z) \int_0^t B_u(z) dg_u.
\end{equation}
(b)\;
$\varphi_n(t)=b_t(w_n)/\sqrt{\alpha'_t}$, where $b_t(w)$ is given by 
\begin{equation}\label{eq:eigf+p}
 b_t(w)=c(w) \Big( C_t(w)-\frac{<1, C(w)>_\alpha}{<1, \tilde G(w)>_\alpha} \tilde G_t(w)\Big)
\end{equation}
with a normalization constant $c$, making $b_t$ of unit norm in the metric of $L^2([0,r], d\alpha)$, i.e.
\[
 1/c^{2}(w)=\int_0^r \big( C_t(w)-\frac{<1, C(w)>_\alpha}{<1, \tilde G(w)>_\alpha} \tilde G_t(w)\big)^2d\alpha_t.
\]
 \end{thm}
\begin{prf}
The entire functions 
\eqref{eq:Gtw+p} are derived from \eqref{eq:Gtwp} by the replacements $A_t(z)\leftrightarrow D_t(z)$ and $B_t(z)\leftrightarrow -C_t(z)$.   
Lemma~\ref{lem:Gtwp} suggests 
\begin{equation}\label{eq:Ginhom+p}
\frac{d^2\tilde G_t(z)}{d\alpha_td\gamma_t}=-z^2 \tilde G_t(z)+\frac{d^2g_t}{d\alpha_td\gamma_t},
\end{equation}
cf. \eqref{eq:Ginhomp}. To confirm this, take the derivative of $\tilde G_t(z)$ with respect to $d\gamma$ to get 
\begin{align*}
g'_t -z \Big(B_t(z)\int_0^t g_udD_u(z) - D_t(z) \int_0^t g_u dB_u(z)\,\Big)
\end{align*}
which coincides with \eqref{eq:Gprimtwp} upon  the required replacements.
Differentiating the latter with respect to $d\alpha$, one  obtains
$\displaystyle\frac{d^2g_t}{d\alpha_td\gamma_t}$ minus $z$ times
\begin{align*}
\frac{dB_t(z)}{d\alpha_t}
\!\!\int_0^t\!\! g_udD_u(z) - \frac{dD_t(z)}{d\alpha_t} 
\!\!\int_0^t\!\! g_udB_u(z) -g_t 
\Big(
B_t(z)\frac{dD_t(z)}{d\alpha_t}- D_t(z)\frac{dB_t(z)}{d\alpha_t}
\Big).
\end{align*}
The last of these tree terms is equal to $z g_t \big(C_t(z)B_t(z)-A_t(z)D_t(z)\big)=-zg_t $, since the determinant of the de Branges matrix is equal to $1$. The sum of the other two terms is equal to $z$ times 
\begin{align*}
-&A_t(z)
\int_0^t g_udD_u(z) + C_t(z) 
\int_0^t g_udB_u(z)\\=&-A_t(z)\Big(g_t D_t(z)-\int_0^t D_u(z)dg_u\Big)
 + C_t(z)\Big(g_t B_t(z)-\int_0^t B_u(z)dg_u\Big)\\=&
 -g_t + \tilde G_t(z).
\end{align*}
So, 
\eqref{eq:Ginhom+p} is confirmed. It provides a key argument for the remainder of the proof. 

The linear operator $\KK$ of Mercer's Theorem~\ref{thm:mercerp} is now determined by means of
\begin{comment}
{\sf
 \begin{align*}
  \frac{db_r}{d\alpha_r}&=\Delta_r(w) \frac{dB_r(w)}{d\alpha_r}\\&+
  \Big( \frac{dD_r(w)}{d\alpha_r}\int_0^r A_u(w)g'_ud\alpha_u -
 \frac{dB_r(w)}{d\alpha_r}\int_0^r C_u(w)g'_ud\alpha_u \Big) \frac{D_r(w)}{w}+g'_r \frac{D_r(w)}{w}\\
 &=w\Delta_r(w) A_r(w)+ 
  \Big( C_r(w)\int_0^r g'_udB_u(w) -
 A_r(w)\int_0^r  g'_udD_u(w) \Big)\frac{D_r(w)}{w}+g'_r \frac{D_r(w)}{w}
 \end{align*}
$g'=dg/d\alpha$.

Boundary condition 
\begin{equation*}
%\label{eq:mestlus}
\frac{db_r}{d\alpha_r}= -w^2 \frac{ g'_r}{\|g'\|^2_\alpha}\int_0^r g_u\,b_ud\gamma_u.
\end{equation*}
Hence
\begin{align*}
w^2\Delta_r(w) A_r(w)+ 
  \Big( C_r(w)\int_0^r g'_udB_u(w) -
 A_r(w)\int_0^r  g'_udD_u(w) \Big)\,{D_r(w)}=0
 \end{align*}
 gives $w$'s.
}
\end{comment}
the  Fredholm integral equation  
\[
 \int_0^r {k(s,t)}\varphi_s ds= \epsilon \varphi_t
\]
with kernel \eqref{eq:barczygp} which, rewritten in terms of \eqref{eq:noyt+p}, takes the form  
\begin{equation*}
%\label{eq:barczy}
 \sqrt{\gamma'_s}\sqrt{\gamma'_t}\;k(s,t)=\gamma_{s\wedge t} - \frac{<1_s, g'>_\gamma <1_t, g'>_\gamma}{\|g'\|^2_\gamma}.
\end{equation*}
With the notation $b_t\pd \sqrt{\gamma'_t}\, \varphi_t$,  rewrite the Fredholm equation as 
\[
 \int_0^r {\sqrt{\gamma'_s}\sqrt{\gamma'_t}}\,{k(s,t)}\,b_s \,d\alpha_s= \epsilon b_t.
\]
By making use of \eqref{eq:mestlup} %and the notations \eqref {eq:noyt} 
we obtain
\begin{align*}\label{eq:mestlus+p}
\int_0^t d\gamma_u\int_u^rb_sd\alpha_s
%&= \epsilon \,b_t+ \frac{<1_t, \kappa>_\alpha}{\|\kappa\|^2_\alpha}\int_0^r <1_s, \kappa>_\alpha \,b_sd\gamma_s\nonumber\\&
= \epsilon \,b_t+ \frac{g_t}{\|g'\|^2_\gamma}<g ,b>_\alpha.
\end{align*}
Set $t=0$ to verify 
the first of boundary conditions $b_0=0$ ($g_0=0$ by assumption). 
Differentiate both sides with respect to $d\gamma$
\begin{equation*}%\label{eq:socalled2cond+}
\int_t^rb_sd\alpha_s= \epsilon \,b'_t+  g'_t\frac{<g,b>_\alpha}{\|g'\|^2_\gamma}
\end{equation*}
where $b'_t={db_t}/{d\gamma_t}$ as before. 
%Set $t=r$ to verify that this  reduces to \eqref{eq:wstwes} if $\epsilon=w^{-2}$. 
Set $t=r$ to verify the second boundary condition.

Differentiate once more with respect to $d\alpha$, this time. We obtain
\begin{equation*}
%\label{eq:mestlus}
-b_t= \epsilon \,\frac{d^2b_t}{d\alpha_t\,d\gamma_t}+ \frac{ d^2g_t}{d\alpha_td\gamma_t}\frac{<g, b>_\gamma}{\|g'\|^2_\alpha} 
\end{equation*}
or
\begin{equation*}
%\label{eq:mestlus}
\frac{d^2b_t}{d\alpha_t\,d\gamma_t}=-w^2\Big( b_t + \frac{d^2g_t}{d\alpha_t\,d\gamma_t}\frac{<g,b>_\alpha}{\|g'\|^2_\gamma}\Big)
\end{equation*}
with $\epsilon=w^{-2}$. 
This second order inhomogeneous differential equation  is 
\eqref{eq:mestmp} with $\alpha$ and $\gamma$ interchanged. Therefore, an adjusted version of 
Theorem \ref{thm:lagrangivitp} determines the eigenvalues and corresponding eigenfunctions for the present linear operator $\KK$.
Indeed, an adaptation entails substituting \eqref{eq:Ginhomp} with
\eqref{eq:Ginhom+p} and two equations \eqref{eq:detertalp} and \eqref{eq:eigfp} determining the eigenvalues and  eigenfunctions
with \eqref{eq:detertal+p} and \eqref{eq:eigf+p}. The proof is complete.
\end{prf}

\subsection{ Fundamental
martingales}\label{s:chemikarhunenip}
%(a)\; For a given spectral measure $\mu$ of property \eqref{eq:sqi}there exists a chain  of de Branges spaces $\HH(E_t)$, $t\ge 0$, of exponential type which is contained isometrically in $L^2(\mu)$. We know this from Theorem~\ref{thm:problem137+}. \\~\\
The preceding results are intended to 
derive the Karhunen-Lo\`{e}ve expansions for certain Gaussian processes associated with double-sided si-processes. A double-sided si-process 
$X$ is defined by the spectral representation \eqref{eq:anamova10p} of its covariance function, where the arguments $s$ and $t$ are allowed to take any value on a whole real axis. Therefore process $X$  can be represented by its even and odd parts
$ X_t= \text{sign}\, (t)\, X^e_{|t|}+ X^o_{|t|}.$
%\end{equation*}
These even and odd processes emerge from the origin and develop to the right,   having  moving average representations in the sense of Theorem~\ref{thm:mavoorxypp}. The Wiener integrals \eqref{eq:mavoorxpp} of this theorem is now simplified to   
\begin{equation*}
%\label{eq:mavoor11}
X^e_{t}=\int_0^t \varphi dM^e\qquad 
X^o_{t}=\int_0^t \psi dM^o,
\end{equation*}
$t\ge 0$, since in Theorem~\ref{thm:convwedgep} and thereafter a chain of de Branges spaces of exponential type  which  {\it determines} the given si-process has a structure function such that $\alpha_t'=1/\gamma_t'$, which means  the corresponding type is   $\tau_t=t$. 
The even and odd fundamental martingales $M^e$ and $M^o$ that appear in the latter display are defined 
with the help of the spectral isometry \eqref{eq:alsnogp} -- they are mutually independent Gaussian processes defined on the same probability space as $X$,  generating the same filtration as $X^e$ and $X^o$, respectively,  with the second order properties 
\begin{align}\label{eq:dubelvar11p}
\EEE (M^e_sM^e_t)&=\int \frac{B_s(\lambda)B_t(\lambda)}{\lambda^2}\mu(d\lambda) =\pi \alpha_{s\wedge t} \nonumber\\
\EEE (M^o_sM^o_t)&=\int \frac{\big(A_s(\lambda)-1\big)\big(A_t(\lambda)-1\big)}{\lambda^2}\mu(d\lambda) =\pi \gamma_{s\wedge t},
\end{align}
cf. \eqref{eq:dubelvarp}. 
The present section obtains KL-series expansions for processes allied to both $M^e$ and $M^o$,  first to the even and then to the odd   one. 
\subsubsection{ Even 
martingales}
In Section~\ref{sub:asspp} we apply the sampling formulas 
%of the de Branges theory, cf.  Theorem~\ref{thm:roparco}, 
to obtain the generalized PW-series for both fundamental martingales. The result is formulated as Theorem~\ref{thm:wcmsp}, item $(i)$.
The expansion for the even fundamental martingale of  this theorem is in a direct relation with the KL-expansion of the following statement. %consequence of Theorem~\ref{thm:wedgemerc}, item $(i)$.

\begin{thm}\label{thm:evenisKLp}
Let $X$ be a double-sided si-process  determined by a chain of de Branges spaces of exponential type so as in Theorem~\ref{thm:convwedgep}. Let $M^e$ be the fundamental martingale which defines the moving average representation of the even part of $X$, with the second order moments
%$\EEE M^e_s\,M^e_t=\pi \alpha_{s\wedge t}$,cf.
\eqref{eq:dubelvar11p}.
Then process $M^e_t/\sqrt{\alpha'_t},$  $0\le t\le r,$   decomposes in  KL-series \eqref{eq:Loevep}   
\[
\frac{1}{\sqrt{\alpha'_t}}M^e_t=\sum_{n=1}^\infty \omega_n \varphi_n(t)\,\xi_n 
\]
where  the $\xi_n$ are i.i.d. Gaussian random variables, $\{\varphi_n(t)= B_t(\lambda_n)/\sqrt{\alpha'_t}\}$ is an orthonormal basis in $L^2[0,r]$, and 
\[
\omega_n=\frac{1}{\lambda_n} \sqrt{\frac{2}{\pi K_r(\lambda_n,\lambda_n)}}
\]
with the positive zeros $\lambda_n$ of $A_r(z)$.

\end{thm}
\begin{prf}
In  virtue of Theorem~\ref{thm:Loevep}, it suffices to verify  the covariance function of  process  $M^e_t/\sqrt{\alpha'_t}$ is of the form \eqref{eq:wedgianp}. But this is clear by the second order property \eqref{eq:dubelvarp} of the even fundamental martingale.
\end{prf}
\noindent
The theorem just proved generalizes the classical  Karhunen--Lo\`{e}ve  expansion \eqref{eq:V.2.17p}. Due to  Theorem~\ref{thm:wedgemercp}, item $(ii)$, it is equally easy to generalize also expansion  \eqref{eq:V.2.18p} of Brownian bridge.

\begin{thm}\label{thm:evenibrKLp}
In the situation of Theorem~\ref{thm:evenisKLp}, the martingale bridge
\begin{equation}\label{eq:evenibrtp}
 B^e_t=M^e_t-\frac{\alpha_t}{\alpha_r}M_r^e,
\end{equation}  
$t\in [0,r]$, divided through $\sqrt{\alpha'_t}$,  
decomposes in  KL-series \eqref{eq:Loevep}  
\[
\frac{1}{\sqrt{\alpha'_t}}B^e_t=\sum_{n=1}^\infty \omega_n \varphi_n(t)\,\xi_n 
\]
where  the $\xi_n$ are i.i.d. Gaussian random variables, $\omega_n$ and $\{\varphi_n\}$ are  of the same form as in Theorem~\ref{thm:evenisKLp} but now $\lambda_n$ are 
the positive zeros  of $B_t(z)$.

\end{thm}
\begin{prf}
By Theorem~\ref{thm:Loevep},  the covariance function of the process in question % $M^e_t/\sqrt{\alpha'_t}$ 
is of the form \eqref{eq:wedgiamp}. 
%But this is clear by the second order property \eqref{eq:dubelvar} of the even fundamental martingale.
\end{prf}
This is reproving of the statement made just after Theorem~\ref{thm:wcmsp} that regards PW-series for  even martingale bridges.\\~\\ 
\noindent 
The martingale bridge  is an even fundamental martingale $M^e$ centered by its conditional expectation, given its value at the right endpoint $r$. Section~\ref{s:marcerip} contains  the subsection on kernel \eqref{eq:barczyp}, the results of which  make available an extension of Theorem~\ref{thm:evenibrKLp} to an even fundamental martingale $M^e$ centered by its conditional expectation, given
the value of a martingale $$N^e_t=\int_0^t \kappa_u dM^e_u$$ at the right endpoint $r$. The kernel $\kappa$ is assumed to be non-constant, otherwise  we would turn back to an already handled  martingale bridge. The Gaussian process so defined
\begin{equation}\label{eq:extmbp}
B^e_t=M^e_t-\EEE (M^e_t| N^e_r)  
\end{equation}
$t\in [0,r]$, has the covariance function
\[
 r(s,t)=\pi \Big(\alpha_{s\wedge t} - \frac{<1_s, \kappa>_\alpha <1_t, \kappa>_\alpha}{\|\kappa\|^2_\alpha}\Big).
 \]
In Section~\ref{s:marcerip} we have introduced the notations \eqref{eq:noytp} and assumed smoothness of the kernel $\kappa$ as to satisfy Theorem~\ref{thm:lagrangivitp}. This theorem provides a key argument for proving Mercer's type Theorem~\ref{thm:barczygp} which implies

\begin{thm}\label{thm:barczy1p}
In the situation of Theorem~\ref{thm:barczygp} process
\eqref{eq:extmbp}, divided through $\sqrt{\alpha'_t}$,  
decomposes in  KL-series \eqref{eq:Loevep}  
\[
\frac{1}{\sqrt{\alpha'_t}}B^e_t=\sum_{n=1}^\infty \omega_n \varphi_n(t)\,\xi_n 
\]
where  the $\xi_n$ are i.i.d. Gaussian random variables,  $\{\varphi_n(t)= b_t(w_n)/\sqrt{\alpha'_t}\}$ is an orthonormal basis in $L^2[0,r]$, and 
\[
\omega_n=\frac{1}{w_n c(w_n)}. 
\]
Here $b_t(w)$ is given by \eqref{eq:eigfp}, the $w_n$
are  positive zeros  of the determinantal equation \eqref{eq:detertalp} and $c(w)$ is a normalization constant in \eqref{eq:eigfp}.
\end{thm}
\begin{prf}
The claim that Theorem~\ref{thm:barczygp} implies the present statement follows from the comparison of the covariance function $r(s,t)$ displayed above with kernel \eqref{eq:barczyp}.  
\end{prf}
\subsubsection{ Odd
martingales}
Preceding results are adapted to odd fundamental martingales with the help of Theorem~\ref{thm:wedgemergp} and Theorem~\ref{thm:barczyp}.  They stem from Theorem~\ref{thm:wedgemercp} and  Theorem~\ref{thm:barczygp}, due to the rotations \eqref{eq:rot1p} and \eqref{eq:rot2p}.
Besides   $\alpha\leftrightarrow\gamma$, the entries $A_t(z)$ and $B_t(z)$ in the upper row of the chain of de Branges matrices are interchanged with the entries $D_t(z)$ and $-C_t(z)$ in the lower row. The role of the reproducing kernel $K_t(w,z)$ is taken over by  kernel \eqref{eq:Qexlap}. In this manner, Theorems~\ref{thm:evenisKLp} -- \ref{thm:barczy1p} are reformulated as follows.    

\begin{thm}\label{thm:evenisKL+p}
Let $X$ be a double-sided si-process  determined by a chain of de Branges spaces of exponential type so as in Theorem~\ref{thm:convwedgep}. Let $M^o$ be the fundamental martingale which defines the moving average representation of the odd part of $X$, with the second order moments
%$\EEE M^e_s\,M^e_t=\pi \alpha_{s\wedge t}$,cf.
\eqref{eq:dubelvar11p}.
Then process $M^o_t/\sqrt{\gamma'_t},$  $t\in[0,r],$   decomposes in  KL-series \eqref{eq:Loevep}   
\[
\frac{1}{\sqrt{\gamma'_t}}M^o_t=\sum_{n=1}^\infty \omega_n \varphi_n(t)\,\xi_n 
\]
where  the $\xi_n$ are i.i.d. Gaussian random variables, $\{\varphi_n(t)= C_t(\lambda_n)/\sqrt{\gamma'_t}\}$ is an orthonormal basis in $L^2[0,r]$, and 
\[
\omega_n=\frac{1}{\lambda_n} \sqrt{\frac{2}{\pi Q_n(\lambda_n,\lambda_n)}}
\]
with the positive zeros $\lambda_n$ of $D_t(z)$.\\
(ii) The martingale bridge 
\begin{equation}\label{eq:oddibrtp}
 B^o_t=M^o_t-\frac{\gamma_t}{\gamma_r}M_r^o,
\end{equation} 
$0\le t\le r$, divided through $\sqrt{\gamma'_t}$,  
decomposes in  KL-series \eqref{eq:Loevep}  
\[
\frac{1}{\sqrt{\gamma'_t}}B^o_t=\sum_{n=1}^\infty \omega_n \varphi_n(t)\,\xi_n 
\]
where  the $\xi_n$ are i.i.d. Gaussian random variables, $\omega_n$ and $\{\varphi_n\}$ are  of the same form as before but now $\lambda_n$ are 
the positive zeros  of $C_t(z)$.

\end{thm}
Recall formula \eqref{eq:Wexla+p} for  kernel $Q(w,z)$ according to which $Q_n(\lambda_n, \lambda_n)=
C_r(\lambda_n) \dot D_r(\lambda_n)$ in item $(i)$ and $Q_n(\lambda_n, \lambda_n)=-\dot C_r(\lambda_n)  D_r(\lambda_n)$ in item $(ii)$.
\\~\\
Let the martingale bridge of the preceding statement be centered by its conditional expectation, given
the value of a martingale $$N^o_t=\int_0^t \kappa_u dM^o_u$$ at the right endpoint $r$. The kernel $\kappa$ is assumed to be non-constant. The Gaussian process so defined
\begin{equation}\label{eq:extmbop}
B^o_t=M^o_t-\EEE (M^o_t| N^o_r)  
\end{equation}
$t\in [0,r]$, has the covariance function
\[
 r(s,t)=\pi \Big(\gamma_{s\wedge t} - \frac{<1_s, \kappa>_\gamma <1_t, \kappa>_\gamma}{\|\kappa\|^2_\gamma}\Big).
 \]
\begin{thm}\label{thm:barczy1+p}
In the situation of Theorem~\ref{thm:barczyp} process
\eqref{eq:extmbop}, divided through $\sqrt{\gamma'_t}$,  
decomposes in  KL-series \eqref{eq:Loevep}  
\[
\frac{1}{\sqrt{\gamma'_t}}B^o_t=\sum_{n=1}^\infty \omega_n \varphi_n(t)\,\xi_n 
\]
where  the $\xi_n$ are i.i.d. Gaussian random variables,  $\{\varphi_n(t)= b_t(w_n)/\sqrt{\gamma'_t}\}$ is an orthonormal basis on $L^2[0,r]$, and 
\[
\omega_n=\frac{1}{w_n c(w_n)}. 
\]
Here $b_t(w)$ is given by \eqref{eq:eigf+p}, the $w_n$
are  positive zeros  of the determinantal equation \eqref{eq:detertal+p} and $c(w)$ is a normalization constant in \eqref{eq:eigf+p}.
\end{thm}
\subsection{Processes associated with an FBM}\label{s:prassfbmp}

(a) In this section we apply the results of the preceding section to obtain KL-expansions for processes associated with a centered Gaussian $H$-self-similar processes $(X^e_t)_{t\ge 0}$  and $(X^o_t)_{t\ge 0}$ of Hurst index $0< H<1$ which are defined in Example~\ref{ex:FBM2p}.
Based on the spectral isometry, it is shown that the martingales  are expressible in terms of   $X^e$ and $X^o$ by \eqref{eq:wclamep}.
The converse relations \eqref{eq:11.2.3eop} are moving average representations.
The second order moment of the martingales are given by \eqref{eq:kvfbmp}.
%Lemma~\ref{lem:MaBessel-} defines 
In this example  the de Branges  matrices   
\begin{equation}\label{eq:MaBesseli11p}
 M_t(z)=
d_t(z)
\begin{bmatrix}
            J_{-H} (tz) &   J_{1-H} (tz)\alpha'(a)\\
           -{J_{H} (tz)}/\alpha'(t) & J_{H-1} (tz)
           \end{bmatrix}
\end{equation}
are defined for $t\ge 0$, where   $z$ is a complex variable and $d_t(z)$ the diagonal matrix
\begin{equation}\label{eq:ddia11p}
 d_t(z)=\text{\rm diag}\big[
            \Gamma(1-H)({tz}/{2})^H,
           \Gamma(H)({tz}/{2})^{1-H}
           \big].
\end{equation}
Check that the determinant is equal to $1$ at hand of the Wronskian  
%\begin{equation}
%\label{eq:Aronskii}
$J_{1-\nu}(z)J_{\nu}(z)+J_{-\nu}(z)J_{\nu-1}(z)=\frac{2}{\pi z}\sin\nu\pi$
%\end{equation}
 \cite[formula 60 on p. 12]{Erd53_2}.

Recall formula \eqref{eq:nulshichp} for  kernel $K_r(\lambda,\lambda)$ on the real diagonal. Formula for kernel $Q_r(\lambda, \lambda)$, cf. \eqref{eq:Wexla+p},  is analogously calculated 
\begin{align}\label{eq:nulshich++p}
 \frac{Q_r(\lambda,\lambda)}{Q_r(0,0)}&=2H\Gamma^2(H) (r\lambda/2)^{2-2H}\nonumber 
 \\
 &\times\Big(J^2_{H}(r\lambda)-\frac{1-2H}{r\lambda}J_{H}(r\lambda)J_{H-1}(r\lambda)+J^2_{H-1}(r\lambda)\Big)
\end{align}
by using %definition \eqref{eq:CDJent-} and 
properties %\eqref{eq:1AppdifB} and \eqref{eq:Hrecnulshi1} 
of the Bessel functions.
\\~\\
\noindent
(b) The general results of the preceding section that regard even and odd fundamental martingales are applied to the present case as follows.
\begin{thm}\label{thm:evenisKL++p}
Let $X$ be a double-sided   FBM of Hurst index $0<H<1$. Then
\\
(i)\; Process $t^{H-\frac12} M^e_t$, $0\le t\le r,$ cf. \eqref{eq:wclamep}, decomposes in  KL-series \eqref{eq:Loevep}  
where  the $\xi_n$ are i.i.d. Gaussian random variables, 
\[
 t^{H-\frac12}\varphi_n(t)=
 \Gamma(1-H)\Big(\frac{t \lambda_n}{2}\Big)^H  
J_{1-H}(t \lambda_n)    
\]
%cf. \eqref{eq:Jent-p}, 
and
\[
 \omega_n=
 \frac{{\sqrt{2}}}
 {\lambda_n\, r^{1-H} \Gamma(1-H) (r\lambda_n/2)^H J_{1-H}(r\lambda_n)}
\]
%cf. \eqref{eq:nulshich+p},
where the $\lambda_n$ are positive roots of the Bessel function $J_{-H}(r\lambda_n)=0$.\\
(ii) \; Process $t^{\frac12-H} M^o_t$, $0\le t\le r,$ cf. \eqref{eq:wclamep}, decomposes in  KL-series \eqref{eq:Loevep}  
where  the $\xi_n$ are i.i.d. Gaussian random variables, 
\[
 \varphi_n(t)=-
 \Gamma(H)\Big(\frac{t \lambda_n}{2}\Big)^{1-H}  
J_{H}(t \lambda_n)     
\]
%cf. \eqref{eq:CDJent-p}, 
and
\[
 \omega_n=
 \frac{{\sqrt{2}}}
 {\lambda_n\,  r^H \Gamma(H) (r\lambda_n/2)^{1-H} J_{H}(r\lambda_n)}
\]
cf. \eqref{eq:nulshich++p},
where the $\lambda_n$ are positive roots of the Bessel function $J_{H-1}(r\lambda_n)=0$.
\end{thm}
In fact, item $(ii)$ is a direct consequence of item $(i)$, for  the martingales in question depend on Hurst index $H$ so that  $ M^e(H)$ $\stackrel{Law}{=}$ $M^o(1-H)$, equality in law.\\~\\ 
Similar results are available for martingale bridges.
\begin{thm}\label{thm:oddisKL++p}
In the situation of Theorem~\ref{thm:evenisKL++p}\\
(i) Process $t^{H-\frac12} B^e_t,$ $0\le t\le r,$ 
cf. \eqref{eq:evenibrtp}, decomposes in  KL-series \eqref{eq:Loevep}  
where  the $\xi_n$ are i.i.d. Gaussian random variables, $\omega_n$ and $t^{H-\frac12}\varphi_n(t)$ is of the same form as in Theorem~\ref{thm:evenisKL++p}, item (i),   
but the
$\lambda_n$ are now positive roots of the Bessel function $J_{1-H}(r\lambda_n)=0$.\\
(ii) Process $t^{\frac12-H} B^o_t,$ $0\le t\le r,$ 
cf. \eqref{eq:oddibrtp}, decomposes in  KL-series \eqref{eq:Loevep}  
where  the $\xi_n$ are i.i.d. Gaussian random variables, $\omega_n$ and $t^{H-\frac12}\varphi_n(t)$ is of the same form as in Theorem~\ref{thm:evenisKL++p}, item (ii),   
but the
$\lambda_n$ are now positive roots of the Bessel function $J_{H}(r\lambda_n)=0$.
\end{thm}
In fact, item $(ii)$ is a direct consequence of item $(i)$, for  the bridges in question depend on Hurst index $H$ so that  $ B^e(H)$ $\stackrel{Law}{=}$ $B^o(1-H)$, equality  in law.

\subsection{Autoregressive processes} \label{s:autorp}
 
In this section we  
discuss the KL decompositions of
{\it autoregressive processes} $Y$ of $n^{\text{th}}$-order  
whose 
spectral measures have densities of the form 
\eqref{eq:muar1p}
where
$\Theta(iz)$ 
is  
a polynomial \eqref{eq:Theta11p} of degree $n$ whose zeros have negative real part.

In the simplest  case of $n=1$ which is  the so-called {\it  Ornstein--Uhlenbeck process},    the covariance function $\EEE(Y_sY_t)=r(s,t) $ is 
of the  exponential form \eqref{eq:covOUp}
with the parameters $\sigma^2, \theta>0$ and the characteristic equation of  the Fredholm type with this kernel
\begin{equation}\label{eq:eqchfp}
 \int_0^r r(s,t)\varphi_l(s)ds=\lambda_l \varphi_l(t)
\end{equation}
is not very difficult to solve. Under the boundary conditions 
\begin{align}\label{eq:mereoup}
 \varphi_l'(0) -\theta \varphi_l(0)=0\nonumber\\ 
 \varphi_l'(r)+\theta \varphi_l(r)=0
\end{align}
one obtains 
\begin{align}\label{eq:klianip}
 \lambda_l &= \frac{\sigma^2}
 {w^2_l+\theta^2}
 \nonumber\\
 \varphi_l(t)&=
 \sqrt{
 \frac{2w_l}
 {r(w_l^2+\theta^2)+2\theta
 }
 }
 \cos (w_l t)+
 \sqrt{
 \frac{2\theta^2}
 {r(w_l^2+\theta^2)+2\theta
 }
 }
 \sin (w_l t)
\end{align}
where the
$w_l$ are positive roots of the determinantal equation
\begin{align}\label{eq:detertaloup}
\begin{vmatrix}
 \theta & -w\\
 \theta-w\tan (rw) & w+\theta\tan (rw)   
    \end{vmatrix}
=0 
\end{align}
arranged in ascending order of magnitude. The proof can be found  in the Ph.D. theses \cite{lim08}, Example 3 on p. 23 (cf. also \cite{cor15}). We shall give more details on  this example  in a separate subsection below.
In what follows we  extend the method of this particular OU case   to  arbitrary  autoregressive processes of order $n$. 

\subsubsection{Mercer's theorem}

We apply   Theorem~\ref{thm:mercerp}  to the autoregressive process $Y$ with  the covariance  function    
\begin{equation}\label{eq:rstThetap}
 r(s,t)=\frac{\sigma^2}{2\pi}\int e^{i\lambda (s-t)} \frac{d\lambda}{|\Theta(i\lambda)|^2}
\end{equation}
with the density of the spectral measure of the form \eqref{eq:muar1p}. 
As before, $\Theta(iz)$
is  
a polynomial \eqref{eq:Theta11p} of degree $n$ and its zeros have negative real part and, moreover, 
 $\phi_k>0$, since process $Y$ is real.

\begin{thm}\label{lem:formnp}
Consider the characteristic equation with kernel \eqref{eq:rstThetap} 
\begin{equation}\label{eq:eqchf+p}
 \int_0^r r(s,t)\varphi_l(s)ds=\lambda_l \varphi_l(t)
\end{equation}
and let the eigenfunctions $\varphi_l(t)$, vanishing outside of interval $[0,r]$,  be $2n$ times continuously differentiable.
 Then eigenfunctions and eigenvalues may be given  the form
 \begin{equation}\label{eq:phinirp}
 \varphi_l(t)=a_l \cos w_l t+ b_l \sin w_l t
\end{equation}
and
\begin{equation}\label{eq:phinirnp}
 \lambda_l=\frac{\sigma^2}{(w_l^2+\phi_1^2)\cdots  (w_l^2+\phi_n^2)}.
\end{equation}
Under the boundary conditions
\begin{align}\label{eq:boundarp}
 \Theta(d/dt) \varphi_l(t)|_{t=0}=0\nonumber\\
 \Theta(-d/dt) \varphi_l(t)|_{t=r}=0
\end{align}
positive numbers $w_l$  are roots of the determinantal equation
\begin{align}\label{eq:detmovap}
 \begin{vmatrix}
 \text{Re\;} \Theta(iw) & 
 \text{Re\;} \Theta(iw)\cos (rw)+ \text{Im\;} \Theta(iw)\sin (rw)
 \\
 \text{Im\;} \Theta(iw)
  &  \text{Re\;} \Theta(iw)\sin (rw)- \text{Im\;} \Theta(iw)\cos (rw)  
    \end{vmatrix}=0
\end{align}
and  the coefficients $a_l$ and $b_l$ are roots of the system of equations
\begin{align}\label{eq:alblsp}
 a \,\text{Re\;} \Theta(iw) +b\, \text{Im\;} \Theta(iw) &=0\nonumber\\
 \int_0^r \big(a\cos (wt)+b\sin (wt) \big)^2\,dt &=1. 
\end{align}

\end{thm}
Note that the eigenvalues are $2\pi$ multiple of the spectral density at points $w_l$.

\begin{prf}
By applying $|\Theta(d/dt)|^2$ to both sides of the characteristic equation \eqref{eq:eqchf+p}, one obtains
\begin{equation}\label{eq:phinisp}
 \sigma^2 \varphi_l(t)=\lambda_l |\Theta(d/dt)|^2 \varphi_l(t).
\end{equation} 
Indeed, 
since 
\[
 |\Theta(d/dt)|^2 e^{-i\lambda t}=|\Theta(i\lambda)|^2e^{-i \lambda t},
\]
the spectral representation \eqref{eq:rstThetap}  turns the left-hand side of the characteristic equation into 
$$
%\label{eq:hatia}
 \frac{\sigma^2}{2\pi} 
 \int_0^r e^{-i\lambda t}\hat\varphi_l(i\lambda) d\lambda
$$
where $\hat\varphi_l(i\lambda)$ is the Fourier transform of $\varphi_l$, 
i.e. 
$\hat\varphi_l(i\lambda)=\int_0^r e^{it\lambda} \varphi_l(t)dt,$
and the required expression is obtained by taking the inverse Fourier transform.

If the eigenfunctions are given the form \eqref{eq:phinirp}, then
we have
\[
 |\Theta(d/dt)|^2 \sin (\lambda t)=\prod_{j=1}^n \Big(-\frac{d^2}{dt^2}+\phi_j^2\Big) \sin (\lambda t) =\prod_{j=1}^n \big(\lambda^2+\phi_j^2\big) \sin (\lambda t)
\]
and the same for cosine, therefore equation \eqref{eq:phinisp} becomes 
\begin{equation*}
%\label{eq:phinis}
 \sigma^2 \varphi_l(t)=\lambda_l |\Theta(iw_l)|^2 \varphi_l(t)
\end{equation*}
which determines  the eigenvalues \eqref{eq:phinirnp}. 

To establish the first of boundary conditions \eqref{eq:boundarp}
apply $\Theta(d/dt)$ to both sides of the spectral representation using $\Theta(d/dt) e^{iz t}=\Theta(iz) e^{izt}$.
Since 
\[
 \Theta(d/dt) r(s,t)=\frac{\sigma^2}{2\pi}\int e^{i\lambda (s-t)} \frac{d\lambda}{\Theta(i\lambda)}
\]
vanishes at $t=0$ by the residue theorem we obtain the first of boundary conditions \eqref{eq:boundarp}.
Similarly, for $0\le s\le r$
\[
 \Theta(-d/dt) r(s,t)=\frac{\sigma^2}{2\pi}\int e^{i\lambda (s-t)} \frac{d\lambda}{\Theta(-i\lambda)}
\]
vanishes at $t=r$ and we obtain the second of boundary conditions \eqref{eq:boundarp}.

With the eigenfunctions of the form \eqref{eq:phinirp}, the boundary conditions give the following system of equations 
\begin{align*}
%\label{eq:boundar+}
& a \,\text{Re\;} \Theta(iw) {+}b\, \text{Im\;} \Theta(iw){=}0\nonumber\\
& a \big(\,\text{Re\;} \Theta(iw) \cos (rw){+} \,\text{Im\;}\sin(rw) \big) {+}b \big(\,\text{Re\;} \Theta(iw) \sin (rw) {-} \,\text{Im\;}\cos(rw) \big){=}0
\end{align*}
whose non-zero solution for $a$ and $b$ requires a vanishing determinant which is \eqref{eq:detmovap}. So, the $w_n$ are roots of this determinantal equation, for their characterisation see below Lemma \ref{lem:rootsp}. 

To calculate the coefficients $a_l$ and $b_l$ make use of the first of boundary conditions and the requirement that the eigenfunctions have to be of unit norm in the metric of $L^2[0,r]$.
\end{prf}
To compute the integral in \eqref{eq:alblsp}, note
\[
 2\varphi^2(t)=\frac{d}{dt} \big((a^2+b^2) t -w^{-2}\varphi(t)\varphi'(t)\big) 
\]
since $\varphi''(t)=-w^2 \varphi(t)$ and $\varphi^2(t)=a^2+b^2 -\big(\varphi'(t)/w\big)^2$. Therefore
%\begin{align*}
 %\frac12\frac{d}{dx}\Big( (a^2+b^2) x +(a^2-b^2) \frac{\sin (wx)\cos (wx)}{w} +ab \frac{\sin^2(wx)-\cos^2(wx)}{w}\Big\\=(a\cos(xw)+b\sin(xw))^2
%\end{align*}
%which gives 
\begin{align}\label{eq:gradop}
 2\int_0^r\varphi(t)^2dt&=(a^2+b^2) r-w^{-2}\varphi(t)\varphi'(t)\big|_0^r\nonumber\\ 
&= (a^2+b^2) r +(a^2-b^2) \frac{\sin (wr)\cos (wr)}{w}+2ab \frac{\sin^2(wr)}{w}.
\end{align}
%use \cite[Table 2.532]{Gre84} which gives an indefinite integral
For more details see concluding subsection where  
Theorem~\ref{thm:arcontip} is presented and proved.

For illustration, we apply the general Theorem~\ref{lem:formnp} to two simplest cases. The case $n=1$ is discussed in the next subsection and $n=2$ in Section~\ref{s:n&a}, Note 6. 

\subsubsection{OU process}
With kernel \eqref{eq:covOUp} the characteristic equation \eqref{eq:eqchfp} becomes
\[
 \frac{\sigma^2}{2\theta} \Big(\int_0^t e^{-\theta (t-s)}\varphi_l(s)ds+\int_t^re^{\theta (t-s)}\varphi_l(s)ds\Big) =\lambda_l \varphi_l(t).
\]
Differentiating both sides yields
\[
 \frac{\sigma^2}{2} \Big(-\int_0^t e^{-\theta (t-s)}\varphi_l(s)ds+\int_t^re^{\theta (t-s)}\varphi_l(s)ds\Big) =\lambda_l \varphi'_l(t).
\]
Evaluate first at $t=0$ and then at $t=r$ to get the boundary conditions \eqref{eq:mereoup}.
Differentiating once more yields
\begin{align*}
 \lambda_l \varphi_l''(t)= \frac{\sigma^2}{2}\Big(\theta \int_0^r e^{-\theta|t-s|}\varphi_l(s)ds-2\varphi_l(t)\Big)=\varphi_l(t)(\theta^2\lambda_l-\sigma^2),
\end{align*}
hence
\[
 \varphi_l''(t)= -w_l^2 \varphi_l(t) 
\]
where $w_l^2=(\sigma^2-\theta^2\lambda_l)/\lambda_l$.
Solving the latter with respect to $\lambda_l$ gives the first equation  \eqref{eq:klianip}. 

The  preceding second order differential equation is solved by a linear combination of the sine and cosine  
$\varphi (t)=a\cos (wt)+b\sin (wt)$. To determine the eigenvalues  $w_l$ of this differential equation use the boundary conditions \eqref{eq:mereoup} which give the system of equations
\begin{align*}
&a\theta-bw=0\\
 & a(\theta-w \tan(rw) )+b(\theta\tan(rw)+w)=0 
\end{align*}
whose non-zero solution is guaranteed by the determinantal equation \eqref{eq:detertaloup}. As for
 the coefficients $a$ and $b$, which are related by $b=(\theta/w) a$, we can use the normalization requirement $\int_0^r\varphi^2(t)dt=1$. Since 
\begin{equation}\label{eq:tangip}
  \tan(rw)=\frac{2\theta w}{w^2-\theta^2}
\end{equation}
 by \eqref{eq:detertaloup}, formula \eqref{eq:gradop} gives 
 \[
  a=\sqrt{ \frac{2w^2}{r(w^2+\theta^2)+2\theta}} \qquad b= \sqrt{\frac{2\theta^2}{r(w^2+\theta^2)+2\theta}}.
 \]
which in tern gives expression \eqref{eq:klianip} for the eigenfunctions.
\\~\\
\noindent
We need to  bring  OU processes  in connection with the chain of de Branges spaces $\HH(E_t)$, $t\ge 0$,  contained isometrically in $L^2(\mu)$ with the spectral measure of the processes \eqref{eq:10ship}.
%\begin{equation}\label{eq:muarnp}
% \mu(d\lambda)= \frac{\sigma^2}{2\pi}\,\frac{d\lambda}{ \theta^2+\lambda^2}.
%\end{equation}
As is said above,  the entire functions \eqref{eq:parelispp} 
%\begin{equation}\label{eq:parelisp} 
%E_t(z)=\sqrt{2\pi/\sigma^2}(\theta-iz)e^{-izt}\qquad t\ge 0
%\end{equation}
are de Branges functions  and all spaces they generate are contained isometrically in $L^2(\mu)$.  
Let us normalize the de Branges function \eqref{eq:parelispp} at the origin. The quotient  
\begin{equation}\label{eq:1eup}
  E^0_t(z)=E_t(z)/E_t(0)=(1-iz/\theta) e^{-itz}
\end{equation}
has  the components 
 $$A_t(z)= \cos (tz) -(z/\theta) \sin (tz) \qquad B_t(z)=\sin (tz) +(z/\theta)\cos (tz).
$$
%and $E(z)E^\sharp(z)=A^2(z)+B^2(z)= 1+z^2/\phi^2$.
Section~\ref{s:stPWep} makes use of sampling points $\lambda_n$ which are positive zeros of $B_r(z)$ and satisfy  equation $\tan(\lambda_n r)=\lambda_n/\theta$.  The present sampling points that come from the identity \eqref{eq:tangip} are  positive zeros of an odd component of another  de Branges function which is $(1-iz/\theta)$ multiple of  $E_r(z) $, i.e. 
\begin{equation}\label{eq:1eu1p}
 E^1_r(z)\pd (1-iz/\theta)^2 e^{-irz}.
\end{equation}
 Its even and odd components are 
\begin{align*}
 A_r^1(z)&= (1-(z/\theta)^2)\cos (rz) -2(z/\theta) \sin (rz)\\
 B_r^1(z)&=(1-(z/\theta)^2)\sin (rz) +2(z/\theta)\cos (rz).
\end{align*}
%Note
%\begin{equation}\label{eq:berti}
% A_r^1(z)=A_r(z)-\frac{z}{\theta} B_r(z)\qquad B_r^1(z)=B_r(z)+%\frac{z}{\theta} A_r(z).
%\end{equation}
%The kernel $K^1(w,z)$ in the space $\HH(E^1)$ is 
%\[
% K^1(w,z)=\big(1+\frac{z\bar w}{\phi^2}\big)K(w,z)+ \frac{1}{\pi\phi} \big(A(z)A(\bar w)+B(z)B(\bar w)\big)
%\]
%with
%\[
% K(w,z)=\frac{(1+z\bar w/\phi^2)\sin (z-\bar w) a+ ((z-\bar w)/\phi)\cos (z-\bar w) a}{\pi(z-\bar w)}.
%\]
Clearly, the positive zeros of $B_r^1(z)$ satisfy
\begin{equation*}
%\label{eq:roots}
 (1-(w/\theta)^2)\sin (rw)+2(w/\theta)\cos (rw)=0
\end{equation*}
which matches \eqref{eq:tangip}. 
It is seen from \eqref{eq:klianip} that at  points $w_l$ which are roots of the determinantal equation \eqref{eq:detertaloup} and in the same time  positive zeros of $B_r^1(z)$, the 
functions $B_t(w_l)=\sin tw_l +(w_l/\phi)\cos tw_l $ become constant multiples of the eigenfunctions, i.e.
\[
\varphi_l(t)=\sqrt{\frac{2\theta^2}{r(w_l^2+\theta^2)+2\theta}}\;
B_t(w_l).
\] 
In this manner 
the 
set of functions $\{B_t(w_1), B_t(w_2),\cdots\}$  form a basis in $L^2[0,r]$. This fact is confirmed  
by the Lagrange formula \eqref{eq:orilagp} for two different roots $w_j$ and $w_k$ 
\begin{align}\label{eq:lagrasn1p}
\int_0^r B_u( w_l) B_u(w_k)\,du= \frac{ w_jB_t(w_k) A_t(w_j)-w_kB_t(w_j)A_t(w_k)}{w_k^2-w_j^2}\Big|_0^r
\end{align}
which is equal to $0$ at point $t=0$, since $A_0(z)=1$ and $B_0(z)=z/\theta$, and at point $t=r$, since 
\begin{equation}\label{eq:idenr0p}
A^2_r(w)=A^2_0(w) \qquad\quad B^2_r(w)=B^2_0(w)
\end{equation}
if $B_r^1(w)=0$. These two identities are obtained by simple algebra. For instance, 
\begin{align*}
 \theta^2B_r^2(w)&=(\theta\sin(rw)+w\cos(rw) )^2\\&= (\theta^2-w^2)\sin^2(tw)+2\theta w \sin(rw)\cos(rw) +w^2=\theta^2B^2_0(w).
\end{align*}

\subsubsection{AR processes {\rm (continued)}}

The relationship of an OU process to a chain of de Branges spaces just discussed does extend to autoregressive processes. 
To show this, we
consider the de Branges function \eqref{eq:parelisap} for $t\ge 0$ which generates a space contained isometrically in $L^2(\mu)$ with the spectral measure   \eqref{eq:muar1p}. Like  in \eqref{eq:1eup}, normalize this de Branges function at the origin   
\[
 E^0_r(z)=E_t(z)/E_t(0)= \Theta(iz)\, e^{-izr}.
\]
This normalized de Branges function has even and odd components
\begin{align}\label{eq:thetaabp}
 A_r(z)&= \text{Re\;} \Theta(iz)\cos (rz) +\text{Im\;} \Theta(iz)\sin (rz)\nonumber\\
 B_r(z)&= \text{Re\;} \Theta(iz)\sin (rz) -\text{Im\;} \Theta(iz) \cos (rz).
\end{align}
Like in \eqref{eq:1eu1p} we consider another de Branges function
\[
 E^1_t(z)= \Theta^2(iz)\, e^{-izt}
\]
with the even and odd components
\begin{align}\label{eq:thetaab1p}
 A_t^1(z)&= \big((\text{Re\;} \Theta(iz))^2-(\text{Im\;} \Theta(iz))^2\big)\cos (tz) +2\text{Re\;} \Theta(iz)\text{Im\;} \Theta(iz) \sin (tz)\nonumber\\
 B_t^1(z)&= \big((\text{Re\;} \Theta(iz))^2-(\text{Im\;} \Theta(iz))^2\big)\sin (tz )-2\text{Re\;} \Theta(iz)\text{Im\;} \Theta(iz) \cos (tz).
\end{align}
The latter, being an odd component of a de Branges function, has a simple interlacing zeros symmetrically  spread about the origin. Clearly, they satisfy the determinantal equation \eqref{eq:detmovap}. Moreover, if $w_l$ is the $l^{\text{th}}$ positive root of the latter, then the odd component \eqref{eq:thetaabp} evaluated at this point $B_t(w_l)$ is a constant multiple of the $l^{\text{th}}$ eigenfunction \eqref{eq:phinirp}, see  Theorem~\ref{thm:arcontip} below. The proof of this theorem is based on the identities \eqref{eq:idenr0p} which hold true in general.
\begin{lem}\label{lem:rootsp}
Let $w$ be a zero of the odd component $B^1_r(z)$ of the form  \eqref{eq:thetaab1p}, i.e.  $B^1_r(w)=0$. Then the components  \eqref{eq:thetaabp} evaluated at  this point  satisfy the identities \eqref{eq:idenr0p}.

If $w_j$ and $w_k$ are two different zero of $B_r^1(z)$, then 
\begin{align}
\label{eq:ortbwp}
 \int_0^r B_t(w_j)B_t(w_k)dt=0.
\end{align}

\end{lem}
\begin{prf}
The proof of \eqref{eq:idenr0p} is as simple as in the case $n=1$. Moreover, the statement stay true even if polynomial $\Theta(iz)$ is substituted by any de Branges function, say $e(z)=a(z)-ib(z)$. Indeed, the even and odd components $a(z)$ and $b(z)$ replace $\text{Re\;} \Theta(iz)$ and $-\text{Im\;} \Theta(iz)$, respectively, and if 
\[
 E^1_t(z)= e^2(z)\, e^{-izt}
\]
then its  even and odd components are
\begin{align*}
%\label{eq:thetaab1}
 A_t^1(z)&= \big(a^2(z)-b^2(z)\big)\cos (tz) -2a(z)b(z) \sin (tz)\nonumber\\
 B_t^1(z)&= \big(a^2(z)-b^2(z)\big)\sin (tz )+2a(z)b(z) \cos (tz).
\end{align*}
Now, if $w$ is a zero of the latter component at the  right endpoint $r$, i.e. $B_r^1(w)=0$, then
\begin{align*}
 A_r^2(w)&{=}a^2(w)\cos^2 (rw) + b^2(w)\sin^2(rw)
 -2a(w)\;b(w)\sin (rw)\cos (rw)\\&
 {=}a(w)^2\cos^2 (rw) + b^2(w)\sin^2(rw)+
 \big(a^2(w)-b^2(w)\big)\sin^2(rw)=A_0^2(w)
\end{align*}
and
\begin{align*}
B_r^2(w)&{=}a^2(w)\sin^2 (rw) + b^2(w)\cos^2(rw)+2a(w)\;b(w)\sin (rw) \cos (rw)\\&
 {=}a^2(w)\sin^2 (rw) + b^2(w)\cos^2(rw)-
 \big(a^2(w)-b^2(w)\big)\sin^2(rw)=B_0^2(w).
\end{align*}
The proof of \eqref{eq:ortbwp} is based on the Lagrange identity \eqref{eq:lagrasn1p}, precisely as in the previous case of $n=1$. 
\end{prf}

The preceding results can be summarized as follows.

\begin{thm}\label{thm:arcontip}
 The eigenfunctions and eigenvalues are of the form
\begin{align*}
\varphi_l(t)=\frac{
B_t(w_l)}{\|B(w_l)\|}\qquad\quad
\lambda_l=2\pi \frac{\mu(d\lambda)}{d\lambda}\Big|_{\lambda=w_l}
\end{align*}
where $w_l$ are positive zeros of $B_r^1(z)$, cf. \eqref{eq:thetaab1p}, and
\[
 2\|B(w_l)\|^2=\big((\text{\rm Re\;} \Theta(iw_l))^2 + (\text{\rm Im\;} \Theta(iw_l))^2\big)r+2\text{\rm Re\;} \Theta(iw_l)\;\text{\rm Im\;} \Theta(iw_l)                                         .
 \]

\end{thm}
The norm is calculated by \eqref{eq:alblsp} with $a=\text{Re\;} \Theta(iw_l)$ and $b=\text{Im\;} \Theta(iw_l)$.

\section{Notes and addenda}\label{s:n&a}
\begin{small}

%\subsection*{Section~\ref{s:PWs}}

\begin{enumerate}

\item{
Regarding the rate of convergence of Paley--Wiener series, discussed in Section~\ref{sub:pwserp}, 
see \cite[Theorem 4.1.4]{zar07}. In Chapter 5 of the latter thesis the results of the theorem are applies to estimate small ball probabilities for si- and stationary processes. Much more far-reaching  developments of theory of small ball asymptotics  can be found in the recent survey   \cite{naz23}.
}

\begin{comment}
 \begin{align*}
 \EEE \|X-X^N\|_p
 \symeq 
 N^{
 -\frac{1+a}{2}
 }
 \sqrt{\ln N}
\end{align*}
for $p\ge 1$ and
 \begin{align*}
 \EEE \|X-X^N\|_{\psi_p}
 \symeq 
 N^{
 -\frac{1+a}{2}
 }
 \Big(
 \sqrt{\ln N}\
 Big)^{1-2/p}
\end{align*}
for $p\ge 2$.
\end{comment}

%\end{enumerate}
%\subsection*{Section~\ref{s:KL}}
%\begin{enumerate}

\item{ 
In
\cite{bar18} Theorem~\ref{thm:barczy1p} is proved in the particular case of the martingale bridge \eqref{eq:extmbp} where $M^e$ is a standard Brownian motion, so that kernel \eqref{eq:barczyp} is simply $k(s,t)=s\wedge t- \int_0^s \kappa_u du  \int_0^t \kappa_u du$ (assuming $\kappa$ is of unit norm, for simplicity). 

}
\item{
The bridging property $B_0^e=B_t^e=0$ of an even  martingale bridge \eqref{eq:evenibrtp}  suggests notion of an {\it inverse bridge} $\overleftarrow {B}_t^e=B_{r-t}^e$ which is a cantered Gaussian process with the covariance function 
\[
 \EEE(\overleftarrow {B}^e_s\overleftarrow {B}_t^e)=\alpha_{r- s \vee t}-\frac{\alpha_{r-s}\alpha_{r-t}}{\alpha_{r}}.
\]
In these terms the following  inverse version of Theorem~\ref{thm:evenibrKLp} holds true.

\begin{thm}\label{thm:evenibrrp}
In the situation of Theorem~\ref{thm:evenisKLp}, the inverse martingale bridge $\overleftarrow {B}_t^e=B_{r-t}^e$, 
$0\le t\le r$, divided through $\sqrt{\alpha'_{r-t}}$,  
decomposes in  KL-series \eqref{eq:Loevep}  
\[
\frac{1}{\sqrt{\alpha'_{r-t}}}\overleftarrow {B}^e_t=\sum_{n=1}^\infty \omega_n \varphi_n(t)\,\xi_n 
\]
where  the $\xi_n$ are i.i.d. Gaussian random variables, 
$\{\varphi_n(t)= B_{r-t}(\lambda_n)/\sqrt{\alpha'_{r-t}}\}$  and 
\[
\omega_n=\frac{1}{\lambda_n} \sqrt{\frac{2}{\pi K_r(\lambda_n,\lambda_n)}}
\]
with the positive zeros $\lambda_n$ of $A_r(z)$.

\end{thm}
}

\item{
The main result of the paper \cite{bar11} applies Theorem~\ref{thm:evenibrrp} to a fractional Brownian motion, more precisely, to an inverse version of an associated even martingale bridge. This is not immediately seen, since a Gaussian process that is studied in 
\cite{bar11}, is defined on interval $[0,r]$  as a Wiener integral with respect to a standard Brownian motion
\[
 X_t^{(\alpha)}=\int_0^t\Big(\frac{r-t}{r-s}\Big)^\alpha\, dW_s
\]
and called {\it alpha-Wiener bridge}.  In the present text symbol $\alpha$ is  used for other purposes and to avoid ambiguities in this note the positive exponent $\alpha$ is substituted by $\frac32-H$. We will discuss only the case $0<H<1$ and  shown that process $X^{(\frac32-H)}$ is associated in a  certain way with a fractional Brownian bridge of Hurst index $H$. More precisely,
a relation will be sought to an even martingale \eqref{eq:wclamep}
whose quadratic variation is given by \eqref{eq:kvfbmp}, with $\alpha_t=t^{2-2H}/(2-2H)$. It is easily seen that
the normalized process
\[
 X_t\pd (r-t)^{\frac12-H}X_t^{(\frac32-H)}
\]
is an {\it anticipative martingale bridge} for a martingale
\[
 M_t^r\pd \int_0^t (r-u)^{\frac12-H} dW_u
\]
where $W$ is a standard Brownian motion (the term is taken over from  \cite{gas07}). By definition, process 
$M_t^r$ on interval $ [0,r]$ defines  the anticipative martingale bridge $M_1(t,r)$ by
\[
 M_1(t,r)=\int_0^t \frac{<M^r>_r-<M^r>_t}{<M^r>_r-<M^r>_t}dM^r_u
\]
where $<M^r>$ is the quadratic variation of the martingale $M^r$ and equals to $ <M^r>_t=\alpha_r-\alpha_{r-t}$, since
\[
 \EEE(M_s^rM_t^r)=\alpha_r-\alpha_{r-s\wedge t}.
\]
Simple algebra shows identity $X_t=M_1(t,r)$ and the following form of the covariance function 
\[
 \EEE(X_sX_t)=\alpha_{r- s \vee t}-\frac{\alpha_{r-s}\alpha_{r-t}}{\alpha_{r}}.
\]
Compare this with the covariance function $\EEE(\overleftarrow {B}^e_s\overleftarrow {B}_t^e)$ of an inverse  even martingale bridge provided in the preceding note. The covariances coincide.  Hence, the present normalized alpha-Wiener bridge $\{X\}_{t\in [0,r]}$ is equal in law to the inverse  bridge $\{\overleftarrow {B}_t^e=B_{r-t}^e\}_{t\in [0,r]}$
where $B^e$ is as in Theorem~\ref{thm:oddisKL++p}. Theorem~\ref{thm:evenibrrp} can be applied with 
$\alpha_t=t^{2-2H}/(2-2H)$ and the components $A_t(z)$ and $B_t(z)$ specified in terms of the Bessel functions as in  Section~\ref{s:prassfbmp}.

}
\item{
Theorem~\ref{thm:evenisKL++p} is proved in \cite{deh03} (cf. \cite{deh06}, Fact 1.2).  
}
\item{

A second order autoregressive process  has  spectral density  of the form
\[
 \mu(d\lambda)= \frac{\sigma^2{d\lambda}}{2\pi(z^2+\phi_1^2)(z^2+\phi_2^2)}
\]
with parameters $\sigma^2>0$, $\phi_1>0$, $\phi_2>0$. According to Lemma~\ref{lem:formnp}, the covariance function $r(s,t)=\EEE(X_sX_t)$ on interval $[0,r]$ has 
the Karhunen--Lo\`{e}ve  expansion 
\[
 r(s,t)=\sum_{n=1}^\infty \lambda_n\varphi_n(s)\varphi_n(t)
\]
with
the eigenfunctions and eigenvalues which are of the form
 \[
 \varphi_n(t)=a_n \cos w_n t+ b_n \sin w_n t
\]
and
\[
 \lambda_n=\frac{\sigma^2}{(w_n^2+\phi_1^2)(w_n^2+\phi_2^2)}
\]
for certain numbers $w_n$ and certain coefficients $a_n$ and $b_n$, both depending on $w_n$. To determine the numbers $w_n$,  consider the boundary conditions
\[
 (d/dt -\phi_1)(d/dt-\phi_2)\varphi_n(0)=0\qquad (d/dt +\phi_1)(d/dt+\phi_2)\varphi_n(1)=0
\]
which for the eigenfunctions of the preceding form imply the system of equations
\begin{align*}
 a_n(\phi_1\phi_2-w^2)&-b_n (\phi_1+\phi_2)w=0\\
 a_n\big((\phi_1\phi_2-w^2)\cos w r&-(\phi_1+\phi_2)w\sin wr\big)\\&+ b_n\big((\phi_1\phi_2-w^2)\sin wr+(\phi_1+\phi_2)w\cos wr\big)=0.
\end{align*}
A non-zero solution with respect to $a_n$ and $b_n$ requires a vanishing determinant
\begin{align*}
&\begin{vmatrix}
 (\phi_1\phi_2-w^2) & (\phi_1\phi_2-w^2)\cos wr-(\phi_1+\phi_2)w\sin wr\\
 -(\phi_1+\phi_2)w & (\phi_1\phi_2-w^2)\sin wr +(\phi_1+\phi_2)w\cos wr
\end{vmatrix}
\\&=
 \big((\phi_1\phi_2-w^2)^2-(\phi_1+\phi_2)^2w^2\big)\sin w+2w(\phi_1+\phi_2)(\phi_1\phi_2-w^2)\cos w=0.
\end{align*}
Let the numbers $w_n$ be positive roots of this determinantal equation (it will be seen in a moment that these are positive zeros of the entire function $B_r^1(z)$ given by \eqref{eq:ertebi+p}).
To determine the numbers $a_n$ and $b_n$, make use of the first of boundary equations
\[
 a_n(\phi_1\phi_2-w_n^2)=b_n (\phi_1+\phi_2)w_n
\]
and the fact that the eigenfunctions are normalized
\[
 \int_0^r \varphi^2_n(t)dt =\int_0^r (a_n \cos w_nt+b_n\sin w_nt)^2dt.
\]
%(see the calculations at the end of the present note).\\~\\
%\noindent
%{\bf The frequency approach.}
To shed more light on these results we consider 
space $L^2(\mu)$ with the present spectral measure which contains isometrically the  part of the chain of de Branges spaces restricted to interval $[0,r]$ and is generated by the de Branges functions
\[
 E_t(z)=(iz-\phi_1)(iz-\phi_2) e^{-izt},
\]
The corresponding odd and even  components are 
\begin{align}
\label{eq:titon+p}
A_t(z)&=(\phi_1\phi_2-z^2)\cos zt-(\phi_1+\phi_2)z\sin zt
 \nonumber\\
 B_t(z)&=(\phi_1\phi_2-z^2)\sin zt+(\phi_1+\phi_2)z\cos zt .
\end{align}
Clearly, $A_0(z)=\phi_1\phi_2-z^2$ and $B_0(z)=(\phi_1+\phi_2)z$. We shell show that at certain points $w$ on the real axis   
$A^2_r(w)=A^2_0(w)$ and $B^2_r(w)=B^2_0(w)$. 
To this end in view, consider
the following de Branges function
\[
 E^1_t(z)=(iz-\phi_1)^2(iz-\phi_2)^2 e^{-izt}
\]
with the even and odd components
\begin{align}\label{eq:ertebi+p}
A^1_t(z)&=\big(A^2_0(z)-B^2_0(z)\big)\cos zt-2A_0(z)B_0(z)\sin zt
\nonumber\\
 B^1_t(z)&=\big(A^2_0(z)-B^2_0(z)\big)\sin zt+2A_0(z)B_0(z)\cos zt 
\end{align}
%(in fact, 
%\[ E^1_t(z)=\Big(\big((\phi_1\phi_2-z^2)^2-(\phi_1+\phi_2)^2z^2\big)-2zi(\phi_1+\phi_2)(\phi_1\phi_2-z^2)\Big)e^{-izt},\] but since $(iz-\phi_1)(iz-\phi_2)=\phi_1\phi_2-z^2-iz(\phi_1+\phi_2)$ we have
%\[E^1_t(z)=(iz-\phi_1)^2(iz-\phi_2)^2 e^{-izt}\] as above). 
As is said above,
the even and odd components $A(z)$ and $B(z)$ of a de Branges function $E(z)=A(z)-iB(z)$ have only simple interlacing real zeros. The following statement makes use of this fact.\\
{\it 
%\begin{lem}
%\label{lem:roots+}
At the right endpoint $t=r$,  let $w_1, w_2, \cdots$ be positive zeros of $B^1_1(z)$  of ascendant degree, i.e. the roots of the equation
\begin{align*}
 \big((\phi_1\phi_2-w^2)^2-(\phi_1+\phi_2)^2w^2\big)\sin wr+2(\phi_1+\phi_2)(\phi_1\phi_2-w^2)w\cos wr=0.
\end{align*}
For any such number $w$
\begin{align*}
%\label{eq:nukii+}
 A_r^2(w)=(\phi_1\phi_2-w^2)^2=A_0^2(w)\qquad B_r^2(w)=(\phi_1+\phi_2)^2w^2=B_0^2(w).
\end{align*}
}
%\end{lem}
%\begin{prf}
The right-hand equations are obvious by definition \eqref{eq:titon+p}. Calculations for the left-hand equations are straightforward
\begin{align*}
 B_r^2(w)&=(\phi_1\phi_2-w^2)^2\sin^2 wr + (\phi_1+\phi_2)^2w^2\cos^2wr\\
 &+2(\phi_1\phi_2-w^2)(\phi_1+\phi_2)w\sin wr\cos wr\\&
 =(\phi_1\phi_2-w^2)^2\sin^2 wr + (\phi_1+\phi_2)^2w^2\cos^2wr\\&-
 \big((\phi_1\phi_2-w^2)^2-(\phi_1+\phi_2)^2w^2\big)\sin^2wr=B_0^2(w)
\end{align*}
and
\begin{align*}
 A_r^2(w)&=(\phi_1\phi_2-w^2)^2\cos^2 wr+ (\phi_1+\phi_2)^2w^2\sin^2wr\\
 &-2(\phi_1\phi_2-w^2)(\phi_1+\phi_2)w\sin wr \cos wr\\&
 =(\phi_1\phi_2-w^2)^2\cos^2 wr + (\phi_1+\phi_2)^2w^2\sin^2wr\\&+
 \big((\phi_1\phi_2-w^2)^2-(\phi_1+\phi_2)^2w^2\big)\sin^2wr=A_0^2(w).
\end{align*}
%\end{prf}
This statement implies the orthogonality of the sequence of functions $B_t(w_n)$. More precisely, we can claim the following\\
{\it 
 Let $w_j$ and $w_k$ be zeros of $B_1^1(z)$. Then 
\begin{align*}
 \int_0^r B_t(w_j)B_t(w_k)dt=0.
\end{align*}
}
This is simply verified. Indeed,
if $f''=-w^2f$ and $g''=-z^2g$, then 
\[
 \int fg dt=\frac{gf'-g'f}{z^2-w^2}
\]
since the derivative of the right-hand side is equal to $(gf''-g''f)/(z^2-w^2)$ and this in turn is equal to $fg$. Hence

\begin{align*}
 \int_0^r B_t(w_j)B_t(w_k)dt=\frac{w_j A_t(w_j) B_t(w_k)-w_kA_t(w_k)B_t(w_j)}{w^2_j-w^2_k}\big|_0^r
\end{align*}
and by preceding statement this is equal to zero.\\

In summary we have 
\begin{thm}\label{thm:abdaphi+p}
 The eigenfunctions and eigenvalues are of the form
\begin{align*}
\phi_n(t)=\frac{B_t(w_n)}{\|B(w_n)\|}=\sqrt{\frac{2}{(a_n+b_n)^2r+2a_nb_n}}\,\big(a_n\cos w_nt+b_n\sin w_nt\big)
 \end{align*}
%\phi_n(t)&=B_t(w_n)/\|B(w_n)\|\\&=\sqrt{\frac{2}{(w_n^2+\phi_1^2)(w_n^2+\phi_2^2)}}\,\big((\phi_1\phi_2-w_n^2)\sin w_nt+(\phi_1+\phi_2)w_n\cos w_nt\big)\end{align*}
and
\[
 \lambda_n=\frac{\sigma^2}{a_n^2+b_n^2}
\]
where $w_n$ are as in Lemma~\ref{lem:rootsp} and $a_n=(\phi_1+\phi_2)w_n$ and $b_n=\phi_1\phi_2-w_n^2$.
\end{thm}
The required norm is calculated as in \eqref{eq:gradop}.

}
\end{enumerate}

\end{small}

\bibliographystyle{plain}
%\bibliography{deBrangesp}
%\backmatter
%\printindex
\input{paper.bbl}

\end{document}

%% file: paper.bbl
\def\cprime{$'$}